\documentclass[arma]{svjour}
\usepackage{graphics}
\usepackage{latexsym}
\usepackage{amssymb,mathrsfs,enumerate}
\usepackage{amsbsy}
\usepackage{amsmath}
\usepackage{epsfig}
\usepackage{graphpap,latexsym,epsf}
\newcommand{\R}{\Bbb{R}}

\newcommand{\N}{\Bbb{N}}

\newcommand{\M}{\Bbb{M}}

\newcommand{\HH}{\mathscr{H}}

\newcommand{\nn}{{\mbox{\boldmath$n$}}}

\renewcommand{\div}{\mathop{\rm div}\nolimits}
\newcommand{\weaksto}{{\rightharpoonup^*}}
\newcommand{\weakto}{\rightharpoonup}
\newcommand{\restr}[1]{\lower3pt\mbox{$|_{#1}$}}
\newcommand{\Restr}[1]{\lower3pt\mbox{$\Big|_{#1}$}}

\newcommand{\conv}{\mathop{\rm co}}

\newcommand{\trepar}{{|\kern-1truept|\kern-1truept|}}
\newcommand{\ltt}{{(\kern-2truept(}}
\newcommand{\rtt}{{)\kern-2truept)}}
\newcommand{\duality}[2]{\left\langle\vphantom{\big(}%
                        {#1},{#2}\right\rangle}
\newcommand{\la}{\langle}
\newcommand{\ra}{\rangle}

\newenvironment{proofx}{\removelastskip\par\medskip  
\noindent{\em Proof.} \rm}{\par}
\newcommand{\nchi}{{\raise.4ex\hbox{$\chi$}}}

\newcommand{\teta}{\vartheta}
\newcommand{\down}{\downarrow}
\newcommand{\up}{\uparrow}
\newcommand{\forae}{\text{for a.e.}}
\newcommand{\aein}{\text{a.e.\ in}}

\newcommand{\Hilbert}{\HH}

\newcommand{\convexSubdifferential}{\partial}

\newcommand{\cvSbd}{\convexSubdifferential}
\newcommand{\FrechetSubdifferential}{\partial}
\newcommand{\FrSbd}{{\FrechetSubdifferential}}
\newcommand{\miniFrSbd}{{\partial_{\kern-1pt F}^{0}}}
\newcommand{\limitingSubderivative}{\partial_{\ell}}

\newcommand{\lmSbd}{\limitingSubderivative}
\newcommand{\stronglimitingSubderivative}{\partial_{s}}
\newcommand{\slmSbd}{\stronglimitingSubderivative}

\newcommand{\piecewiseConstant}[2]{\overline{#1}_{\kern-1pt#2}}

\newcommand{\piecewiseLinear}[2]{#1_{\kern-1pt#2}}

\newcommand{\pwM}[2]{\widetilde{#1}_{\kern-1pt#2}}
 \def\trait #1 #2 #3 {\vrule width #1pt height #2pt depth #3pt}
\newcommand{\pwN}[2]{#1_{\kern-1pt#2}}
 \def\trait #1 #2 #3 {\vrule width #1pt height #2pt depth #3pt}

 \def\fin{\hfill
         \trait .3 5 0
         \trait 5 .3 0
         \kern-5pt
         \trait 5 5 -4.7
         \trait 0.3 5 0
 \medskip}

\newcommand{\wlim}{w\kern-2pt-\kern-6pt\mathop{\rm lim}\limits}
\newcommand{\nnu}{{\mbox{\boldmath$\nu$}}}

\newcommand{\eps}{\varepsilon}

\newcommand
{\pairing
  }[4]{ \sideset{_{ #1 }}{_{ #2 }}  {\mathop{\langle #3 , #4  \rangle}}}
\newcommand{\variabile}{\,\cdot\,}

\newcommand{\disp}{\displaystyle}
\newcommand{\cx}{\mathcal{X}}
\newcommand{\dcx}{d_{\mathcal{X}}}
\newcommand{\dcxw}{d_{X}^{w}}
\newcommand{\sfl}{\mathscr S}
\newcommand{\att}{\mathcal{A}}
\newcommand{\rest}{\mathcal{Z}(\sfl)}
\newcommand{\semif}{\mathcal{S}}
\newcommand{\geneset}{\mathcal{G}}
\newcommand{\geneatt}{\mathcal{A}}
\newcommand{\esseci}{\mathcal{E}}
\newcommand{\ovesse}{\overline{\esseci}}
\renewcommand{\cx}{X}
\renewcommand{\sfl}{G}
\renewcommand{\dcx}{d_\cx}

\renewcommand{\att}{A}
\renewcommand{\rest}{Z(\sfl)}
\begin{document}
\title{Attractors for gradient flows of non convex
functionals and applications}
\author{Riccarda Rossi \and Antonio Segatti \and Ulisse Stefanelli} 
%
%
\date{Received: date / Revised version: date}
%
\maketitle
\begin{abstract}
 This paper addresses the long-time behaviour of gradient flows
of non convex functionals in Hilbert spaces.
Exploiting the notion of  \emph{generalized semiflows} by
\textsc{J.~M. Ball}, we provide some sufficient conditions for the
existence of a global attractor. The abstract results are
applied to various classes of non
convex evolution problems. In particular, we discuss the long-time behaviour
of solutions of quasi-stationary phase field models and prove the
existence of a global attractor.
\end{abstract}
\section{Introduction}
The aim of this paper is to address the long-time behaviour of strong solutions of the gradient flow equation
\begin{equation}
\label{eq:gflows} \tag{GF}
    u'+ \slmSbd \phi (u)\ni 0
    \quad \text{a.e.\ in } \ (0,+\infty), \quad  u(0)=u_0,
  \end{equation}
associated with the {\it (strong) limiting subdifferential} $\,\slmSbd \phi:\Hilbert\to 2^\Hilbert\,$ of a functional
\begin{equation}\label{eq:basic-prop-phi}
 \phi: \Hilbert\to (-\infty,+\infty] \quad \text{proper and lower semicontinuous,}
\end{equation}
possibly non convex, defined on a separable Hilbert space $\,\Hilbert\,$ with scalar product $\la \cdot, \cdot \ra$ and norm $\,
|\cdot|$. The strong limiting subdifferential $\, \slmSbd \phi\,$ of $\, \phi\,$ is a suitably generalized gradient notion (see below), related to the sequential strong closure in $\, \Hilbert \times \Hilbert \,$ of the graph of the {\it Fr\'echet subdifferential} $\,\partial \phi\,$ of $\, \phi$. The latter is defined, letting $\, D(\phi):=\{u \in \Hilbert : \phi(u) < + \infty\}$, as
\begin{equation}
  \label{eq:def-frechet-subdif}
 \xi \in \partial \phi (v)\quad \text{iff} \quad v \in D(\phi), \ \  \liminf_{w\to
  v}\,\frac{\phi(w)-\phi(v)-\la\xi,{w-v}\ra}{|w-v|}\ge0.
\end{equation}

Existence and approximation results for \eqref{eq:gflows} have been
obtained in \cite{Rossi-Savare-Proc,Rossi-Savare04} for the
non-autonomous situation $\, u'+ \slmSbd \phi (u)\ni f$, where $\, f
\in L^2_{loc}(0,+\infty;\Hilbert) \,$ and an initial datum $\, u_0
\in D(\phi)\,$ is given. The arguments of \cite{Rossi-Savare04}  are
based on the theory of
 \emph{Minimizing Movements} \cite{Ambrosio95,DeGiorgi93} and of \emph{Curves of
 Maximal Slope}
  \cite{Ambrosio-Gigli-Savare04,Cardinali-Colombo-Papalini-Tosques97,DeGiorgi-Marino-Tosques80,Marino-Saccon-Tosques89}, as well as
  on Young measures in Hilbert spaces. A remarkable result of \cite{Rossi-Savare04} is that solutions of \eqref{eq:gflows} fulfil the  {\it energy identity}
\begin{equation}
  \label{eq:11}
  \phi(u(t))+\int_s^t|u'(r)|^2\,dr=\phi(u(s))\quad \forall 0 \leq s \leq t <+\infty.
\end{equation}

The main issue of this paper
is to show that, under suitable assumptions, the set of solutions of \eqref{eq:gflows} admits a \emph{global attractor}. Equation \eqref{eq:11} entails that the functional $\, \phi \,$ decreases along trajectories. Hence, we shall focus our attention on the metric phase space $\, (X,\dcx)\,$ given by
$$
{X}:=D(\phi), \quad
\dcx(u,v):= |u-v| + |\phi(u)-\phi(v)| \  \ \forall u,v \in
{X}.
$$
 Indeed, we define this
  phase space in terms of the functional $\, \phi$, which turns out to be a {\it Lyapunov function} for the system (see
\cite{Rocca-Schimperna04,Segatti04} for some analogous choices).

Due to the possible non convexity of the functional $\, \phi$, uniqueness for \eqref{eq:gflows} may genuinely fail. Hence,
\eqref{eq:gflows} does not generate a semigroup, and we cannot rely
on the well-established theory of \cite{Temam88} for the study of the long-term
dynamics of the solutions. In recent years, several approaches
have been developed in order to
 address the asymptotic behaviour of solutions of
differential problems without  uniqueness. Without any claim of
completeness, we may refer  the reader to, e.g.,  the results by
\textsc{Sell} \cite{Sell73,Sell96},  \textsc{Chepyzhov \& Vishik}
\cite{Chep-Vish95}, \textsc{Melnik \& Valero} \cite{Mell-Vall00}, to
the survey by
 \textsc{Caraballo, Mar\'\i n-Rubio \& Robinson}
 \cite{Caraballo03}, and to the work of \textsc{J.~M. Ball}
 \cite{Ball97,Ball04}.

In particular, we will
 focus here on the theory of \emph{generalized
 semiflows} proposed in \cite{Ball97}. A {generalized
 semiflow} is a family of  functions on $\,[0,+\infty)\,$ taking values
 in the phase space and complying with  suitable existence, stability for time translation,
 concatenation, and upper semicontinuity axioms (see Section \ref{sez:2}).
Within this setting,  it is possible to introduce a suitable notion of
{global attractor} and to
characterize the existence of such an attractor in terms of
boundedness and compactness properties.

The main results of this paper state that, under suitable
assumptions, the solution set to \eqref{eq:gflows} is a {generalized
semiflow} in the space $\, (X,\dcx)$ (Theorem \ref{teor:1}), and
that it possesses a global attractor (Theorem \ref{teor:2}). The key
point in our proofs involves passing to the limit in the energy
identity \eqref{eq:11} by means of compactness results for Young
measures in Hilbert spaces.

A large part of the paper is devoted to a discussion on
the applications of the aforementioned abstract results
to evolution problems with a gradient flow structure.
 First of all, we show the existence of a global
 attractor in the case of $\, \phi \,$ being a suitable
  perturbation of a convex functional. In fact, our results apply
  to $\, C^1 \,$ perturbations as well as to dominated
  concave perturbations of convex functionals (see Section \ref{applications} below).

Secondly, we investigate the long-time
behaviour of a class of solutions of
the so-called {\it quasi-stationary phase field} system
\begin{equation}
\label{e:corazon}
\begin{cases}
  \begin{aligned}
      \partial_t (\teta +\chi) -\Delta \teta
      &=0,
      \\
      F'(\chi) &=\teta,
      \end{aligned}
      \end{cases}
\end{equation}
in $\,\Omega \times (0,+\infty)$, where $\,\Omega\,$ is a bounded
domain and $\, F' \, $ is the G\^ateaux derivative of a functional
$\, F$, (possibly neither smooth nor  convex). The model
\eqref{e:corazon} arises as a suitable generalization of the
(formal) quasi-stationary asymptotics of the standard parabolic phase
field model \cite{Caginalp86}, which describes the phase transition
in an ice-water system. In this connection, $\, \teta \,$ is the
relative temperature of the system, while the order parameter $\,
\chi \,$ yields the local proportion of the liquid versus the solid
phase. The usual choice for $\, F \,$ is
\begin{equation}\label{cagi}
 F(\chi):=\frac12\int_\Omega |\nabla \chi|^2dx + \frac14\int_\Omega(\chi^2 - 1)^2dx.
\end{equation}

The existence of solutions of some initial and boundary value
problem for \eqref{e:corazon} with $F$ as in  \eqref{cagi} was
firstly proved by {\sc Plotnikov \& Starovoitov} in
 \cite{Plotnikov-Starovoitov93}. The latter paper addresses the case of homogeneous Dirichlet conditions on $\,\teta\,$ and homogeneous Neumann
  conditions on $\,\chi\,$,  and exploits a compactness method and a non
  standard unique continuation result. Let us mention that the latter technique heavily relies on the
  precise form of \eqref{cagi} and cannot be easily extended to a more general situation.
   A second result in the direction of the existence of a solution of
   \eqref{e:corazon}-\eqref{cagi}
    is due to {\sc Sch\"atzle} \cite{Schatzle00}. The argument devised in \cite{Schatzle00}
  for
  proving existence for \eqref{e:corazon}-\eqref{cagi}, supplemented
  with homogeneous Neumann-Neumann boundary conditions
  on both $\,\teta\,$ and $\,\chi$, exploits some spectral analysis results
  and the analyticity of $\,\chi \mapsto (\chi^2 -1)^2/4$.
  Once again, this technique is especially tailored to
   the form of \eqref{cagi} and cannot be reproduced
   for general functionals $\, F$.
   We may observe (see, e.g., \cite{Visintin96}) that,
    indeed, \eqref{e:corazon} stems as the formal
    gradient entropy flow for the phase field system.
     The latter {gradient flow} approach to the
     problem of existence of solutions of
   \eqref{e:corazon} has been fully considered in detail by {\sc Rossi \& Savar\'e}
   in \cite{Rossi-Savare-Proc,Rossi-Savare04}. In
particular, the existence results in
\cite{Rossi-Savare-Proc,Rossi-Savare04} provide a unified frame and
 extend the previous aforementioned contributions on existence
results for quasi-stationary phase fields.

The gradient flow structure of \eqref{e:corazon} is enlightened by introducing the variable $\,u:= \teta+\chi$.
Following \cite{Rossi-Savare04}, one can
 rigorously prove that \eqref{e:corazon}, along with the boundary conditions
$\,u-\chi =\partial_n \chi =0 \,$ on $\, \partial \Omega\,$ for instance, may be
interpreted  as the {gradient flow} equation  in the Hilbert
space $\,H^{-1}(\Omega)\,$
 of the
functional
 $\,\phi: H^{-1}(\Omega) \to (-\infty,
+\infty]\, $defined by
\begin{equation}
 \label{eq:funz-qstat}
 \phi(u):= \inf_{\chi \in H^1 (\Omega)}\left( \frac12 \int_\Omega |u-\chi|^2\,dx + F(\chi)\right), \quad \text{$
D(\phi):=L^2
 (\Omega).$}
\end{equation}
 Namely, in \cite{Rossi-Savare04} it has been shown that the solutions of \eqref{eq:gflows},
  with the choice \eqref{eq:funz-qstat} for $\, \phi$, provide a family of solutions of \eqref{e:corazon}, supplemented with homogeneous Dirichlet-Neumann conditions.

In Section \ref{sez:6} we show that the solutions of \eqref{e:corazon} arising from the gradient flow of the functional $\, \phi \,$ \eqref{eq:funz-qstat} indeed form a generalized semiflow, which admits a global attractor. Let us stress that this gradient flow approach does not provide the description of the long-term behaviour of the \emph{whole} set of
solutions of \eqref{e:corazon}, but is rather concerned with a proper subclass of solutions. Moreover, we present some
 result on the long-time behaviour of solutions in the weakly coercive case
  of Neumann-Neumann boundary conditions. The latter situation
is more delicate, since \eqref{e:corazon} fails to have a
gradient flow structure. However,
the existence of solutions may be deduced by suitably approximating
the system by means of more regular problems of gradient flow type. The latter
 approximation procedure has been in fact
 detailed in \cite{Rossi-Savare04} and is here
reconsidered from the point of view of the long-time dynamics. In
particular, in the weakly coercive case, the set of solutions of
\eqref{e:corazon} obtained as mentioned above fails to be a
generalized semiflow. Nevertheless, by slightly extending
\textsc{Ball's} theory (see Section \ref{sez:2}), in
Section \ref{s:5.2.2}
 we  prove the existence of a suitable notion
of weak global attractor for the weakly coercive problem as well.
Indeed, denoting by $\, \geneatt_\lambda \,$ for $\, \lambda\in (0,1)\,$
the family of global attractors for the approximate problems and by
 $\geneatt$ the weak global attractor for the (weakly coercive) limit
problem, we also prove in Section \ref{sec:approx-attractor} \, the
convergence of $\geneatt_\lambda$ to $\geneatt$, as the
approximation parameter $\lambda\down 0$,
 with respect to a suitable  Hausdorff
semidistance.
\paragraph{\bf Plan of the paper.} We present some introductory material in
Section \ref{preliminaries}. In particular, Section \ref{egf}
concerns the existence of solutions of \eqref{eq:gflows} and reports
a result from \cite{Rossi-Savare04},
 while in Section \ref{sez:2} we  recall some results on
 \textsc{Ball's}
 theory on generalized semiflows and develop additional material,
in the direction of studying a weak semiflow structure and a weak
notion of attractor. Section \ref{sez:4} contains the statement and
the proof of our main abstract results (Theorems \ref{teor:1} and
\ref{teor:2}). The ensuing Sections
\ref{applications}-\ref{sez:3-qstat-evol} are devoted to
applications. In particular, Section \ref{applications} deals with
the long-time behaviour of solutions of gradient flows of suitably
perturbed convex functionals. We consider both the case of $\, C^1
\,$ perturbations and that of (suitably dominated) concave
perturbations. Moreover, some PDE examples are provided within these
classes of problems. Section \ref{sez:3-qstat-evol} is focused on
the long-time behaviour of solutions of the quasi-stationary phase
field model \eqref{e:corazon}. Since our approach to the long-time
behaviour of \eqref{e:corazon} is substantially based on the
gradient flow strategy developed in \cite{Rossi-Savare04}, we shall
briefly recall  the techniques and results of the latter paper in
Section \ref{sez:3-qstat-evol}. Then, Theorems \ref{teor:1} and
\ref{teor:2} are applied to the quasi-stationary
problem \eqref{e:corazon} in Section \ref{sez:6}. 
\paragraph{\bf Acknowledgment.} The authors would like to thank Prof.
Giuseppe Savar\'{e} for some va\-lua\-ble and inspiring
conversations.

\section{Preliminary results}\label{preliminaries}

\subsection{Existence for gradient flows of non convex functionals}\label{egf}

In this section we gain some insight into
an existence result for \eqref{eq:gflows}
that has been obtained in  \cite{Rossi-Savare04}.
To this aim, let us start by reviewing the
 results on gradient flows in the convex case.
 Given  $\, T >0 \,$ and $\, f: (0,T) \to \Hilbert $, we consider the problem
\begin{equation}
\label{eq:cauchy-generico}
    u'(t)+ \partial \phi (u(t))\ni f(t)
    \quad \forae \, t \in  (0,T), \quad  u(0)=u_0.
  \end{equation}
 When $\,\phi\,$
 is a \emph{convex} functional, the Fr\'echet subdifferential of $\, \phi \,$
coincides with the
 subdifferential $\cvSbd \phi$ of $\phi$ in the sense of Convex
 Analysis (so we shall use the same notation for both
 subdifferential notions). The latter is defined by
\begin{equation}
  \label{eq:def_subdifferential_convex}
  \xi\in \cvSbd\phi(v)\quad\text{iff}\quad
  v\in D(\phi), \quad
  \phi(w)-\phi(v)-\duality\xi{w-v}\ge0
  \quad\forall\,w\in \Hilbert.
\end{equation}
The literature on existence, uniqueness,
regularity, and approximation of solutions of
\eqref{eq:cauchy-generico} is well-established and dates back to the
early 70s (see the seminal references
\cite{Brezis71,Brezis73,Crandall-Pazy69,Komura67}). In particular,
it is well-known that, if $\,u_0 \in D(\phi)\,$ and $\,f \in
L^2(0,T;\Hilbert)$, then the Cauchy problem
\eqref{eq:cauchy-generico} admits a {unique} solution $\,u \in H^1
(0,T;\Hilbert)$, which complies with the energy identity
\begin{equation}
  \label{eq:110}
  \phi(u(t))+\int_s^t|u'(r)|^2\,dr=\phi(u(s)) +\int_s^t \la f(r), u'(r) \ra
  dr\quad \forall 0 \leq s \leq t \leq T.
\end{equation}
Indeed, relation
\eqref{eq:110} follows from the \emph{chain rule}
property of convex subdifferentials, i.e.,
\begin{equation}
  \label{eq:10}
  \begin{gathered}
    \text{if}\ \
    u\in H^1(0,T;\Hilbert),\ \xi\in L^2(0,T;\Hilbert),\
    \xi(t)\in \partial \phi(u(t))\ \forae\ t\in (0,T),\\
    \text{then}\quad \phi\circ u \in AC(0,T),\quad
    \tfrac d{dt}\phi(u(t))=\duality{\xi(t)}{u'(t)}\ \forae \ t \in (0,T).
  \end{gathered}
\end{equation}
In fact,
the \emph{strong-weak closure} of $\cvSbd \phi$
in the sense of graphs, i.e.,
\begin{equation}
  \label{eq:strong-weak-closure}
 u_n \to u,  \ \xi_n \weakto \xi \  \text{in $\Hilbert$,} \
\xi_n \in \FrSbd \phi(u_n) \ \forall n \quad \Rightarrow \quad
\phi(u_n) \to \phi(u), \, \xi \in \FrSbd \phi(u),
\end{equation}
the elementary  continuity property
\begin{equation}
\label{eq:continuity-convexity} u_n \to u, \quad \sup_n |\FrSbd
\phi^\circ(u_n)| <+\infty \quad \Rightarrow \phi(u_n) \to \phi(u),
\end{equation}
(where we use the notation $\, |A^\circ|:=\inf_{a \in A} |a|$ for
all
 non-empty sets $\, A \subset
\Hilbert$),
 and the chain rule
 \eqref{eq:10}
play a crucial role in the   proof of the existence
 of solutions of \eqref{eq:cauchy-generico}.
  Furthermore,
the long-time behaviour of \eqref{eq:cauchy-generico}  from the
point of view of the theory of universal attractors, see e.g. {\sc Temam}
\cite{Temam88}, is quite well-understood, even in the non autonomous
case (see also \cite{siyk98}).

Let us now turn  to the case of a proper, lower semicontinuous, and
non convex functional $\,\phi$, cf. \eqref{eq:basic-prop-phi}. One
shall observe that, even in the non convex case, for any $\,u \in
D(\FrSbd \phi) $ the Fr\'echet subdifferential $\, \FrSbd \phi(u)
\,$ is  a convex subset of  $\Hilbert$. On the other hand, the
elementary example $\,\phi(x):= \min \{(x-1)^2, (x+1)^2\}\,$ (with
$\,\FrSbd \phi(x):= 2(x+1) \,$ for $\,x <0$, $\,\FrSbd \phi(x):=
2(x+1)\, $ for $\,x>0\,$, but $\,\FrSbd \phi(0)=\emptyset$) shows
that, unlike the convex case (see \eqref{eq:strong-weak-closure}),
the graph of the Fr\'echet subdifferential of a non convex
functional may not be strongly-weakly closed.

Therefore, following \cite{Rossi-Savare04}, we define the
\emph{strong limiting subdifferential} $\,\slmSbd\phi\,$ of $\, \phi
\,$ at  a point $\,v\in D(\phi)$  as the set of the
   vectors  $\,\xi \,$ such that  there exist sequences
\begin{equation}
  \label{eq:def-limiting-subdif}
v_n,\xi_n\in\Hilbert\quad\text{with }\quad \xi_n\in
\FrSbd\phi(v_n),\
  v_n\to v,\
  \xi_n \to \xi, \  \phi(v_n) \to \phi(v),
\end{equation}
as $\,n\to+\infty$. Furthermore, we define  the
 \emph{weak limiting subdifferential}
 $\,\lmSbd \phi\,$ of $\,\phi\,$ at  $\,v\in D(\phi)\,$ as the set
 of all vectors $\,\xi\,$ such that
  there exist sequences
\begin{equation}
  \label{eq:def-weak-limiting-subdif}
v_n,\xi_n\in\Hilbert\quad\text{with}\quad \xi_n\in
\FrSbd\phi(v_n),\
  v_n\to v,\
  \xi_n \weakto \xi, \  \sup_n \phi(v_n) <+\infty.
\end{equation}
 Of course,  $\lmSbd \phi$ and $\slmSbd \phi$ reduce to the
subdifferential $\FrSbd \phi$ of $\phi$ in the sense of Convex
Analysis whenever $\phi $ is convex, due to
\eqref{eq:strong-weak-closure} and \eqref{eq:continuity-convexity}.

 Note that
 the strong limiting subdifferential $\slmSbd \phi$
of $\phi$   fulfils this closure property:
\begin{equation}
\label{eq:closure-properties-limiting}
\begin{gathered}
 \forall\, \{u_k\}, \{\xi_k\}  \ \text{such that}
\ \
  u_k \to u, \, \xi_k \to \xi, \,   \phi(u_k) \to \phi(u), \,
\text{as $k \up+ \infty$,}\\ \quad  \xi_k
 \in \slmSbd \phi(u_k) \ \forall k \in \N,
 \text{then} \ \xi \in \slmSbd \phi(u).
 \end{gathered}
\end{equation}
Instead,
 $\lmSbd \phi$  is not strongly-weakly closed in the
sense of graphs. Actually, $\lmSbd \phi$ can be characterized as a
version of the  \emph{strong-weak} closure of $\slmSbd \phi$, as
 the following result
 shows.
\begin{lemma}
\label{lemma:closure-properties-limiting} Let $\phi: \Hilbert \to
(-\infty, +\infty]$ comply with \eqref{eq:basic-prop-phi}. Then, for
any $u \in \Hilbert$
\begin{gather}
 \xi \in  \lmSbd\phi(u)
\ \Longleftrightarrow \ \exists \, \{u_k\}, \{\xi_k\} \subset
\Hilbert: \ u_k \to u, \, \xi_k \weakto \xi, \, \sup_k \phi(u_k)
<+\infty,\nonumber\\
 \, \xi_k
 \in \slmSbd \phi(u_k) \, \forall k \in \N\label{eq:closure-properties-weak-limiting}
\end{gather}
i.e., $\lmSbd \phi$ coincides with  the (sequential) strong-weak
closure of $\slmSbd \phi$ along sequences with bounded energy.
\end{lemma}
 \begin{proofx} The left-to-right implication in
 \eqref{eq:closure-properties-weak-limiting} follows immediately from
 the definition of $\lmSbd \phi$,
  noting that $\FrSbd
 \phi(u) \subset \slmSbd \phi(u) $ for any $u \in \Hilbert$.
  In
 order to prove the converse implication, we recall that in separable
 Hilbert spaces (more in general, in reflexive spaces and dual of
separable spaces, cf. \cite[Chap. 3]{Brezis83}), it is possible to
 introduce a norm $\trepar\cdot\trepar$,
  and thus a metric, inducing weak convergence on
every bounded set.
 Thus,
  let us fix a sequence $\{(u_k, \xi_k)\} $ as in
 \eqref{eq:closure-properties-weak-limiting}: necessarily, there
exists $M \geq 0$ such that $|\xi_k| \leq M$, and $\xi_k \weakto
 \xi$ may be rephrased as $\trepar \xi_k -\xi \trepar \to 0.$
  In order to prove that the limit pair $(u,\xi)$ fulfils $\xi \in
 \lmSbd \phi(u),$ we are going to construct by a diagonalization
 procedure a sequence $\{(v_k, \omega_k)\} \subset \Hilbert \times
 \Hilbert$ such that
\begin{equation}
  \label{eq:final-aim}
 v_k \to u, \ \omega_k \weakto \xi \ \text{as $k \up +\infty$}, \quad
 \sup_k \phi(v_k) <+\infty, \quad \omega_k \in \FrSbd \phi(v_k) \
 \forall k \in \N.
  \end{equation}
 Note that the relation
  $\xi_k \in \slmSbd \phi(u_k) $ for all $k\in \N$ can
 be rephrased in the following way: for any $k \in \N$ there exist
 sequences  $\{u_{n}^k\}, $ $\{\xi_{n}^k\}, \subset \Hilbert $ with
 $$
 u_{n}^k \to u_k, \ \xi_{n}^k \to \xi_k, \ \phi(u_{n}^k) \to
 \phi(u_k) \quad \text{as $n \up +\infty$, \  and} \quad \xi_{n}^k
 \in \FrSbd \phi(u_{n}^k) \ \forall n \in N.
 $$
 In particular, for any $k \in \N$ we may find  $n(k) \in \N$ such
 that
 $$
 |u_{n(k)}^k - u_k| + |\xi_{n(k)}^k - \xi_k| + |\phi(u_{n(k)}^k) -
 \phi(u_k)| \leq \frac{1}{k}.
 $$
 Then, let us set  $v_k:= u_{n(k)}^k $ and $\omega_k :=\xi_{n(k)}^k$.
 Obviously, $\omega_k \in \FrSbd \phi(v_k),$ $v_k \to u$ as $k \up
 +\infty$, and  $\sup_k \phi(v_k) <+\infty. $
  On the
 other hand, we remark that the sequence $\{\omega_k\}$ lies in a
 bounded set of $\Hilbert$, since for all $k \in \N$
  $|\omega_k | \leq |\omega_k -\xi_k| + |\xi_k| \leq 1+ M$.
  Therefore, \eqref{eq:final-aim} follows by noting that
 \begin{equation*}
  \trepar \omega_k - \xi \trepar \leq \trepar \omega_k -\xi_k
  \trepar + \trepar \xi_k -\xi
  \trepar \leq \sqrt{2}|\omega_k -\xi_k| + \trepar \xi_k -\xi
  \trepar \to 0 \quad \text{as $k \up +\infty.$}
 \end{equation*}
 \end{proofx}

 Under the assumption that $\lmSbd \phi$ satisfies a  chain rule
property analogous to the chain rule \eqref{eq:10} of the
subdifferential of Convex Analysis, in \cite{Rossi-Savare04}
existence and approximation results have been obtained for
\eqref{eq:gflows}, supplemented with some  initial datum $u_0 \in
D(\phi)$ and source term $f$. \label{sez:3-exist} Let us now recall
one of the existence results proved in \cite{Rossi-Savare04}.
\begin{theorem}
\label{teor:exist-gflows}
 Suppose that  $\phi : \Hilbert \to
(-\infty,+\infty]$ complies with \eqref{eq:basic-prop-phi}, with the
\emph{coercivity} assumption
\begin{equation}
    \label{eq:coerc-ass}
    \exists\,\kappa \geq  0:\quad
    v \mapsto \phi(v)+\kappa |v|^2\quad
    \text{{has compact sublevels,}}
  \tag{\textsc{comp}}
  \end{equation}
and with the \emph{chain rule} condition
 \begin{equation}
    \label{eq:chain-rule-energy}
    \begin{gathered}
      \text{\emph{if} $v\in H^1(a,b;\Hilbert)$, $\xi\in L^2(a,b;\Hilbert)$, $\xi\in
        \lmSbd\phi(v)$ a.e.\ in $(a,b)$,}\\
      \text{and $\phi\circ v$ is bounded, {then} $\phi\circ
        v\in AC(a,b)$ and}
      \\      \frac d{dt}\phi(v(t)) =\la \xi(t),v'(t)\ra
      \quad\text{for a.e. }\ t \in  (a,b).
    \end{gathered}
        \tag{\textsc{chain}}
    \end{equation}
Then, for any $u_0 \in D(\phi)$, $T>0$   and $f \in
L^2(0,T;\Hilbert)$
the Cauchy problem
$$
u'(t)+ \slmSbd \phi (u(t))\ni f(t)
    \quad \text{a.e.\ in }(0,T), \quad  u(0)=u_0,
$$
admits a solution $u \in H^1(0,T;\Hilbert)$. Moreover, one has the \emph{energy
identity}
\begin{equation}
      \label{enid2}
      \int_{s}^{t}
      |u' (\sigma)  |^{2} \, d \sigma + \phi(u (t))
       =    \phi(u(s)) + \int_{s}^{t} \la f(\sigma),u' (\sigma) \ra \, d\sigma
      \quad\forall 0 \leq s\leq t\leq T .
    \end{equation}
\end{theorem}

The chain rule \eqref{eq:chain-rule-energy}, which is indeed
classical in the convex case \eqref{eq:10}, holds true in a variety
of non convex situations as well. First of all,
\eqref{eq:chain-rule-energy} is fulfilled by $\,C^1\,$ perturbations
of convex functionals. In particular, letting $\phi  = \phi_1 +
\phi_2$, where $\, \phi_1 \,$ is convex and $\, \phi_2 \,$ is
$\,C^1$, and exploiting Lemma \ref{c1perturb}, one readily checks
that $\,\lmSbd\phi=\partial\phi_1 + D\phi_2\,$ and
\eqref{eq:chain-rule-energy} follows. A second class of  functionals
complying with the chain rule \eqref{eq:chain-rule-energy} is
provided by do\-mi\-nated concave perturbations of (convex)
functionals. In particular, in \cite[Thm. 4]{Rossi-Savare04} it is
proved  that \eqref{eq:chain-rule-energy} holds for
  all proper, lower semicontinuous functionals $\,\phi\,$
admitting the decomposition
\begin{equation}
  \begin{gathered}
    \phi=\psi_1 -\psi_2 \text{ in }D(\phi),
    \quad\text{with}\quad
    \text{$\psi_1:D(\phi)\to\R$ l.s.c.\ and satisfying
      \eqref{eq:chain-rule-energy}},\\
    \psi_2:\conv\big(D\phi\big)\to \R
    \text{ convex and l.s.c.\ in $D(\phi)$},
    \quad
    D(\lmSbd\psi_1)\subset
    D(\FrSbd\psi_2),
    \end{gathered}
\label{eq:decomposition}
\end{equation}
(where $\,\conv\big(D\phi\big)\,$ denotes the convex hull of $\,D(\phi)$),
and fulfilling
\begin{equation}
        \label{eq:17}
        \begin{gathered}
          \forall M \geq 0 \quad \exists \ \rho<1,\, \gamma \geq 0 \quad
          \text{such that}\quad
          \sup_{\xi_2\in\FrSbd\psi_2(u)}
          |\xi_2|\le \rho
          |(\lmSbd\psi_1(u))^\circ|+\gamma
          \\
          \quad
          \text{for every $ u \in
            D(\lmSbd \psi_1)$ with $\max(\phi(u),|u|) \leq M $}.
      \end{gathered}
      \end{equation}
 Namely, if $\,\psi_1\,$ is itself convex, we are
requiring the domain of $\,\FrSbd \psi_1\,$ to be included in
$\,D(\FrSbd\psi_2)$. This in fact implies that $\,\FrSbd \psi_1\,$
somehow dominates $\,\FrSbd\psi_2$.





\subsection{Generalized semiflows}
\label{sez:2}
 For the reader's convenience, we recall here
 the main definitions and results on the theory of
\emph{attractors} for \emph{generalized semiflows},
 closely
following  \cite{Ball97}.
 Our final aim is to apply Ball's theory to
the Cauchy problem \eqref{eq:gflows},  and slightly extend it in
view of applications.
 \par\noindent {\bf Notation.} Let $\,(\cx,
\dcx)\,$ be a (not necessarily complete) metric space. We recall
that the \emph{Hausdorff semidistance} or \emph{excess} $\,e(A,B)\,$
of two non-empty subsets $A, \, B \subset \cx$  is given by
\linebreak$\,e(A,B):= \sup_{a \in A} \inf_{b \in B} \dcx(a,b)$. For
all $\eps>0$, we also denote by $\, B(0,\eps)\,$ the ball $\, B(0,\eps) :=\{x \in X \ : \ \dcx(x,0) < \eps\}$,
and by $N_\eps (A):= A + B(0,\eps)$ the $\eps$-neighborhood of  a
subset $A$.
\begin{definition}[Generalized semiflow]
\label{def:generalized-semiflow} A \emph{generalized semiflow}
$\,\sfl\,$ on $\,\cx\,$ is a family of maps $\,g:[0,+\infty) \to
\cx\,$ (referred to as ``solutions"), satisfying:
\begin{description}
\item[(H1)] \emph{\bf (Existence)} for any $\,g_0 \in \cx\,$ there exists at
least one $\,g \in \sfl$ with $\,g(0)=g_0,$
 \item[(H2)]
\emph{\bf (Translates of solutions are solutions)} for any $\,g \in
\sfl\,$ and $\,\tau \geq 0$, the map $\,g^\tau (t):=g(t+\tau),$ $\,t
\in [0,+\infty),$ belongs to $\,\sfl$,
 \item[(H3)]
\emph{\bf (Concatenation)} for any $\, g \,$, $\, h \in \sfl\,$ and $\,t
\geq 0\,$  with $\,h(0)=g(t)$, then  $\,z \in \sfl$, $\, z\,$ being
the map defined by $ z(\tau):=  g(\tau)\,$ if $\,0 \leq \tau \leq t,$ and
$\, h(\tau-t)\,$ if $\, t <\tau$.
\item[(H4)] \emph{\bf (Upper-semicontinuity w.r.t. initial data)} If
$\,\{g_n \} \subset \sfl\,$ and $\,g_n (0) \to g_0, $ then there
exists a subsequence $\,\{g_{n_k}\}\,$ of $\,\{g_n \}\,$ and $\,g
\in \sfl\,$ such that $\, g(0)=g_0$ and  $\,g_{n_k}(t) \to g(t)\,$
for all $\,t \geq 0.$
\end{description}
\end{definition}

 The
application of the theory of generalized semiflows  to suitable
 classes of differential problems is often delicate. Indeed, one
usually  needs to choose carefully the correct notion of solution of
the problem in  order to check  the validity of the properties
(H1)-(H4). This process  may  not be straightforward whenever one
considers some suitably weak notion of solvability. On the one hand,
solutions have indeed to be weak enough to comply with
 (H1) (assumption (H2) is generally easy to meet in actual situations). On
the other hand, the notion of solution has to be \emph{robust}
enough in order to fulfil (H4). This robustness may turn out to be
in conflict with (H3). For instance, this may occur when the
existence of weak solutions of a differential problem  is
proved by approximation  (like e.g. for the solutions of the
quasi-stationary phase field  Problem \ref{quasi-stationary} in the
weakly coercive case, cf. Theorem \ref{teor:exist-w-coerc}). Then,
one is naturally led to define the candidate semiflow  as the set of
all solutions which are limits in a suitable  topology of sequences
of
 approximate solutions. Axioms (H1) and (H2)
 will be trivially checked, and,   if the aforementioned topology is strong
 enough, one can hopefully verify (H4) as well. However, due to this
 approximation procedure, the concatenation in  (H3) may not hold
 (the approximating sequences may
  not have the same indices). This is particularly the case of the
 set of \emph{limiting energy solutions} of Problem \ref{quasi-stationary} in
the weakly coercive case (cf. Definition
\ref{def:limiting-energy-solu}).

Therefore, in the setting of the phase space $(\cx, \dcx) $
 we aim at partially extending the standard theory of generalized semiflows
 to the
case of a  non-empty set $\geneset $ of  functions $g:[0,+\infty)
\to X$,
 complying with (H1), (H2), (H4), but not necessarily with (H3).
In this framework, we shall    introduce  a weakened
 notion of attractor,  for
objects which are slightly more general than semiflows. Before
moving on, let us explicitly stress that we do not claim originality
for the notion of {\em weak generalized semiflow} we present
below.
 Indeed, the
possibility of studying the long-time dynamics of differential
systems by considering (multivalued) solution operators fulfilling
\eqref{eq:inclu} has been recently considered in \cite{Mell-Vall98,Mell-Vall00}. In particular, this  multivalued approach has also
been applied  to the standard phase field system by {\sc Kapustyan,
Melnik \& Valero} \cite{Kapp-Mell-Vall03}.
\begin{definition}[Weak generalized semiflow]
 \label{def:w-generalized-semiflow} We say
that a non-empty family $\geneset$  of maps
 $g : [0,+\infty) \to X$ is  a \emph{weak generalized semiflow}
 on $\,\cx\,$ if $\,\geneset$ complies with the properties
  {\rm (H1), (H2),} and~{\rm (H4)}.
\end{definition}
\noindent
{\bf Continuity property (C4).} We say that a (weak) generalized
semiflow fulfills $(C4)$ if
 for any $\,\{g_n \} \subset \geneset\,$ with  $\,g_n (0) \to g_0, $ there exists
a subsequence $\,\{g_{n_k}\}\,$ of $\,\{g_n \}\,$ and $\,g \in
\geneset\,$ such that $\,g(0)=g_0\,$ and $\, g_{n_k} \to g\,$
uniformly on the compact subsets of $\,[0,+\infty),$

\noindent {\bf Orbits, $\,\omega$-limits, and attractors.} Given a
{weak generalized semiflow $\,\geneset\,$ on $\,\cx\,$, we may
introduce for every $\,t \geq 0 $ the operator $\,{T}(t): 2^\cx \to
2^\cx\,$ defined  by
 \begin{equation}
 \label{eq:operat-T}
  {T}(t)E:=\{ g(t) \ : \
g \in \geneset \ \  \text{with} \ \ g(0) \in E\}, \quad E \subset
\cx.
 \end{equation}
 The family of operators $\{{T}(t)\}_{t \geq 0}$
  fulfils  the following
 property
 \begin{equation}
   \label{eq:inclu}
   { T}(t+s) B \subset { T}(t){ T}(s) B \quad \forall s,t \geq 0 \quad\forall \, B \subset
   X,
 \end{equation}
 and in general does not define a \emph{semigroup}  on
the power set $\,2^\cx$. Note that  \eqref{eq:inclu} improves to a
semigroup relation when $\geneset$ is a {\em generalized semiflow}}.
Given  a solution $\, g \in \geneset$, we introduce the
\emph{positive orbit of $\,g\,$} as the set $\,\gamma^+ (g):= \{
g(t) \ : \ t \geq 0\}$, while its \emph{$\,\omega$-limit}
$\,\omega(g)\,$ is defined by
$$
\omega(g):= \{ x \in \cx \ : \ \exists \{t_n\}, \ t_n \to +\infty, \
\text{such that} \ \ g(t_n) \to x \}.
$$
We say that $\,w: \R \to \cx\,$ is  a \emph{complete orbit} if, for
any $\,s \in \R$, the translate
 map $\,w^s \in \geneset\,$ (cf. (H2)).
Moreover, we may consider the \emph{positive orbit of a subset $\,E
\subset \cx$}, i.e. the set $\, \gamma^+ (E):=\cup_{t \geq 0} T(t)E=
\cup \{\gamma^+ (g) \ : \ g \in \geneset, \ g(0) \in E\}, $ and, for
every $\,\tau \geq 0$, we define $\, \gamma^\tau (E):=\cup_{t \geq
\tau} T(t)E= \gamma^+ (T(\tau)E).$ Finally, the
\emph{$\omega$-limit} of $E$ is defined as
\begin{gather}
\omega(E):=\big\{
 x \in \cx \ : \ \exists \{g_n\} \subset \geneset \
 \text{such that $\,\{g_n (0)\} \subset E$,} \nonumber\\
 \qquad\qquad \text{$\,\{g_n (0)\}\,$ is bounded,
and} \ \  \exists t_n \to +\infty \ \text{with $\, g_n (t_n) \to x$}
\big\}.\nonumber
\end{gather}
Given subsets $\,U, E \subset \cx$, we say that $\,U \,$
\emph{attracts} $\,E\,$ if $\,e(T(t)E,U) \to 0\,$ as $\,t \to
+\infty$. Further, we say that $\,U \,$ is \emph{positively
invariant} if $\,T(t)U \subset U\,$ for every $\,t \geq 0$, that
$\,U\,$ is \emph{quasi-invariant} if for any $\,v \in U\,$ there
exists a complete orbit $\,w\,$ with $\,w(0)=v\,$ and $\,w(t) \in
U\,$ for all $\,t \in \R$, and finally that $\,U\,$ is
\emph{invariant} if $\,T(t)U = U\,$ for every $\,t \geq 0\,$
(equivalently, if it is both positively and quasi-invariant).
\begin{definition}[Weak  Global Attractor and Global Attractor.]
\label{def:w-attract}
 Let $\geneset$ be  a {\em weak generalized semiflow}.
  We say that a non-empty set $\geneatt$
is a \emph{weak global attractor} for $\geneset$ if it is compact,
quasi-invariant, and attracts all the bounded sets of $X$.
We say that  a set $\,\att \subset \cx\,$ is a \emph{global
attractor} for a  generalized semiflow $\,\sfl\,$ if $\,\att\,$ is
compact, invariant, and attracts all the bounded sets of $\,\cx$.
\end{definition}
 \noindent The price of dropping the semigroup property
for ${T}$ consists in the fact that the notion of weak attractor
introduced above will be quasi-invariant but will fail to be
invariant. Moreover, we may observe that  a weak global attractor
(if existing), is minimal in the set of the closed subsets of $X$
attracting all bounded sets, hence it is unique, cf.
\cite{Mell-Vall98}.
 \par\noindent {\bf Compactness and dissipativity
properties.}
 Let $ \geneset $ be a \emph{weak generalized semiflow.} We say that $
 \geneset$ is
\begin{description}
\item \emph{\bf eventually bounded} if for every bounded $\,B \subset
\cx\,$ there exists $\,\tau \geq 0\,$ such that $\, \gamma^\tau
(B)\,$ is bounded,
\item \emph{\bf point dissipative} if there exists a bounded
set $\,B_0 \subset \cx\,$ such that for any $\,g \in \geneset\,$
there exists $\,\tau \geq 0\,$ such that $\,g(t) \in B_0 \,$ for all
$\,t \geq \tau,$
\item \emph{\bf compact} if for any sequence $\,\{g_n \} \subset \geneset\,$
with $\,\{g_n (0)\}\,$  bounded,  there exists a subsequence
$\,\{g_{n_k} \}\, $ such that $ \,\{g_{n_k}(t) \} \,$ is convergent
for any $\,t
>0.$
\end{description}
The  notions we have just  introduced are not independent one from
each other  cf.
\cite[Prop. 3.1 \& 3.2]{Ball97}.

\noindent {\bf Lyapunov function.}  The notion of \emph{Lyapunov
function} can be introduced starting from the following definitions:
we say that a complete orbit $\,g \in \geneset\,$ is
\emph{stationary} if there exists $\,x \in \cx\,$ such that
$\,g(t)=x\,$ for all $\,t \in \R\,$ - such $\,x\,$ is then called a
\emph{rest point}.  Note that the set of rest points of  $\,\geneset
$, denoted by $\,Z(\geneset)$, is closed in view of {(H4)}. A
function $\,V: \cx \to \R\,$ is said to be a \emph{Lyapunov
function} for $\,\geneset\,$ if: $\,V\,$ is continuous, $\,V(g(t))
\leq V(g(s))\,$ for all $\,g \in \geneset\,$ and $\,0 \leq s \leq
t\,$ (i.e., $\,V\,$ decreases along solutions), and, whenever the
map $\,t \mapsto V(g(t))\,$ is constant for some complete orbit $\,
g\,$, then $\, g \,$ is a stationary orbit.

Finally, we say that a global attractor $\,\att\,$ for $\,\sfl\,$ is
\emph{Lyapunov stable} if for any $\,\eps>0\,$ there exists
$\,\delta>0\,$ such that for any $\,E \subset \cx\,$ with
$\,e(E,\att) \leq\delta, $ then $\,e(T(t)E,\att) \leq\eps\,$ for all
$\,t \geq 0.$
\par\noindent {\bf Existence of the global attractor.} We recall the
main results from {\sc Ball}  \cite{Ball97} (cf. Thms. 3.3, 5.1, and
6.1 therein), which provide criteria  for the existence of a global
attractor $\,\att\,$ for a generalized semiflow $\,\sfl$. More
precisely, Theorem \ref{thm:ball1} gives a characterization of
$\,\att$, whereas Theorem \ref{thm:ball2} states a sufficient
condition for the existence of $\,\att\,$ in the case in which
$\,\sfl\,$ also admits
 a Lyapunov function.
\begin{theorem}
\label{thm:ball1} An eventually bounded, point dissipative, and compact
 generalized semiflow $\,\sfl\,$ has a global
attractor.
 Moreover, the attractor $\,\att\,$ is unique, it is
the maximal compact invariant subset of $\,\cx$, and it  can be
characterized as
\begin{equation}\label{eqn:attrattore}
\att= \cup \{\omega(B) \ : \ \text{$\,B \subset \cx$,
bounded}\}=\omega(\cx).
\end{equation}
Besides,  if all elements of $\,\sfl \,$ are continuous functions in
 $(0,+\infty)$ and \emph{(C4)} is fulfilled, then
$\,\att\,$ is Lyapunov stable.
\end{theorem}
\begin{theorem}
\label{thm:ball2} Assume  that $\,\sfl\, $ is eventually bounded and
compact,  admits a Lyapunov function $\,V$, and that the sets of its
rest points $\,\rest\,$ is bounded. Then, $\,\sfl\,$ is also point
dissipative, and thus admits a global attractor $\,\att.$ Moreover,
$\omega(u) \subset \rest$ for all trajectories $ u \in \sfl$.
 \end{theorem}
\noindent {\bf Existence of the weak global attractor.}
\begin{theorem}
\label{teor:exist-w-att} Let $\geneset$ be a weak generalized semiflow.
 Moreover, assume that
  $\geneset$ is eventually bounded, point dissipative and compact.
  Then, $\geneset $
  possesses a unique weak global attractor $\geneatt$. Moreover,
   $\geneatt$  can be characterized as
\begin{equation}\label{struttura-attrattore}
\geneatt =\left\{\xi\in X: \mbox{ there exists a bounded complete orbit } w:
\;w(0)=\xi\right\}.
\end{equation}
\end{theorem}
\noindent Clearly, one can replace $0$ in formula
\eqref{struttura-attrattore} with any $s\in \mathbb{R}$.

 Concerning the first part
of the statement, it is sufficient to check that the argument
developed in \cite{Ball97} for the proof of Theorem \ref{thm:ball1}
goes through without the  concatenation condition (H3). As for
\eqref{struttura-attrattore}, the fact that
 the global attractor is generated by all
 complete bounded trajectories is well-known for semigroups and semiflows (cf. \cite{Temam88}),
 and, up to our knowledge, it has been observed
in some generalized framework  in \cite{Chep-Vish95,elmounir-simondon00}. Note that this characterization also holds for
the  global attractors of  the  standard generalized semiflow
constructed in Theorems \ref{thm:ball1} and \ref{thm:ball2}. As
already mentioned, we shall apply the weak global attractor
machinery to a class of differential problems for which the
existence of solutions is proved by means of an approximation
argument. In this framework, in  Section
\ref{sec:approx-attractor}
the structure formula
 \eqref{struttura-attrattore} will
play a basic role in the proof that the sequence of global
attractors of the approximate problems converges in a suitable sense
to the weak global attractor of the limit problem (see also
 \cite{Segatti06} for an analogous approximation result).

 \emph{Proof.} Arguing as in \cite[Thm. 3.3]{Ball97}, one
 has to preliminarily show the following two facts: their proof simply
 consists in repeating
  the arguments of \cite[Lemmas 3.4, 3.5]{Ball97},
  which are valid independently of (H3).
 \begin{description}
 \item[{\bf Claim 1.}] \emph{If $\geneset$ fulfills} (H1), (H2), (H4)\emph{
 and is asymptotically compact, then for any non-empty and bounded
 set $B \subset X,$
  $\omega(B)$ is non-empty, compact, quasi-invariant, and attracts
  $B$}.
  \item[{\bf Claim 2.}] \emph{If $\geneset$ fulfills} (H1), (H2),
  (H4), \emph{it
  is asymptotically compact and point dissipative,
  then there exists a bounded set $\mathcal{B}$ such that for any
  compact set $K \subset X$ there
  exist $\tau=\tau(K) >0$ and $\eps=\eps(K) >0$
  with $T(t)(N_\eps (K)) \subset \mathcal{B} $ for all $t \geq \tau(K).$}
 \end{description}

 Hence, let us define $\geneatt:= \omega (\mathcal{B})$ where ${\cal B}$
  is exactly the bounded set of Claim 2. Owing to Claim 1,
 $\geneatt$ is non-empty, compact, quasi-invariant, and attracts
 $\mathcal{B}. $ Let us now fix any bounded set $B$ and consider its
 compact  $\omega$-limit $K:=\omega (B)$, which attracts $B$ by
 Claim 1. Using Claim 2, one readily exploits \eqref{eq:inclu} and
 adapts the proof of \cite[Thm. 3.3]{Ball97} in order to  infer that
 $\mathcal{B}$ attracts $B$ as well. Thus, also $\geneatt$ attracts
 $B$ and, being $B$ arbitrary among bounded sets, we have
  checked that $\geneatt $ is the weak global attractor.

 Let us now  prove \eqref{struttura-attrattore}. To this aim, we fix
 $\xi\in\geneatt$. Then, the \emph{quasi-invariance} of $\geneatt$
 entails that there exists a complete orbit $w$ such that $w(0)=\xi$
 and $w(t)\in \geneatt$ for any $t$. In particular, $w$ is also
 bounded since $\geneatt$ is bounded and we shave shown one inclusion
 in \eqref{struttura-attrattore}. To prove
  the converse
 inclusion, consider any bounded and complete orbit $w$ in $\geneset$
 and  set $\mathcal{O}:=\left\{w(t),\;\; t\in \mathbb{R}\right\}$.
 The set $\mathcal{O}$ is clearly bounded in the phase space and
 quasi-invariant, and the following chain of inclusions holds
 \begin{equation}
 \mathcal{O}\subset T(t)\mathcal{O}\subset \omega(\mathcal{O})\subset \geneatt.
 \end{equation}
 In fact, the first inclusion is due to the quasi-invariance of
 $\mathcal{O}$, while the second one holds since the $\omega$-limit
 set of any bounded set attracts the set itself, Finally, the last
 inclusion  follows from \eqref{eqn:attrattore}. Thus, we conclude
 that, for any bounded and complete orbit $w$ of $\geneset$, $w(0)\in
 \geneatt$. which clearly implies \eqref{struttura-attrattore}. \fin

Finally,  by adapting the proof  of \cite[Thm. 5.1]{Ball97}, one may
obtain the analogue of Theorem \ref{thm:ball2} for weak global
attractors, namely
\begin{theorem}
\label{teor:Lyapun-w-attractor}
 Let $\geneset$ be an eventually bounded and compact weak generalized semiflow.
Moreover, suppose that $\geneset$ admits a  Lyapunov  function, and
that there exists a non-empty subset $\mathcal{D}$ of $X$ such that
\begin{eqnarray}
&& \label{eq:invariance-2} {T}(t) \mathcal{D} \subset \mathcal{D}
\quad \forall t \geq 0,
\\
&& \label{eq:rest-points-2} \text{the set $ Z(\geneset) \cap
\mathcal{D} $ is bounded in $\cx$}.
\end{eqnarray}
 Then, $\geneset$
possesses a unique weak global attractor $\geneatt$ in
$\mathcal{D}$. Furthermore, for any trajectory $u \in \gamma^+
(\mathcal{D})$ we have $\omega(u) \subset Z(\geneset)$ and
the weak
global attractor $\geneatt$ complies with
\eqref{struttura-attrattore}.
\end{theorem}
Indeed, Theorem \ref{teor:Lyapun-w-attractor} directly corresponds
to Theorem \ref{thm:ball2} with the choice $\, \mathcal{D}=X$. On
the other hand, the need for restricting  the natural phase space
$\, X \,$ to a proper subset $\, \mathcal{D} \,$ is well motivated
by applications and the reader is referred to Section \ref{sez:6}
for some example in this direction.

\section{Main results}
\label{sez:4}
 In view of the assumption $u_0 \in D(\phi)$
  in the existence Theorem \ref{teor:exist-gflows},
   we are naturally led to work in the phase
space
\begin{equation}
\label{eq:phase-space} \cx:=D(\phi), \quad  \text{with} \quad
\dcx(u,v):= |u-v|+ |\phi(u)-\phi(v)| \quad  \forall u,v \in \cx.
\end{equation}
Note that $(\cx, \dcx)$ is not,  in general, a \emph{complete}
metric space.

 For the sake of simplicity, we will assume that $0\in
D(\phi)$ and $\phi(0)=0$, but it is clear that this assumption is
not at all restrictive,  since with a proper translation we can deal
with the general case in which $0\notin D(\phi)$. Hence,  a subset
$B\subset \cx$ is $\dcx$-bounded iff it is contained in a
$\dcx$-ball $B(0,R)$ for some  $R>0$, i.e.
\begin{equation}
\label{e:boundedness-spazio-fasi}
 |u| + |\phi(u)| \leq R \quad
\forall u \in B.
\end{equation}
\begin{definition}
\label{def:semiflow} We  denote by $\mathcal{S}$ the set of all
functions $\, u : [0,+\infty) \rightarrow \Hilbert \,$ such that $\,u \in H^1(0,T;\Hilbert)\,$ for all $T>0$ and
\begin{equation}
\label{eq:gflow-infty} u'(t)+ \slmSbd \phi (u(t))\ni 0
    \quad \forae \ t \in (0,+\infty).
\end{equation}
\end{definition}
\begin{remark}\rm
\label{rem:sorgente-costante} We  could include a  \emph{constant}
source term $f \in \Hilbert$  in \eqref{eq:gflow-infty} by
 replacing $\phi$ with the functional
$\phi_f$ defined by $\phi_f (v):=\phi(v) -\la f,v \ra\,$ for all $\, v \in \Hilbert$.
\end{remark}
\begin{theorem}[Generalized semiflow]
\label{teor:1} Let $\phi$ comply with the assumptions
\eqref{eq:basic-prop-phi}, \eqref{eq:coerc-ass}, and
\eqref{eq:chain-rule-energy} of Theorem \ref{teor:exist-gflows}. In
addition, assume that
\begin{equation}
\label{ass:bounded-below} \exists K_1, K_2 \geq 0 : \quad \phi(u)
\geq -K_{1}|u| -K_2 \quad \forall u \in \Hilbert.
\end{equation}
 Then,
 $\semif$ is a generalized semiflow
on $\cx$, whose elements are continuous functions on $[0,+\infty)$ and comply with {\rm
(C4)}.
\end{theorem}
In order to study the long-time behaviour of our gradient flow equation, we assume an
 additional continuity property of the potential $\phi$, that is
\begin{equation}
\label{eqn:cont}
v_n\rightarrow v,
\;\;\;\sup_n\big(\vert(\partial_\ell\phi(v_n))^\circ\vert,\phi(v_n)\big)<+\infty
\;\;\Rightarrow \phi(v_n)\rightarrow\phi(v).
 \tag{\textsc{cont}}
\end{equation}
Note that \eqref{eqn:cont} is readily fulfilled by lower
semicontinuous convex functionals (cf.
\eqref{eq:continuity-convexity}). Let
  $\,\{ T(t)\}_{t \geq 0} \,$ be the family of operators
   \eqref{eq:operat-T}
associated with the generalized semiflow $\semif$.
    We have
\begin{theorem}[Global attractor]
\label{teor:2} Let $\phi$ fulfil
\eqref{eq:basic-prop-phi}, \eqref{eq:coerc-ass},
\eqref{eq:chain-rule-energy}, \eqref{eqn:cont},
 and
\begin{equation}
\label{ass:strong-coercivity} \liminf_{|u| \rightarrow
+\infty}\phi(u) =+\infty.
\end{equation} Further, let $\mathcal{D} $
be a non-empty subset of $ \cx$ satisfying
\begin{eqnarray}
&& \label{eq:invariance} {T}(t) \mathcal{D} \subset \mathcal{D}
\quad \forall t \geq 0,
\\
&& \label{eq:rest-points} \text{the set $ Z(\semif) \cap \mathcal{D}:= \{u \in D(\partial_s \phi)\,:\, 0 \in \partial_s \phi(u)\} \cap \mathcal{D}
$ is bounded in $\cx$}.
\end{eqnarray}
Then, there exists a unique attractor $\att$ for $\semif$ \emph{in
$\mathcal D$}, given by
$$
\att: = \cup \left\{ \omega(D) \, : \, D \subset \mathcal D \
\text{bounded} \right\}.
$$
Moreover, $\att $ is Lyapunov stable.
\end{theorem}

With respect to applications, let us stress that assumptions \eqref{eq:invariance}-\eqref{eq:rest-points} are of course to be checked for all current choices of the functional $\, \phi$. In order to fix ideas, let us remark that, in the convex case, \eqref{eq:rest-points} follows for instance from \eqref{ass:strong-coercivity}.
\paragraph{The fundamental theorem of Young measures
for weak topologies.} Before developing the proof of Theorems
\ref{teor:1}, \ref{teor:2}, we report a  compactness result for
Young measures
  in the  framework of the weak topology,
  which shall play a
crucial role in the sequel. Hence,  for the reader's convenience
 let us   recall the definition of {\em (time-dependent) parametrized
(or Young) measures}. Denoting by
 $\mathcal{L} $ the
$\sigma$-algebra of the Lebesgue measurable subsets of $(0,T)$ and
by
  $\mathscr B(\Hilbert)$ the
 Borel $\sigma$-algebra of $\Hilbert$,
we define
  a \emph{parametrized  (Young) measure} in  $
  \Hilbert$
  to be  a family
  $\nnu:=\{\nu_t\}_{t \in (0,T)} $ of Borel probability measures
  on $\Hilbert$
  such that for all $B \in \mathscr{B}(\Hilbert)$ the map
$
 t \in (0,T) \mapsto \nu_{t}(B) $ is
 ${\mathcal
L}\mbox{-measurable}$. We denote by $\mathcal{Y}(0,T;\Hilbert)$ the
set of all parametrized measures. The following result has been
proved in \cite{Rossi-Savare04} (cf. Thm. 3.2 therein), as a
consequence
 of
the so-called fundamental compactness
 theorem  for Young measures,
\cite[Thm. 1]{Balder84} (see  also \cite{Ball89}).
\begin{theorem}
  \label{teor:young-measures-weak}
  Let $\{ v_n\}_{n\in\N}$ be a \emph{bounded} sequence
  in $L^p(I;\Hilbert)$, for some $p>1$.
  Then,
  there exists a subsequence $k\mapsto v_{n_k}$ and
  a parametrized measure
  $ \nnu=\{ \nu_{t} \}_{t \in I} \in \mathcal{Y}(I;\Hilbert)$
  such that for a.e.\ $t\in I$
  \begin{equation}
    \label{concentration}
    \mbox{$ \nu_{t} $ is
      concentrated on
      the set
      $ L(t)$ of the \emph{weak limit points} of $ \{v_{n_k}(t) \}$,}
  \end{equation}
  \begin{equation}
    \label{basicrel2}
  \int_I\left( \int_{\Hilbert}       |\xi|^p \, d \nu_{t}(\xi) \right)
  dt\le \liminf_{k \rightarrow \infty}
  \int_I|v_{n_k}(t)|^p \,dt<+\infty.
\end{equation}
Moreover, setting $v(t):=\int_\Hilbert \xi\,d\nu_t(\xi),$ we have
\begin{equation}
  \label{eq:35}
  v_{n_k}\weakto v\text{ in }L^p(I;\Hilbert)  \ \text{if $p<\infty$}
  \quad\text{and} \ \
v_{n_k}\weaksto v\text{ in }L^{\infty}(I;\Hilbert) \ \text{if
$p=\infty$}.
\end{equation}
\end{theorem}

 Henceforth, we will   denote by $ C$
any positive constant coming into play throughout  the following
proofs, pointing out the occurring exceptions.
\subsection{Proof of Theorem \ref{teor:1}}
\label{sez:4-2}
 It follows from Theorem \ref{teor:exist-gflows}
  that for any $u_0 \in
\cx$ there exists $u:(0,+\infty) \to \Hilbert $ fulfilling
$u(0)=u_0$ and \eqref{eq:gflow-infty}. Moreover, the energy identity
$$
 \int_{0}^{t}
      |u' (\sigma)  |^{2} \, d \sigma + \phi(u (t))    =   \phi(u_0)
      \quad\forall t \in [0,+\infty)
$$
yields that $u' \in L^2(0,T;\Hilbert)\,$ for all $\, T>0 $, whence
$u \in H^1(0,T;\Hilbert)\,$ for all $\, T>0$: therefore, $\semif$
complies with { (H1)}. It is easy to check that $\semif$ satisfies {
(H2)} and { (H3)} as well. Besides, the elements of $\semif$
are continuous functions on $[0,\infty)$:
 in fact, $u \in C^0([0,T];\Hilbert)\,$ for all $\,
T >0\,$ and, in view of \eqref{eq:chain-rule-energy},  $\phi \circ u
\in AC(0,T)$ for any $T>0$.
\par\noindent{\bf Proof of (H4).}  Let
us fix a sequence $\{u_0^n \} \subset D(\phi)$ converging to $u_0
\in D(\phi)$ w.r.t. the metric of $\cx$, i.e.
\begin{equation}
\label{eq:convergenze-dati}
 |u_0^n - u_0| + |\phi(u_0^n) -\phi(
u_0)| \to 0 \quad \text{as $n \up +\infty$,}
\end{equation}
and let $u_n \in H^1(0,T;\Hilbert)\,$ for all $\, T >0\,$ be the
corresponding  sequence of solutions in $\semif$. We split the proof
of (H4) into steps.
\par\noindent{\bf A priori estimates on $\{u_n\}$.}
\emph{For any $T>0$
there exists a positive constant $C_{T}$, only
depending on $u_0$ and $T$, such that}
\begin{eqnarray}
&&
 \label{eq:aprio-1}
 \| u_n \|_{L^{\infty}(0,T;\Hilbert)} + \| u_{n}'
 \|_{L^{2}(0,T;\Hilbert)}\leq C_{T},
\\
&&
 \label{eq:aprio-2}
 \sup_{[0,T]} |\phi(u_{n}(t))| \leq C_{T}.
 \end{eqnarray}
 Indeed, it follows from the energy identity  and from
 \eqref{eq:convergenze-dati}
 that
\begin{equation}
\label{eq:enid-enne}
 \int_{0}^{t}
      | u_{n}'(\sigma)  |^{2} \, d \sigma + \phi(u_n (t))   =    \phi(u_0^n ) \leq C
 \end{equation}
for any $n \in
 N$ and $t \in (0,+\infty)$.
 On the other hand, for any fixed $T>0$ and
  $t \in (0,T],$
 \begin{eqnarray}
 &&
\frac{1}{2} \int_{0}^{t}
      | u_{n}'(\sigma)  |^{2}
       \geq \frac{1}{2t} |u_{n}(t) - u_0^n|^2
\nonumber\\
&&\geq \frac{1}{4T} |u_{n}(t)|^2 - \frac{1}{2T} |u_0^n|^2 \geq K_1
|u_{n}(t)| -TK_{1}^2  - \frac{1}{2T} |u_0^n|^2.
 \nonumber
\end{eqnarray}
 Therefore, \eqref{eq:enid-enne} yields
$$
\frac{1}{2} \int_{0}^{t}
      | u_{n}'(\sigma)  |^{2} + K_1
|u_{n}(t)| + \phi(u_n (t))  \leq C + \frac{1}{2T} |u_0^n|^2 +
TK_{1}^2,
$$
for any $t \in (0,T]$,
 whence \eqref{eq:aprio-1} (in view of \eqref{ass:bounded-below}),
 as well as \eqref{eq:aprio-2}.
\par\noindent{\bf Convergence results for $\{u_n\}$.} \emph{There exist a subsequence
 $\{u_{n_k}\}$,
a function \\
$u \in H^1(0,T;\Hilbert)\,$ for all $\, T >0$,  and a limit Young
measure $\nnu=\{\nu_t\}_{t \in (0,+\infty)} \in
\mathcal{Y}(0,+\infty;\Hilbert) $
associated with $\{u_{n_k}'\}$,
 such that}
\begin{eqnarray}
&& \label{eq:conve-1} u_{n_k} \weakto u \quad \text{in
$H^{1}(0,T;\Hilbert)$} \quad \forall T>0,
\\
&&
 \label{eq:conve-2} u_{n_k} \to u \quad \text{in
$C^{0}([0,T];\Hilbert)$} \quad \forall T>0,
\\
&& \label{eq:concentrazione-crucial} \text{$\nu_t$ is concentrated
on $-\lmSbd \phi(u(t)),$ for a.e. $ t \in (0,+\infty),$ and}
\\
&& \label{eq:basic-relazione} \int_0^T\left( \int_{\Hilbert}
|\xi|^2 \, d \nu_{t}(\xi) \right)
  dt\le \liminf_{k \rightarrow \infty}
  \int_0^T|u_{n_k}'(t)|^2 \,dt<+\infty \quad \forall T>0.
\end{eqnarray}
  The
estimates \eqref{eq:aprio-1},  \eqref{eq:aprio-2}, and the
assumption \eqref{eq:coerc-ass} yield that for any fixed
${T}>0$ there exists a compact set $\mathscr{K}({T})
\subset \Hilbert$ such that\linebreak $\,
\cup_{n \in \N} \left\{ u_n (t) \ : \ t \in [0,{T}] \right \}
\subset \mathscr{K}({T})$.
Hence, taking into account   the estimate \eqref{eq:aprio-1} and
applying the generalized Ascoli theorem \cite[Lemma 1]{Simon86}, we conclude that there exist a
subsequence $u_n$  (which we do not relabel) and a limit $u \in
H^{1}(0,{T};\Hilbert)$ fulfilling \eqref{eq:conve-1} and
\eqref{eq:conve-2}  on $(0,{T})$. On the other hand, using the
compactness result for Young measures Theorem
\ref{teor:young-measures-weak},   up to a further extraction
 we also find
 a limit Young measure  $\nnu \in
\mathcal{Y}(0,{T};\Hilbert)$   such that  the
lower-semicontinuity relation
 \eqref{basicrel2} with $p=2$,  and the concentration property
 \eqref{concentration}   hold for $\nnu$ and
$\{u_{n}'\}$. Note that \eqref{concentration} yields relation
 \eqref{eq:concentrazione-crucial} on $(0,{T}) $ for all $T>0$.
Indeed,
  the set
$L(t)$ of the weak limit points of $u_{n}'(t)$ fulfils
$$
 L(t) \subset
-\lmSbd \phi(u(t)) \quad \forae \ t \in (0,T)
$$
  in
 view of  the convergence \eqref{eq:conve-2} for $\{u_{n}(t)\}$, of
the a priori estimate \eqref{eq:aprio-2} for $\{\phi(u_{n}(t))\}$,
   and of Lemma
 \ref{lemma:closure-properties-limiting}.
Then, by a diagonal argument, we extend  the maps $t \mapsto u(t) $
and $t \mapsto \nu_t $
 to $(0,+\infty)$, finding   that $u \in H^1(0,T;\Hilbert)\,$ for all $\, T>0$,  $\nnu \in
\mathcal{Y}(0,+\infty;\Hilbert)$ and fulfils
\eqref{eq:concentrazione-crucial}, and we extract a subsequence
$u_{n_k}$ for which \eqref{eq:conve-1}, \eqref{eq:conve-2}, and
\eqref{eq:basic-relazione}  hold.
\par\noindent{\bf Passage to the limit.}
 \emph{The limit function $u $ belongs to $\semif$, fulfils $u(0)=u_0$,
and}
\begin{equation}
\label{eq:final-convergences} \dcx(u_{n_k}(t), u(t)) \to 0 \quad
\text{as $k \up +\infty$} \quad \forall t \geq 0.
\end{equation}
Let us now fix an arbitrary $t >0$: taking the $\liminf$ of both
sides of \eqref{eq:enid-enne} in view of
\eqref{eq:convergenze-dati},
\eqref{eq:conve-1}-\eqref{eq:basic-relazione},  and of the lower
semicontinuity of $\phi$, we find
\begin{equation}
\label{eq:pass-lim-1}
 \frac12 \int_{0}^{t}  |u'(s)|^2 ds +
 \frac12  \int_{0}^{t} \left(\int_{\Hilbert} |\xi|^2 \, d \nu_{s}(\xi) \right)
  ds
      + \phi(u(t)) \leq
      \phi(u_0 ).
      \end{equation}
On the other hand,
 \eqref{eq:concentrazione-crucial} and \eqref{eq:basic-relazione}
ensure that
$$
 \int_0^t |(\lmSbd\phi(u(s)))^\circ|^2\,ds<+\infty.
$$
Hence, {by \cite[Thm.~3.3, Prop.~3.4]{Rossi-Savare04}}, there
exists a
 selection $ \xi(\cdot) \in \lmSbd
\phi(u(\cdot)) $ in $L^2(0,t;\Hilbert)$. Also,
\eqref{eq:pass-lim-1} yields that $\phi \circ u$ is bounded on
$(0,T).$
 Thus, we  may apply the chain rule \eqref{eq:chain-rule-energy} and
 conclude that $\phi
\circ u \in AC(0,T). $ Owing to
\cite[Thm.~3.3]{Rossi-Savare04}, the limit Young measure $\nnu$,
fulfilling \eqref{eq:concentrazione-crucial}, also complies with  a
chain rule formula, yielding that (here we set $ \eta(s):=
\int_\Hilbert \xi d \nu_{s}(\xi)$):
\begin{equation}
\label{eq:chain-rule-1-bis} \phi(u_0)-\phi(u(t)) = \int_0^t \la
u'(s), \eta(s) \ra  ds.
\end{equation}
Note that, owing to Theorem \ref{teor:young-measures-weak} and, in particular, \eqref{eq:35}, we have
that $\, u' = \eta \,$ almost everywhere.
 Combining \eqref{eq:pass-lim-1} and
\eqref{eq:chain-rule-1-bis}, we deduce
\begin{gather}
  \frac12\int_0^t  \int_\Hilbert |\xi-u'(s)|^2 d
 \nu_{s}(\xi) \nonumber\\
= \frac12 \int_0^t|u'(s)|^2ds + \frac12 \int_0^t \int_\Hilbert
|\xi|^2 d
 \nu_{s}(\xi) - \int_0^t \left(u'(s),\int_\Hilbert \xi d
 \nu_{s}(\xi)\right)ds \leq 0,\nonumber
\end{gather}
 whence
 $
\nu_s = \delta_{u'(s)} $ $ \forae \ s \in (0,t).
 $
Therefore, \eqref{eq:concentrazione-crucial} yields
$$
u'(s) \in -\lmSbd \phi(u(s)) \quad \forae \ s \in (0,t).
$$
Being $t$ arbitrary, we infer that $u $ solves
\begin{equation}
\label{eq:cauchy-limiting}
 u'(t)+ \lmSbd \phi (u(t))\ni 0
    \quad \forae \ t \in (0,+\infty).
\end{equation}
The initial condition $u(0)=u_0$ of course ensues from
 \eqref{eq:convergenze-dati} and \eqref{eq:conve-2}. Finally, let us
take the $\limsup_{k \up +\infty}$ of \eqref{eq:enid-enne}. By
\eqref{eq:convergenze-dati},
\begin{gather}
\limsup_{k \up +\infty} \left( \int_{0}^{t}
      | u_{n_k}'(\sigma)  |^{2} \, d \sigma + \phi(u_{n_k} (t)) \right)
      \leq
      \limsup_{k \up +\infty} \phi(u_0^{n_k})
= \phi(u(0))\nonumber\\
   =   \phi(u(t)) + \int_{0}^{t}
      | u'(\sigma)  |^{2} \, d \sigma,\nonumber
\end{gather}
where we have used that $u$ fulfils \eqref{eq:cauchy-limiting} and
applied  the chain rule \eqref{eq:chain-rule-energy} to the selection
$-u'(t) \in \lmSbd \phi(u(t)).$ Therefore, we deduce that for every
$t>0$
\begin{eqnarray}
&&
 \int_{0}^{t}
      | u'(\sigma)  |^{2} \, d \sigma \leq \liminf_{k \up +\infty} \int_{0}^{t}
      | u_{n_k}'(\sigma)  |^{2} \, d \sigma \nonumber\\
&&\qquad\qquad\leq
      \limsup_{k \up +\infty} \int_{0}^{t}
      | u_{n_k}'(\sigma)  |^{2} \, d \sigma \leq \int_{0}^{t}
      | u'(\sigma)  |^{2} \, d \sigma,
      \nonumber
      \\
      &&
 \phi(u(t)) \leq \liminf_{k \up +\infty} \phi(u_{n_k} (t)) \leq
 \limsup_{k \up +\infty}  \phi(u_{n_k} (t)) \leq \phi(u(t)).
 \nonumber
\end{eqnarray}
Finally,
\begin{gather}
 u_{n_k} \to u \quad \text{strongly in
$H^1 (0,T;\Hilbert)$ for any $T>0$,} \nonumber\\
\quad u_{n_k} (t) \to u(t)
\quad \text{in $\cx$ for any $t>0,$}\label{eq:final-convergences-2}
\end{gather}
whence it is easy to infer that $u $ in fact solves
\eqref{eq:gflow-infty}. We can conclude that $\semif$ fulfils~(H4).
\par\noindent{\bf Conclusion of the proof.}
Note that we can in fact improve \eqref{eq:final-convergences-2}. By subtracting the energy identity for $u$ from the energy
identity for $u_{n_k}$, we get that
\begin{gather}
\left| \phi(u_{n_k}(t)) -\phi(u(t)) \right| \leq
|\phi(u_0^{n_k})-\phi(u_0)| + \int_0^t \left |  |u_{n_k}'(s)|^2 -
 |u'(s)|^2 \right| ds \nonumber\\
\leq |\phi(u_0^{n_k})-\phi(u_0)|
 +  \int_0^t \left(|u_{n_k}'(s))| +
|u'(s)|\right) |u_{n_k}'(s)) -u'(s)| ds\nonumber\\
 \leq
|\phi(u_0^{n_k})-\phi(u_0)|  + C \|u_{n_k}'-u'\|_{L^2
(0,t;\Hilbert)},\label{eq:conv-unif-cpt}
\end{gather}
and the above right-hand side goes to zero, as $k \up +\infty$,
uniformly in $t$ on the compact subsets of $[0,+\infty)$,
 so that we may conclude that $\,
\phi \circ u_{n_k} \to \phi \circ u \,$ {uniformly on  the
compact subsets of $[0,+\infty).$
We have thus also proved
 the continuity property { (C4)}.
 \fin
\subsection{Proof of Theorem \ref{teor:2}}
\label{sez:4-3} Note that any of the trajectories $u \in \semif$ complies
with the energy identity \eqref{enid2} (with $f \equiv 0$). In
particular, $\phi(u(t)) \leq \phi(u(s))$ for all $0 \leq s \leq t
<+\infty $. Let us check that $\phi$
is a
 Lyapunov function for $\semif$. In fact, $\phi$ is trivially
 continuous w.r.t. the topology of $\cx$. Let now $w: \R \to \cx$
 be a complete orbit for $\semif$. Note that,
 by Definition \ref{def:semiflow}, $w \in H_{\text{loc}}^1(\R;\Hilbert)
 $, and  that it fulfils the energy identity
\eqref{enid2} on $\R$.  Suppose that  the function $t \in \R \mapsto
\phi(w(t))$ is constant: hence,
$$
\int_s^t |w'(\sigma)|^2 d \sigma = \phi(w(s))-\phi(w(t))=0 \quad
\forall  \, s,t \in \R, \ s \leq t.
 $$
 Thus, $w'(t)=0$ for a.e. $t \in \R$; as  $w$ is absolutely continuous, we
 deduce that $w$ is a stationary orbit.

 Therefore, in view of Theorem \ref{thm:ball2} and of the
  assumptions \eqref{eq:invariance}-\eqref{eq:rest-points},
   it is sufficient to show  that $\semif$ is
\emph{eventually bounded} and
\emph{compact}.
\par\noindent{\bf Eventually boundedness.} In order to check that
 $\semif$ is eventually bounded, we fix  a ball  $B(0,R)$ centered at $0$
   of radius $R$ in $\cx$:
we will show that there exists $R'>0$ such that the evolution of the ball
$B(0,R)$ is contained in the ball $B(0,R')$.
  Indeed,
      let $u \in \semif$ be a trajectory
   starting from some $u_0 \in B(0,R)$, cf.
\eqref{e:boundedness-spazio-fasi}.
  By the energy identity,
\begin{equation}
\label{e:ad-event-bounded-1}
 \phi(u(t)) \leq \int_0^t |u'(s)|^2 ds +  \phi(u(t)) \leq \phi(u_0)
\leq R \quad \forall t \geq 0.
\end{equation}
Therefore, taking into account our coercivity assumption
\eqref{ass:strong-coercivity}, we deduce that
\begin{equation}
\label{e:ad-event-bounded-2}
 |u(t)| \leq R'' \quad
\forall t \geq 0,
\end{equation}
for some $\, R''>0\,$ and the eventual boundedness follows with $R' := R + R''$.

\paragraph{Compactness.} In order to  verify that $\semif$ is compact, we  consider  a
sequence $u_n\in\semif$ such that $u_n(0)$ is bounded in $\cx$: we
will show that
\begin{gather}
\text{
 there exists a subsequence $u_{n_k}$  such that}\nonumber\\
\text{
$u_{n_k}$ is convergent in $\cx$ \emph{for all t} $>0$.}\label{e:compactness} 
\end{gather}
In fact, since  $u_n(0)$ is bounded in $\cx$, we may   argue
 as in the proof
of (H4) in Theorem \ref{teor:1}, and obtain
that there exists a subsequence $u_{n_k}$ and a limit function $u \in
H^1(0,T;\Hilbert)\,$ for all $\, T>0\,$ such that the a priori bounds \eqref{eq:aprio-1}-\eqref{eq:aprio-2} and
 the convergences (\ref{eq:conve-1})-(\ref{eq:basic-relazione})
hold. However, unlike in the proof of Theorem \ref{teor:1}, we
cannot directly conclude (\ref{eq:pass-lim-1}) anymore, since now we
only have
 $u_{n}(0)\rightarrow u_{0}$ in $\Hilbert$. Actually, we will prove
 \eqref{e:compactness} by
 combining {the assumed
  (cf. (\ref{eqn:cont}))  continuity of $\phi$ along the
  sequences with equibounded slope with Helly's
compactness principle for monotone functions with respect to the
pointwise convergence (for the proof of this  result, the reader is
referred to, e.g., \cite[Chap.~4]{Ambrosio-Gigli-Savare04}}.
 Indeed,  thanks to the energy identity
\begin{equation}
      \label{eqn:energy a livello n}
      \int_{s}^{t}
      |u'_n (\sigma)  |^{2} \, d \sigma + \phi(u_n (t))  =   \phi(u_n(s))
      \quad\forall s,t \in [0,T],
    \end{equation}
  the function $t\mapsto
\phi(u_n(t))$ is non-increasing. Thus, Helly's Theorem 
applies,  and we obtain that there exists a function
$\varphi:[0,+\infty)\rightarrow (-\infty,+\infty]$, which is
non-increasing, such that
\begin{equation}\label{eqn:lim1}
\varphi(t):=\lim_{k\up +\infty}\phi_{n_k}(u(t)) \;\;\;\forall t\geq 0,
\end{equation}
for a proper subsequence $n_{k}$ of $n$. Now, by
\eqref{eq:gflow-infty} and
 (\ref{eq:aprio-1}),
we have $$
 \sup_{k\in \N} \int_0^T \left|
(\partial_\ell\phi(u_{n_k}(t)))^\circ \right|^2 \,
 dt < +\infty.
$$ Hence, by  Fatou's Lemma,
\begin{equation}\label{eqn:fatou}
\liminf_{k\up +\infty} \left|  (\partial_\ell\phi(u_{n_k}(t)))^\circ
\right|^2 <+\infty  \quad \forae \ t \in (0,T).
\end{equation}
Therefore, for almost any $t$ we can select a proper subsequence
$n_{k_{\lambda}}$ of $n_{k}$ (note that, at this stage, the latter extraction depends on $t$) such that
 $ (\partial_\ell\phi(u_{n_{k_\lambda}}(t)))^\circ$ is bounded as $\lambda\up +\infty$.
 Also in view of (\ref{eq:aprio-2})
  and (\ref{eq:conve-2}), we can now exploit
 (\ref{eqn:cont}) and conclude that
\begin{equation}\label{eqn:lim2}
\lim_{\lambda\up +\infty}\phi(u_{n_{k_{\lambda}}}(t))= \phi (u(t)).
\end{equation}
Actually, the extraction of the  subsequence in (\ref{eqn:lim2})
does not in fact depend on $t$,  since, by the lower semicontinuity
of $\phi$, we have  $\forae \ t \in (0,T)$
\begin{equation}
\label{eqn:lim3}\liminf_{k\up +\infty}\phi(u_{n_{k}}(t))\le
 \liminf_{\lambda\up+\infty}\phi (u_{n_{k_{\lambda}}}(t))
 =\phi(u(t))\le \liminf_{k\up +\infty}\phi(u_{n_{k}}(t)),
\end{equation}
yielding
\begin{equation} \label{eqn:phi-varphi}
\phi(u(t)) = \lim_{k\up +\infty}\phi(u_{n_k}(t)) =\varphi(t) \quad \forae  \ t \in (0,T).
\end{equation}
 In the next lines, we will actually show that
\eqref{eqn:phi-varphi} holds
 \emph{for all} $t >0$,
thus concluding \eqref{e:compactness} thanks to (\ref{eqn:lim1}). To
this aim, we will use the same technique  devised for proving the
upper semicontinuity property (H4). First,
 we take the $\liminf$ as $k\up+ \infty$ of both sides of
(\ref{eqn:energy a livello n}). In view of the convergences
(\ref{eq:conve-1})-(\ref{eq:basic-relazione}), (\ref{eqn:lim1}), and of
the fact that $\varphi(t)=\phi(u(t))$ for almost every $t>0$, we
obtain
\begin{equation}
\label{eqn:pass-lim-comp-1}
 \frac12 \int_{s}^{t}  |u'(\sigma)|^2 ds +
 \frac12  \int_{s}^{t} \left(\int_{\Hilbert} |\xi|^2 \, d \nu_{\sigma}(\xi) \right)ds  + \phi(u(t))
       \leq
      \phi(u(s) )
      \end{equation}
      for all $t\in (0,T]$ and for a.e. $0<s\le t$.
Now, by arguing exactly as in the proof of (H4) (with the sole
difference that all the time integrals are now considered between
$s$ and $t$, with $s>0$, since we do not have the convergence in
$\cx$ for the sequence of the initial values $u_n(0)$), we deduce
that the limit function $u$ solves
\begin{equation}
 u'(t)+ \lmSbd \phi (u(t))\ni 0
    \quad \forae \ t \in (0,+\infty).
\end{equation}
Thus, in view of (\ref{eq:chain-rule-energy}), the function $u$ also
verifies the energy identity on the interval $(s,t)$, with $\, 0< s
\leq t \leq T\,$ (see \cite[Theorem 3]{Rossi-Savare04}). In
particular, this means that the map $t \mapsto \phi(u(t))$ is
continuous and non-increasing,  and thus $\varphi(t)=\phi(u(t))$
\emph{for any} $t>0$, as desired.
 \fin


\section{Applications: perturbations of convex functionals}\label{applications}
In this section, we apply our abstract theory to some  concrete
examples of parabolic partial differential equations. More
precisely, we will deal with the long-time dynamics of gradient
flows of  various non convex perturbations of convex functionals.
First, we shall consider the case of $C^1$ perturbations. Secondly,
we apply our abstract results to gradient flows of functionals
$\phi$ given by the difference of two convex and lower
semicontinuous functionals.
\subsection{$C^1$ perturbations of convex functions}
We consider functionals $\phi:\Hilbert \to (-\infty,+\infty]$ of the type
\begin{equation}\label{appl1}
\begin{gathered}
 \phi=\phi_1+\phi_2, \ \text{with}  \
\text{$\phi_1$ proper,  l.s.c.,}\nonumber\\
\text{and convex on $D(\phi_1)\subset
\Hilbert$}, \ \phi_2\in C^1(\Hilbert).
\end{gathered}
\end{equation}
 The problem of the existence of solutions of gradient
  flow equations for functionals $\phi$ of this form
  has been addressed in \cite{Marino-Saccon-Tosques89}
   (see also
   the lectures notes \cite{Ambrosio95}
    and \cite{Ambrosio-Gigli-Savare04}).
     The uniqueness of solutions is an
     open problem owing to the possible
     non convexity of the perturbation $\phi_2$.
     Here, we prove that the set of all solutions of
\begin{equation}\label{eqn:C1perturbation}
u'(t)+\partial\phi_1(u(t))+D\phi_2(u(t))=0  \;\mbox{ for a.e. }t\in (0,+\infty)
\end{equation}
is a \emph{generalized semiflow} in the phase space $\cx=D(\phi_1)$, endowed with the metric
\begin{equation}\label{metricaappl1}
\dcx(u,v):=\vert u-v\vert+\vert\phi_1(u)-\phi_1(v)\vert \;\;\forall u,v\in \Hilbert.
\end{equation}
 Moreover, we show that this
 generalized semiflow has a unique global
  attractor in the phase space $D(\phi_1)$.
  The following proposition is a consequence
  of Theorems \ref{teor:1}
   and \ref{teor:2}, with the choice $\mathcal{D}=\Hilbert$.
\begin{proposition}[Global attractor for $C^1$-perturbations of convex functions]
\label{prop:c-1-pertub} Let $\phi:\Hilbert\to (-\infty,+\infty]$ be
as in \eqref{appl1}. Suppose that $\phi$ complies with the
assumptions \eqref{eq:coerc-ass} and \eqref{ass:strong-coercivity}.
Moreover, we assume that
\begin{eqnarray}
&\forall T>0
 \ \ \disp v\in H^{1}(0,T;\Hilbert)
 \Rightarrow D\phi_2(v)\in L^{2}(0,T;\Hilbert), \;\;\label{eqn:cond-chain}\\
&\disp  \ the \ set \;\; \left\{v\in \Hilbert :
\partial\phi_1(v)+D\phi_2(v)\ni 0\right\} \mbox{ is bounded in }
D(\phi_1)\label{eqn:cond-restpoints}.
\end{eqnarray}
Then, the set of all solutions
in $H^{1}(0,T;\Hilbert),\;\forall T>0$,
 of \eqref{eqn:C1perturbation} is a \emph{generalized semiflow} on $(D(\phi_1),\dcx)$ (see
 \eqref{metricaappl1}) and
  possesses a unique global attractor. Moreover, the attractor is Lyapunov stable.
\end{proposition}
Preliminarily, we need the following
\begin{lemma}
\label{c1perturb} Let $\phi_1: \Hilbert \to (-\infty, +\infty]$ be a
proper, lower semicontinuous functional, and let $\phi_2: \Hilbert
\to (-\infty, +\infty]$ be continuous and  G\^ateau differentiable,
with $D \phi_2: \Hilbert \to \Hilbert $ demicontinuous. Set
$\phi:\phi_1 + \phi_2.$ Then,
\begin{equation}
\label{somma-sottodifferenziali} \lmSbd \phi(u) = \lmSbd \phi_{1}
(u) +D\phi_2 (u) \quad \forall \,  u \in \Hilbert.
\end{equation}
The same conclusion holds for $\slmSbd \phi$.
\end{lemma}
For the proof of this lemma, we refer the interested reader to \cite{rossi-segatti-stefanelli-preprint}.\\

\emph{Proof of Proposition \ref{prop:c-1-pertub}.} \ First of all we
note that, since $\partial_s \phi(v)=\partial\phi_1(v)+D\phi_2(v)$
 for all $v\in \Hilbert$ thanks to Lemma \ref{c1perturb}, the set of all
 solutions of \eqref{eqn:C1perturbation} coincides with set
 of all  solutions of \eqref{eq:gflow-infty}, with $\phi=\phi_1+\phi_2$.
 In order to apply our Theorems \ref{teor:1} and \ref{teor:2},
  we only need to verify the validity of the chain
   rule \eqref{eq:chain-rule-energy}, of the
    continuity condition \eqref{eqn:cont}, and of \eqref{eq:rest-points}.
    For any given functions
     $v, \xi$  like
     in the hypothesis of \eqref{eq:chain-rule-energy},
     condition \eqref{eqn:cond-chain} and the fact that
     $\phi_2\in C^{1}(\Hilbert)$ entail the validity
      of the chain rule for $\phi_2$.
       Moreover, since $D\phi_2(v)\in L^{2}(0,T;\Hilbert)$,
        by \eqref{eqn:cond-chain} we also have
         by comparison that $\partial\phi_1(v)\ni \xi-D\phi_2(v)\in L^{2}(0,T;\Hilbert)$.
 Since the chain rule \eqref{eq:chain-rule-energy} holds for $\phi_1$ by convexity,
 we conclude that it holds as well  for $\phi$.
The continuity property \eqref{eqn:cont} easily follows
from the continuity of $\phi_2$ and from the convexity
of $\phi_1$ (see \eqref{eq:continuity-convexity}).
Finally, the condition \eqref{eqn:cond-restpoints} on
the solutions of the stationary equation ensures
the validity of \eqref{eq:rest-points}. Indeed, again thanks to Lemma \ref{c1perturb}, we have
\begin{equation}\label{eqn:c1rest-points}
Z(\mathcal{S})=\left\{v\in \Hilbert:
\partial_s\phi(v)\ni 0\right\}=\left\{v\in \Hilbert : \partial\phi_1(v)+D\phi_2(v)\ni 0\right\},
\end{equation}
and the latter set is bounded by assumption. Thus, the assertion
follows. \fin
\subsection{Dominated concave perturbations of convex functions}
In this section, we apply our results to gradient flows of
functionals $\phi$ given by
\begin{equation}\label{otani}
\begin{gathered}
\phi=\psi_1-\psi_2, \ \text{with} \ \text{$\psi_1$, $\psi_2$ proper,
l.s.c.}\nonumber\\
\text{ and convex on $D(\psi_i)\subset \Hilbert, i=1,2.$}
\end{gathered}
\end{equation}
{Of course,  $ D(\phi) = D(\psi_1)\cap D(\psi_2)$.}
 The starting point of our analysis is the following Lemma, which sheds light
 on the structure of the limiting subdifferential for a functional $\phi$
  as in \eqref{otani}, and states a sufficient condition for the validity
  of the chain rule \eqref{eq:chain-rule-energy} (its proof is
  to be found in \cite[Lemma 4.8 and Lemma 4.9]{Rossi-Savare04}).
\begin{lemma}[Subdifferential decomposition and chain rule]\label{lemma:subdiffdecomp}
Let $\phi:\Hilbert\to (-\infty,+\infty]$ fulfil \eqref{otani}, \eqref{eq:coerc-ass}, and
\begin{gather}
\forall M\ge 0, \;\;\exists \rho<1, \gamma\ge 0 \; \;
 \mbox{ such that } \sup_{\xi\in \partial\psi_2(u)}\vert \xi\vert\le \rho\vert(\partial\psi_1(u))^{\circ}\vert+\gamma\nonumber\\
 \mbox{ for every } u\in D(\partial\psi_1)\mbox{ with } \max(\phi(u),\vert u\vert)\le M.\label{eqn:cond-chain2}
\end{gather}
Then,  every $g\in \partial \phi(u)$ with $\max(\phi(u),\vert
u\vert)\le M$ can be decomposed as
\begin{eqnarray}\label{eqn:decomposition}
g=\lambda^1-\lambda^2, \;\;\;\lambda^i\in \partial\psi_i(u),
\end{eqnarray}
where $\rho,\gamma$ are given in terms of $M$ by
\eqref{eqn:cond-chain2}; moreover, $\phi$ satisfies the chain rule
\eqref{eq:chain-rule-energy}.
\end{lemma}
As a consequence of this lemma,  we have that the solutions of the
gradient flow equation \eqref{eq:gflow-infty} with
$\phi=\psi_1-\psi_2$ (whenever they exist) indeed solve the equation
\begin{equation}\label{eqn:valero}
u'(t)+\partial\psi_1(u(t))-\partial\psi_2(u(t))\ni 0 \; \mbox{ for a.e. } t\in (0,+\infty).
\end{equation}
We note that  the existence of a  global attractor for the solutions
of equations of the form \eqref{eqn:valero} has been addressed by
\textsc{Valero} \cite{Valero01}. The approach in \cite{Valero01} is
different from the present one since it is based on the abstract
theory, developed by \textsc{Melnik \& Valero} in
\cite{Mell-Vall98}, of attractors for multivalued semiflows. Our
result is however sharper although less general. On the one hand, we
do not focus on the whole class of solutions of equation
\eqref{eqn:valero}, but rather on the set of solutions which can be
obtained from the gradient flow approach. Secondly, \textsc{Valero}
\cite{Valero01} tackles the problem in the phase space
$\overline{D(\psi_1)}$, endowed with the metric of $\Hilbert$. On
the other hand,  our analysis is performed in the phase space given
by the domain of the potential $\phi$, that is $\,D(\phi)$, endowed
with metric $\dcx$ (see \eqref{eq:phase-space}), which is stronger
than that of $\Hilbert$. In the following concrete examples of PDEs,
taken from \cite{Valero01}, the difference between Valero's results
and the present ones will be clarified.

 Henceforth, we fix the Hilbert space $\,\Hilbert\,$ to be $\,\Hilbert:= L^{2}(\Omega),$
$\,\Omega\,$ being a bounded domain of $\,\mathbb{R}^d\,$ with
smooth boundary $\,\partial\Omega$;
 we shall denote by $| \cdot|$ the norm in $L^2 (\Omega)$.
 In particular,  all
subdifferentials are computed with respect to the metric in $\,
L^2(\Omega)$. Let $\,x=(x_1,\ldots,x_d)\,$ represent the variable in
$\,\Omega$.
\paragraph{\bf Example 1.}
We shall consider gradient flow solutions of the equation
\begin{equation}\label{esvalero1}
\begin{gathered}
u'(t)-\Delta_p u(t) -\vert u(t)\vert^\alpha u(t)=0 \;\;\mbox{ for a.e. } t\in (0,+\infty),\\
\text{where} \;\;\;\Delta_p
u:=\sum_{i=1}^{d}\frac{\partial}{\partial x_i} \bigg(\bigg\vert
\frac{\partial u}{\partial x_i}\bigg\vert^{p-2}\frac{\partial
u}{\partial x_i}\bigg),
\end{gathered}
\end{equation}
and $\,p>2\,$ and $\,\alpha\,>0$ fulfil
\begin{eqnarray}
\left\{
\begin{array}{ll}
 2+\alpha<p, \\
 2(1+\alpha)\le\frac{dp}{d-p}, \mbox{ if }d>p.
\end{array}
\right.\label{eqn:condizionialpha}
\end{eqnarray}
Then, let us consider the following functionals defined in
$\,L^{2}(\Omega)$:
\begin{eqnarray}
\psi_{1}(v)&:=&
\disp \frac{1}{p}\sum_{i=1}^{d}\int_{\Omega}\bigg\vert \frac{\partial v}{\partial x_{i}}\bigg\vert^{p}dx \quad\mbox{ if } v\in W^{1,p}_{0}(\Omega), \nonumber\\
&& \qquad \text{and} \ \ \psi_{1}(v):=
+\infty\quad  \mbox{otherwise},\label{eqn:psi1}\\
\psi_2(v)&:=&
\disp \frac{1}{\alpha+2}\int_{\Omega}\vert v\vert^{\alpha +2}dx \quad \mbox{ if } v\in L^{\alpha +2}(\Omega), \nonumber\\ && \qquad \text{and} \ \ \psi_2(v):= +\infty\quad \mbox{otherwise}.\label{eqn:psi2}
\end{eqnarray}
It is clear that, with these choices of $\,\psi_1\,$ and $\,\psi_2$,
equation \eqref{esvalero1} can be rewritten in the form of
\eqref{eqn:valero}. Consequently, we let $\phi=\psi_1-\psi_2$,  and
we study the long-time behaviour of the solutions of the gradient
flow for $\phi$ in the framework of the phase space
\begin{eqnarray}\label{otani2}
\left\{
\begin{array}{ll}
\cx:=D(\phi)=W^{1,p}_{0}(\Omega)=\mathcal{D} \mbox{ (see } \eqref{eqn:condizionialpha}),\mbox{ with } \\
\disp\dcx(u,v):=\|u-v\|_{L^{2}(\Omega)}\\
+\bigg\vert\disp\frac{1}{p} \big(\|\nabla u\|_p^{p}-\|\nabla
v\|_{p}^{p}\big) -\disp\frac{1}{2+\alpha} \big(\|
u\|^{2+\alpha}_{2+\alpha}-\|
v\|^{2+\alpha}_{2+\alpha}\big)\bigg\vert,
\ \   u, v\in W^{1,p}_{0}(\Omega).   \\
\end{array}
\right.
\end{eqnarray}
Thanks to \eqref{eqn:condizionialpha}, it  is not difficult to prove
that $\phi$ is lower semicontinuous, has compact sublevels in
$\,L^{2}(\Omega)$, and satisfies the coercivity condition
\eqref{ass:strong-coercivity} (which clearly entails
\eqref{ass:bounded-below}). Concerning the lower semicontinuity, we
have to prove that given a sequence
$\{u_{n}\}_{n=1}^{+\infty}\subset W^{1,p}_{0}(\Omega)$ such that
$u_n\rightarrow u$ in $L^{2}(\Omega)$, one has $\liminf_{n\to
+\infty}\phi(u_n)\ge\phi(u)$. Without loss of generality, we may
suppose that $\sup_n\phi(u_n)<+\infty$. Hence, the first relation in
\eqref{eqn:condizionialpha} gives the following chain of
inequalities for a positive constant $C$ independent of $n$
\begin{eqnarray}
&\disp C\ge \phi(u_n)\ge \,\|u_n\|^{p}_{W^{1,p}_{0}(\Omega)}-
C\|u_n\|^{2+\alpha}_{W^{1,p}_{0}(\Omega)} &\disp \ge
C\|u_n\|^{p}_{W^{1,p}_{0}(\Omega)}-C,\label{semicontval1}
\end{eqnarray}
where we have also used the Young inequality. Thus, by standard weak
compactness results and  possibly extracting some not relabeled
subsequence, $u_{n}\weakto u$ in $W^{1,p}_{0}(\Omega)$ and,
 by the compact embedding $W^{1,p}_{0}(\Omega)\subset L^{2+\alpha}(\Omega)$
 (see \eqref{eqn:condizionialpha}),
  we have $u_n\rightarrow u$ strongly in $L^{2+\alpha}(\Omega)$.
  The lower semicontinuity of $\phi$ is now an
  easy consequence of this strong convergence and
  of the lower semicontinuity of the norms w.r.t.
  the weak convergence. The estimate \eqref{semicontval1}
   also shows that the sublevels of $\phi$, being bounded in $W^{1,p}_{0}(\Omega)$,
    are compact in $L^{2}(\Omega)$.
    The coercivity condition \eqref{ass:strong-coercivity}
     easily follows from \eqref{semicontval1} by noting that, since $p>2$ by assumption, we have
$$\disp \phi(v)\ge C\|u_n\|^{p}_{W^{1,p}_{0}(\Omega)}-C\disp \ge C| u|
-C \;\;\forall \;v\in W^{1,p}_{0}(\Omega).$$
In order to apply Theorems \ref{teor:1} and \ref{teor:2}
 and find that the set of all the solutions of
\begin{equation}\label{eqn:nostravalero}
u'(t)+\partial_s(\psi_1-\psi_2)(u(t))\ni 0 \mbox{ in } \Hilbert
\quad \mbox{ for a.e. } t\in (0,+\infty)
\end{equation}
generates a \emph{generalized semiflow} on $(\cx,d_{\cx})$ (see
\eqref{otani2}), possessing a unique global attractor,  we still
need to check that $\phi=\psi_1-\psi_2$ complies with the chain rule
\eqref{eq:chain-rule-energy} and with \eqref{eq:rest-points} (in
fact, \eqref{eq:invariance}  is valid since we have chosen
$\mathcal{D}=X$). As for  proving
 the validity of the chain rule for
$\phi$, we have to check that the subdifferentials of $\psi_1$ and
$\psi_2$ comply with  condition \eqref{eqn:cond-chain2}. To this
aim, we note that for any $u\in D(\partial\psi_1)$ with
$\max(\phi(u),\vert u\vert)\le M$, we have that
$\|u\|_{W^{1,p}_0(\Omega)}\le \gamma$, where $\gamma$ is a positive
constant depending on $M$ and on $\Omega$. Thus,
\eqref{eqn:cond-chain2} follows with any choice of $\rho\in (0,1)$
by simply noting that
$\vert\partial\psi_2(u)\vert=\|u\|^{\alpha+1}_{2(\alpha+1)}$, and
recalling that $W^{1,p}_{0}(\Omega)\subset L^{2(\alpha+1)}(\Omega)$
with continuous injection (see \eqref{eqn:condizionialpha}). In
order to prove \eqref{eq:rest-points} for $\phi$, we have to check
(see \eqref{otani2}) that  the set
\begin{eqnarray}\label{eqn:rest-points-valero1}
\left\{v\in \Hilbert: \partial_s(\psi_1-\psi_2)(v)\ni 0\right\}
\mbox{ is bounded in } (W^{1,p}_{0}(\Omega),d_{\cx}).
\end{eqnarray}
Note that Lemma \ref{lemma:subdiffdecomp} and the definition of $\psi_1$ and $\psi_2\,$ entail that
\begin{gather}
\disp\left\{v\in \Hilbert: \partial_s(\psi_1-\psi_2)(v)\ni 0\right\}\nonumber\\
\subseteq \left\{v\in W^{1,p}_{0}(\Omega):
 \Delta_p v-\vert v\vert^{\alpha}v=0 \
  \mbox{ a.e. in }\Omega\right\}.\label{eqn:rest-points-valero2}
\end{gather}
Then, we only need to prove that the latter set is bounded in
$(W^{1,p}_{0}(\Omega),d_{\cx})$.  This follows by simply testing in
$L^{2}(\Omega)$ the equation $\Delta_p v-\vert v\vert^{\alpha}v=0$
with $v$ and performing the same computations as
 for proving
\eqref{semicontval1}. This produces a bound in $W^{1,p}_{0}(\Omega)$
for the solutions of the aforementioned stationary equation, which
entails the bound in the phase space \eqref{otani2} by using the
embedding $ W^{1,p}_0(\Omega)\subset L^{\alpha+2}(\Omega)\,$ again.
We have thus proved the following.
\begin{proposition}\label{cor:condompert1}
Let  $\alpha$ and $p$ satisfy \eqref{eqn:condizionialpha} and
$\phi=\psi_1 - \psi_2$ with $\psi_1$ and $\psi_2$ as in
\eqref{eqn:psi1}-\eqref{eqn:psi2}.  Then, the solutions of the
gradient flow equation
\begin{equation}\label{eqn:nostrovalero11}
u'(t)+\partial_s(\psi_1-\psi_2)(u(t))\ni 0\quad  \mbox{ for a.e. }t
\mbox{ in } (0,+\infty)
\end{equation}
generate a \emph{generalized semiflow} in
$(W^{1,p}_{0}(\Omega),d_{\cx})$  which possesses a unique global
attractor. This attractor is also Lyapunov stable.
\end{proposition}
\noindent{\bf Example 2.} \ In this example, still taken from
\cite{Valero01}, we consider gradient flow solutions of
\begin{equation}\label{esvalero2}
u'(t)-\Delta u(t)-f(u(t))\in \lambda H(u(t)-1) \;\;\mbox{ for a.e. } t\in (0,+\infty),
\end{equation}
 where $H$  is the Heaviside graph,
  i.e. the maximal multivalued monotone
   graph in $ \mathbb{R}\times\mathbb{R}$  given by
\begin{equation}
H(v):= 1  \mbox{ if } v>0, \quad H(v):= [0,1]  \mbox{ if } v=0, \ \ \ \text{and} \ \ H(v)= 0  \mbox{ if } v<0,
\end{equation}
and  $\lambda$ is a non-negative constant. Finally, $f:\mathbb{R}\to\mathbb{R}$ is a  non-decreasing continuous function such that
\begin{equation}\label{eqn:f}
\vert f(s)\vert\le k_1+k_2\vert s\vert, \mbox{ with } k_1\ge 0  \mbox{ and } 0\le k_2<\lambda_1,
\end{equation}
with $\lambda_1$ the first eigenvalue of the Laplace operator with Dirichlet boundary conditions.
We introduce the following  functionals, defined in $L^{2}(\Omega)$,
\begin{eqnarray}\label{eqn:psi12}
\psi_{1}(v)&:=&
\disp \frac{1}{2}\int_{\Omega}\vert \nabla v\vert^{2}dx \quad \mbox{ if } v\in H^{1}_{0}(\Omega),\quad \psi_{1}(v):=
+\infty \mbox{ otherwise},\\
\label{eqn:psi22}
\psi_2(v)&:=&\int_{\Omega}F(u)dx+\lambda\int_{\Omega}( u-1)^+ dx,
\end{eqnarray}
where $F'=f$. As in Example $1$,  we aim to consider the dynamics of
the gradient flow for the functional $\phi=\psi_1-\psi_2$ in the
phase space
\begin{eqnarray}\label{otani3}
\left\{
\begin{array}{ll}
\disp\disp\cx:=D(\phi)=H^{1}_{0}(\Omega)=\mathcal{D}, \mbox{ with } \\
\disp\dcx(u,v):=\|u-v\|_{L^{2}(\Omega)}\\+\bigg\vert\disp\frac{1}{2}\big(\|\nabla
u\|_2^{2}-\|\nabla v\|_{2}^{2}\big)
\nonumber\\
-\disp\bigg(\int_\Omega (F(u)-F(v))dx +\lambda\int_\Omega ( u-1)^+ dx
-\lambda\int_\Omega ( v-1)^+ dx\bigg) \bigg\vert,\nonumber\\
  u, v\in H^{1}_{0}(\Omega).\nonumber
\end{array}
\right.
\end{eqnarray}
The functional $\phi$ is lower semicontinuous in $L^{2}(\Omega)$.
 In fact, if we are given a sequence
  $\left\{u_n\right\}_{n=1}^{+\infty}\subset H^{1}_{0}(\Omega)$
  with $u_n\rightarrow u$ in $L^{2}(\Omega)$ and $\disp\sup_n\phi(u_n)<+\infty$,
   then the growth condition on $f$ entails that
\begin{equation}
F(u_n)\le C +C\vert u_n\vert^{2} \ \ \mbox{ a.e. in } \Omega.
\end{equation}
Thus, by a variant of the Dominated  Convergence Theorem (see, e.g.,
\cite[Theorem 4]{evans-gariepy}), there holds $F(u_n)\rightarrow
F(u)$ in $L^{1}(\Omega)$. The lower semicontinuity of $\phi$ now
descends  from the lower semicontinuity of norms w.r.t. the weak
convergence. Moreover, $\phi$ has compact sublevels in
$L^{2}(\Omega)$. In fact, by the Poincar\'e inequality, combined
with the growth condition on $f$ (recall that $k_2<\lambda_1$), we
have that
\begin{eqnarray}\label{semicontval2}
& \phi(u)\ge\disp\frac{1}{2}\|\nabla v\|^{2}_2-\int_{\Omega}F(v)dx-\lambda\int_{\Omega}(v-1)^{+}dx\nonumber\\
&\disp \ge  C\|\nabla v\|^{2}_2-C \ \  \mbox{ for a given } C>0 \mbox{ and } \forall v\in H^{1}_{0}(\Omega).
\end{eqnarray}
Thus, the sublevels of $\phi$ are bounded in $H^{1}_{0}(\Omega)$,
 which is clearly compact
 in $L^{2}(\Omega)$. Note that \eqref{semicontval2}
  entails the coercivity assumption \eqref{ass:strong-coercivity}
   (and thus \eqref{ass:bounded-below}).
   Again, in order to apply our Theorems \ref{teor:1} and \ref{teor:2},
   we need to check \eqref{eq:chain-rule-energy} and
\eqref{eq:rest-points} (\eqref{eq:invariance} is again trivial since
we take $\mathcal{D}=X$). The chain rule
\eqref{eq:chain-rule-energy} easily follows from Lemma
\ref{lemma:subdiffdecomp}. In fact, the subdifferential
$\partial\psi_2$ of $\psi_2$ is simply
$$\partial\psi_2(v)=f(v)+\lambda H(v-1)
 \quad \forall \, v \in D(\partial \psi_2)=L^{2}(\Omega),$$
thanks to the growth condition  \eqref{eqn:f}, while the
subdifferential of $\psi_1$ is clearly $-\Delta $, with domain
$H^{2}(\Omega)\cap H^{1}_{0}(\Omega)$. Thus, condition
\eqref{eqn:cond-chain2} easily follows with $\gamma=M$ and with any
choice of $\rho\in (0,1)$. Finally, the condition on the rest points
follows from the same argument used in
\eqref{eqn:rest-points-valero1}-\eqref{eqn:rest-points-valero2}. In
this case, the analogue of the stationary equation in
\eqref{eqn:rest-points-valero2} is the following
\begin{eqnarray}\label{eqn:rest-points-valero3}
\disp   -\Delta v -f(v)\in \lambda H(v-1)  \ \ \mbox{ a.e. in } \Omega, \quad
\disp v=0  \ \ \mbox{ a.e. on  } \partial \Omega.
\end{eqnarray}
We thus have
\begin{proposition}\label{cor:condompert2}
 Let us consider
  the functional $\phi=\psi_1 - \psi_2$,
   with $\psi_1$ and $\psi_2$ as in \eqref{eqn:psi12}-\eqref{eqn:psi22}.
 Then, the solutions of the gradient flow equation
\begin{equation}\label{eqn:nostrovalero12}
u'(t)+\partial_s(\psi_1-\psi_2)(u(t))\ni 0 \quad \mbox{ for a.e. }t
\mbox{ in } (0,+\infty)
\end{equation}
generate a  \emph{generalized semiflow} in $(H^{1}_0(\Omega),d_\cx)$
which possesses a unique global attractor. The attractor is also
Lyapunov
 stable.
\end{proposition}


\noindent{\bf Example 3.}\ We are interested in the study of  the
long-time behaviour for the gradient flow solutions of the following
equation
\begin{equation}\label{esvalero3}
u'(t)-\Delta u(t)+\partial I_{K}(u(t))+f_1(u(t))-f_2(u(t))\ni 0
\quad \mbox{ for a.e. t in } (0,+\infty).
\end{equation}
In the latter equation, the symbol $\partial I_{K}$ represents the subdifferential of the indicator function of the closed and convex set $K$ (see definition \eqref{eqn:ostacolo} below), while $f_i:\mathbb{R}\to \mathbb{R}$ (i=$1,2$) are two non-decreasing continuous functions which satisfy the following growth and compatibility conditions
\begin{equation}\label{eqn:f1-f2}
\begin{gathered}
\mbox{ there exist }  \ \ 0\leq  k_1, k_2, k_3<\lambda_1, \ k_4\ge 0, \;\; \varepsilon>0  \ \
\mbox{ such that }\\
\vert f_1(s)\vert \le k_1(\vert s\vert ^{{d}/({d-2})}+1)  \quad \forall s\in \R \quad \mbox{ if } d\ge 3\\
\vert f_2(s)\vert\le k_2+k_3\vert s\vert,\\
 (f_1(s)-f_2(s))s\ge (-\lambda_1-\varepsilon)s^2-k_4,
\end{gathered}
\end{equation}
where $\lambda_1$  is the first eigenvalue of the Laplace operator
with Dirichlet boundary conditions. We then denote by $F_1$ and
$F_2$ the primitives of $f_1$  and $f_2$ respectively. Consequently,
$F_1$ and $F_2$ are differentiable convex functions in $\mathbb{R}$
such that $ F_i'=f_i$, $i=1,2$. Without loss of generality, we
assume that $F_1(0)=0$. We will only consider  the case in which
$\Omega$ is a bounded domain of $\mathbb{R}^{d}$ with $d\ge 3$, the
one-dimensional and the two-dimensional cases being easier. The only
difference is in the growth condition imposed on $f_1$,  which may
be weakened. More precisely, in two dimensions we can deal with a
function $f_1$ growing at most like a polynomial with order $\nu$,
$\nu$ being any real number $1\le \nu<+\infty$. In one dimension, we
do not need any additional growth condition.

\noindent
Now, let us consider the following functionals on $L^{2}(\Omega)$
\begin{eqnarray}\label{eqn:psi13}
\psi_{1}(v)&:=&
\disp \frac{1}{2}\int_{\Omega}\vert \nabla v\vert^{2}dx +\int_{\Omega}F_1(v)dx \quad \mbox{ if } v\in K, \nonumber\\
&&\qquad \text{and} \ \ \psi_1(v):=
+\infty \mbox{ otherwise},\\
\label{eqn:psi23}
\psi_2(v)&:=&\int_{\Omega}F_2(v)dx,
\end{eqnarray}
where $K$ is the following closed and convex subset
\begin{equation}\label{eqn:ostacolo}
K:=\left\{v\in H^{1}_{0}(\Omega): v(x)\ge 0, \mbox{ a.e. in }\Omega\right\}.
\end{equation}
It is not difficult to show (see \cite[Prop. 2.17]{Brezis73}) that the subdifferential of $\psi_1$ in $L^{2}(\Omega)$ has the following expression
\begin{equation}\label{structsubdiff}
\begin{gathered}
w^{}\in \partial\psi_1(u)\Leftrightarrow  w^{}\in -\Delta u +\partial I_{K}(u) +f_1(u),\\
\text{with} \ \ \ D(\partial\psi_1)=H^{2}(\Omega)\cap K.
\end{gathered}
\end{equation}
 Thus, it is clear that \eqref{esvalero3}
  could be rewritten
  in the form of \eqref{eqn:valero},
  with $\psi_1$ and $\psi_2$ as in \eqref{eqn:psi13}-\eqref{eqn:psi23}.

Again, we are interested  in the long-time dynamics of the gradient
flow for the functional $\phi=\psi_1-\psi_2$, in the framework of
the phase space
\begin{eqnarray}\label{otani4}
\left\{
\begin{array}{ll}
\disp\cx:=D(\phi)= K=\mathcal{D} \mbox{ with } \\
\disp\dcx(u,v):=\| u-v\|_{L^{2}(\Omega)}\\
+\bigg\vert\disp \frac{1}{2}\big(\| \nabla u\|_2^{2}-\| \nabla
v\|_{2}^{2}\big)\nonumber \disp \nonumber\\
+\bigg(\disp\int_\Omega (F_1(u)-F_1(v))dx
-\int_\Omega (F_2(u)-F_2(v))dx\bigg) \bigg\vert.
\end{array}
\right.
\end{eqnarray}
Thus, we have to check  the validity of the hypotheses of Theorems
\ref{teor:1} and \ref{teor:2}. The proof of the lower semicontinuity
follows exactly the same lines of Example 2, with some minor
modifications due to presence of the two extra terms $I_K(v)$ and
$\int_\Omega F_1(v)dx$ in the definition of \eqref{eqn:psi13}.
However, these two terms, being positive, could be easily handled by
lower semicontinuity.  Moreover,
 arguing as in \eqref{semicontval2}
we find that $\phi$ complies with \eqref{eq:coerc-ass} and
\eqref{ass:strong-coercivity} with respect to the norm of
$L^{2}(\Omega)$. The chain rule property
\eqref{eq:chain-rule-energy} is an easy consequence of Lemma
\ref{lemma:subdiffdecomp}, of  the structure of the subdifferential
of $\psi_1$ (see \eqref{structsubdiff}), and of the growth condition
on the perturbation $f_2$, which entails that the subdifferential
$\partial\psi_2$ of $\psi_2$ is simply given by
$\partial\psi_2(v)=f_2(v)\,$ for all $\, v \in L^2(\Omega)$. Thus,
in order to check the validity of \eqref{eqn:cond-chain2} and of
\eqref{eq:rest-points}, we proceed as in the former examples. In
this case, the analogue of the stationary equation in
\eqref{eqn:rest-points-valero2} reads
\begin{equation}\label{statio3}
-\Delta v+\partial I_K(v)+f_1(v)-f_2(v)\ni 0 \quad \mbox{ a.e. in
}\Omega.
\end{equation}
Hence, by simply testing  \eqref{statio3} with $v$ (which belongs to
$K$, being a solution of \eqref{statio3}) and using the last
condition in \eqref{eqn:f1-f2}, we get a bound for $v$ in
$H^{1}_{0}(\Omega)$. Thus, by the growth condition \eqref{eqn:f1-f2}
on $f_1$
 we have a similar bound for the growth of $F_1$. More
precisely, there holds
\begin{equation}
\label{cresciF}
\vert F_1(v)\vert\le C(\vert v\vert +\vert v\vert ^{{(2d-2)}/{(d-2)}}) \ \  \mbox{ for some}  \  \ C>0.
\end{equation}
Thus, since $H^{1}_{0}(\Omega)$ is continuously embedded in $L^{{(2d-2)}/({d-2})}(\Omega)$, by \eqref{cresciF} $F_1$ maps bounded sets in $H^{1}_{0}(\Omega)$ to bounded sets of $L^{1}(\Omega)$. To conclude that the set of the rest points is bounded, it remains to show the boundedness of the solutions of the stationary inclusion \eqref{statio3}, also with respect to the $F_2$-part of the metric $\dcx$ (see \eqref{otani4}). But this is simpler, thanks to the linear growth of $f_2$, which entails a quadratic growth for its primitive $F_2$.

We thus have the following.
\begin{proposition}\label{cor:condompert3}
Let $K$ be as in \eqref{eqn:ostacolo} and the functional $\phi$ be given by $\phi=\psi_1-\psi_2\,$ with $\psi_1$ and $\psi_2$ as in \eqref{eqn:psi13}-\eqref{eqn:psi23}.
 Then, the solutions of the gradient flow equation
\begin{equation}\label{eqn:nostrovalero13}
u'(t)+\partial_s(\psi_1-\psi_2)(u(t))\ni 0  \ \ \mbox{ for a.e. }t
\mbox{ in } (0,+\infty)
\end{equation}
generate a \emph{generalized semiflow} in $K$  which possesses a
unique global attractor. The attractor is also Lyapunov stable.
\end{proposition}

Note that, as we have already mentioned,  our choice of the phase
space brings to a more regular attractor, attracting with respect to
the $W^{1,p}$-norm (Example $1$) and to the $H^{1}_{0}$-norm
(Examples $2$-$3$), whereas in \cite{Valero01} the attraction holds
with respect to the $L^{2}$-metric. Furthermore, our phase space
keeps track of the constraint imposed on the
 unknowns.

\section{Applications: long-time behaviour of quasi-stationary evolution systems}
\label{sez:3-qstat-evol} \noindent
{\bf General Setup.} The
functional setting we  deal with  features
 a standard Hilbert triplet
\begin{equation}
\label{eq:hilb-triplet}
  V\subset H\equiv H'\subset V',\quad
  \text{with dense and compact inclusions.}
\end{equation}
We denote by $ \pairing{V'}{V}{\variabile}{\variabile} $ the duality
pairing between $ V' $ and $V$ and by $ (\variabile,\variabile)_{H}
$ the scalar product in $H$, recalling that $\,\pairing{V'}{V}{u}{v}  =( u,v)_H\quad\forall\, u\in H,v\in V$.
Furthermore, let $a:V\times V\to\R$ be a \emph{non-negative,
symmetric, and continuous} bilinear form, and let $A:V\to V'$ be the
continuous linear operator associated with $a$, i.e.
\begin{equation}
\label{eq:operator-A}
 \pairing{V'}{V}{Au}{v}:=a(u,v)\quad\forall\,
u,v,\in V.
\end{equation}
We also
 consider a
proper functional $F:H\to [0,+\infty]$ whose sublevels
\begin{equation}
  \label{eq:gflows_RosSav_2004:39}
  \big\{\chi\in H:F(\chi)\le s\big\}
  \quad\text{are \emph{strongly} compact in }H,
\end{equation}
and we denote by
 by $\FrSbd F:H\to 2^H$ the Fr\'echet  subdifferential of $F$ in $H$, namely
\begin{gather}
  \theta\in \FrSbd F(\chi)\quad\nonumber\\
  \Leftrightarrow\quad
  \chi\in D(F)\subset H,\quad
  \liminf_{\|\eta-\chi\|_H\to0}
  \frac{F(\eta)-F(\chi)-(\theta,\eta-\chi)_H}{\|\eta-\chi\|_H}\ge0.\nonumber
\end{gather}

We aim to investigate the long-time behaviour of (a class of
solutions of)  the following evolution system, coupling a diffusion
equation with a quasi-stationary condition, for which an existence
result  was obtained in \cite[Sect. 5]{Rossi-Savare04}.
\begin{problem}
  \label{quasi-stationary}
  Given $T >0$   and $u_0\in H$,
  find a pair  $u,\chi:(0,T)\to H$,
  with $u(t)-\chi(t)\in V$ $\forae \ t \in (0,T), $
  which satisfies at a.e.\ $t\in (0,T)$ the system
  \begin{equation}
  \label{eq:qstat-evolution}
  \left\{
    \begin{aligned}
      u'(t)+A(u(t)-\chi(t))&=0\quad\text{in }V',
      \\
      \chi(t)+\FrSbd F(\chi(t))&\ni u(t)\quad\text{in $H$},\\
      u(0)&=u_0.
    \end{aligned}
  \right.
\end{equation}
\end{problem}

\noindent
 In Section \ref{subsez:coercive}, we briefly summarize for the reader's convenience
 the
 techniques developed  in \cite[Sec.~5]{Rossi-Savare04} for
 Problem
\ref{quasi-stationary}. Hence,  we distinguish  the two following
 cases:
\begin{enumerate}
\item  the form $a$ is \emph{coercive,}
  i.e., \ there exists a constant  $\alpha>0$ such that
  \begin{equation}
    \label{eq:coercive-case}
    a(u,u)\ge\alpha \|u\|_V^2\quad\forall\, u\in V,
  \end{equation}
\item
$a$ is \emph{weakly coercive,}
  namely, \ there exist  $\lambda,\,\alpha_{\lambda}>0$ s.t.
\begin{equation}
    \label{eq:weakly-coercive-case}
    a(u,u)+\lambda \| u\|_H^2 \ge\alpha_{\lambda} \|u\|_V^2\quad\forall\, u\in
    V.
  \end{equation}
\end{enumerate}
  In fact, whenever $\, a \,$ is weakly coercive, \emph{for all} $\lambda>0$ it is possible to find a constant
   $\alpha_{\lambda}>0$ fulfilling \eqref{eq:weakly-coercive-case}.

\subsection{Existence results for Problem~\ref{quasi-stationary}}
 \label{subsez:coercive}
\subsubsection{The coercive case: existence  by a gradient flow
approach}
\label{s:5.1.1}
\noindent
 Assume that \eqref{eq:coercive-case} holds. Then, we
endow  $V$ with the norm $\| v \|_V^2:= a(v,v)$   for all $ v \in
V$  and $A$ turns out to be an isometry between the spaces $V$
and $V'$. Let us introduce the functional $\phi: V' \to \R \cup
\{+\infty\}$ defined by
\begin{align}
  \label{eq:right-functional}
  \phi(u):=&
  \inf_{\chi\in H}  \mathscr F(u,\chi),
  \quad
  \mathscr F(u,\chi): = \begin{cases}
    \frac12\|u-\chi\|_H^2+F(\chi)&\text{if }u,\chi\in H,\\
    +\infty&\text{otherwise.}
  \end{cases}
\end{align}
Clearly,  $D(\phi)=H$; further, for $u \in D(\phi) $ we denote by $
M(u)$  the set of the elements $\chi \in H$ attaining the minimum in
\eqref{eq:right-functional}, i.e.
\begin{equation}
\label{eq:argmin}
 M(u):=\big\{\chi\in H:\mathscr F(u,\chi)=\phi(u)\big\}.
\end{equation}
Note that $M(u) \neq \emptyset$ for all $u \in H$,  since $\, F \,$
is l.s.c. and has compact sublevels. {Further, the following
formula
\begin{equation}
\label{decomp-phi}
  \phi(u)=\frac12\|u\|_H^2-
  \sup_{\chi\in H}\left((u,\chi)_H-
  \Big(\frac12\|\chi\|_H^2+F(\chi)\Big)\right)
\end{equation}
shows that $\phi$ is  in fact  a concave perturbation of a quadratic
functional
 (cf. with
\eqref{eq:decomposition}).}

Proposition \ref{prop:subdif-repre} and Corollary \ref{coro:subdif}
below
 (which we recall from  \cite[Sec.5]{Rossi-Savare04})
  ensure
that
  Problem
  \ref{quasi-stationary} may be interpreted as the Cauchy problem
\begin{equation}
\label{eq:nostro-cauchy}
 u'(t)+ \slmSbd \phi (u(t))\ni 0
    \quad \text{a.e.\ in }(0,T), \quad  u(0)=u_0,
\end{equation}
for the functional $\phi$  in the
Hilbert space
\begin{equation}
  \label{eq:hilbert-coercive}
  \begin{gathered}
    \Hilbert:=V',
\quad\text{endowed with the scalar product}\\
    \la
    u,v\ra_\Hilbert:=a(A^{-1}u,A^{-1}v)=\pairing{V'}{V}{u}{A^{-1}v}= \pairing{V'}{V}{v}{A^{-1}u}  \quad
\forall u,v \in V'.
\end{gathered}
\end{equation}
Note that,
 in this framework, the Fr\'echet
and the (strong and weak) limiting subdifferentials of $\phi$ have
to be considered with respect to the scalar product
\eqref{eq:hilbert-coercive}.
\begin{proposition}
\label{prop:subdif-repre} The functional $\phi: \Hilbert \to
[0,+\infty]$ defined by \eqref{eq:right-functional} has $D(\phi)=H$,
is lower semicontinuous on the Hilbert space $\Hilbert$
\eqref{eq:hilbert-coercive}, and complies with \eqref{eq:coerc-ass}
 and
\eqref{eqn:cont}. Moreover,
  for every $u\in H$,
\begin{equation}
\label{min-set-phi-nostro}
    \chi\in M(u)\quad\Rightarrow\quad
    \chi+\FrSbd F(\chi)\ni u,
 \end{equation}
while for every $ u \in D(\lmSbd \phi) $
  \begin{equation}
    \label{eq:limit-subdif-repre}
    \xi\in\lmSbd \phi(u)\quad
    \Rightarrow\quad
    \exists\, \chi\in M(u):\quad
    u-\chi\in V,\quad
    \xi= A(u-\chi),
  \end{equation}
  and the same result holds for $\slmSbd \phi$.
\end{proposition}
\begin{corollary}[Gradient flows solve the system]
\label{coro:subdif} Suppose that $u_0 \in H$. Then, any solution
$u\in H^1(0,T;\Hilbert)$ of the Cauchy problem
\eqref{eq:nostro-cauchy} in the Hilbert space
\eqref{eq:hilbert-coercive}
   fulfils
\begin{equation}
\label{eq:solution-pair}
  \begin{gathered}
  \text{$u\in L^\infty(0,T;H)$ and
 there exists $\chi\in L^\infty(0,T;H)$ with}
 \\
    u-\chi\in L^2(0,T;V),\quad
    \chi(t)\in M(u(t))\quad\forae\ t\in (0,T),
    \\
\text{and  the pair $(u,\chi)$
  solves the system \eqref{eq:qstat-evolution}.}
  \end{gathered}
\end{equation}
\end{corollary}

In view of the above results, in \cite{Rossi-Savare04}
 the existence of solutions of Problem \ref{quasi-stationary}
is deduced from the general Theorem~\ref{teor:exist-gflows},
 applied to the
 Cauchy problem
\eqref{eq:nostro-cauchy} (with the choice
\eqref{eq:right-functional} for $\phi$). As a consequence, the
following result has been obtained (see
\cite[Thm.~5.8]{Rossi-Savare04})
\begin{theorem}
\label{exis-qstat-evol} In the setting of
 \eqref{eq:hilb-triplet}, \eqref{eq:operator-A},
 \eqref{eq:coercive-case},
suppose that $F$
 complies with \eqref{eq:gflows_RosSav_2004:39} and
 \emph{either} with
  \begin{equation}
    \label{eq:chain-rule-1}
\begin{gathered}
\forall \, M\ge0 \ \exists\, \rho<1, \,\gamma\ge0
\ \text{such that this
  \emph{a priori estimate} holds:}
\\
    \left.
      \begin{gathered}
        u\in V,\quad \chi\in M(u),\\
        \max(\|u\|_H,F(\chi))\le
        M
      \end{gathered}
      \right\}
      \quad\Rightarrow\quad
    \chi\in V,\quad
    \|A\chi\|_{V'}\le\rho \|Au\|_{V'}+\gamma.
\end{gathered}
  \end{equation}
\emph{or} with
\begin{gather}
\label{eq:intermediate}
\begin{gathered}
 \text{there exists a Banach space $W$ such that} \ \
V\subset
W\subset H \\
\ \ \text{with   continuous inclusions,}
\\
\text{$H$ satisfies the interpolation property} \quad
  (W,V')_{1/2,2}\subset H,
    \end{gathered}
\\
\text{and for every  $M\ge0$ there exists   $C>0$ such that
  this \emph{a priori estimate} holds:}
\nonumber
\\
    \label{eq:stima-intermediate}
    \left.
      \begin{gathered}
        u-\chi\in V,\quad \chi\in M(u)\\
        \max(\|u\|_H,F(\chi))\le
        M
      \end{gathered}
      \right\}
      \quad\Rightarrow\quad
      \|\chi\|_W\le C\left(1+\|A(u-\chi)\|_{V'}\right).
    \end{gather}
Then,  for every $u_0\in H $ and $T>0$,
  Problem \ref{quasi-stationary} admits  a solution $(u,\chi)$, with
$u\in H^1(0,T;V')\cap L^\infty(0,T;H)$,
 $\chi\in L^\infty(0,T;H)$,
  $u-\chi\in L^2(0,T;V)$, fulfilling $u(0)=u_0$, the system
 \begin{equation}
  \label{eq:energy-solution}
    \begin{cases}
      u'(t)+A(u(t)-\chi(t))=0\quad\text{in }V' \quad \forae \ t \in
      (0,T),
      \\
      \chi(t)\in M(u(t))\quad\text{in $H$} \quad \forall t \in
      (0,T),
    \end{cases}
\end{equation}
and the
 \emph{energy identity}
  \begin{equation}
  \label{eq:gflows_RosSav_2004:86}
     \int_s^t a(u(r)-\chi(r)) \,dr+ \mathscr F(u(t),\chi(t)) =      \mathscr F(u(s),\chi(s))
 \quad \forall \,  0\leq s \leq t \leq T.
  \end{equation}
\end{theorem}
\noindent
 Let us stress that
Theorem~\ref{exis-qstat-evol} yields the existence of a special
class of solutions of Problem \ref{quasi-stationary}, satisfying in particular the
energy identity \eqref{eq:gflows_RosSav_2004:86}. 

\subsubsection{The weakly  coercive case: existence by an
approximation argument}
\label{s:5.1.2}
 In \cite[Sec.
5]{Rossi-Savare04}, it has been  shown that, in the
 setting of
\eqref{eq:hilb-triplet}-\eqref{eq:gflows_RosSav_2004:39} and
\eqref{eq:weakly-coercive-case}, the same conclusions of Theorem
\ref{exis-qstat-evol} hold. The proof of
 this result
  is
performed
 by
an approximation technique which we briefly recall. In fact, this
procedure
 has inspired our approach to the study of the long-time behaviour
of the solutions of Problem \ref{quasi-stationary} in the weakly
coercive case (cf. Section \ref{s:5.2.2} later on).

For any $\lambda
>0$ we consider the \emph{coercive} bilinear forms
$ a_\lambda(u,v):=a(u,v)+\lambda (u,v)_H\quad
    \forall\, u,v\in V$
and the related operators $A_\lambda: V \to V'$. Theorem
\ref{exis-qstat-evol}  yields the existence
of a solution pair $(u_\lambda, \chi_\lambda)$ to the Cauchy problem
\begin{equation}
\label{probl:lambda}
    \begin{cases}
      u_{\lambda}'(t)
      +A_\lambda(u_{\lambda}(t)-\chi_{\lambda}(t))=0\quad\text{in }V' \quad \forae
      \ t \in (0,T),
      \\
      \chi_{\lambda}(t)\in M(u_{\lambda}(t)) \quad\text{in $H$}
\quad \forall t \in (0,T),
      \\
      u_{\lambda}(0)=u_0,
    \end{cases}
\end{equation}
 fulfilling for any $T>0$ the  energy identity
\begin{gather}
      \int_s^t a_\lambda (u_{\lambda}(r)-\chi_{\lambda}(r)) \,dr+
       \mathscr F(u_{\lambda}(t),\chi_{\lambda}(t))
       =  \mathscr F(u_{\lambda}(s),\chi_{\lambda}(s))\nonumber\\
 \quad \forall \, 0\leq s \leq t \leq
      T.
      \nonumber
\label{eq:identity-lambda}
  \end{gather}
  Then, it is possible to show that the  sequences $\{u_\lambda\} \subset
  H^1 (0,T;V') \cap L^\infty (0,T;H)$ and $\{\chi_\lambda \} \subset L^\infty
  (0,T;H)$  in fact approximate a solution of Problem
  \ref{quasi-stationary}.
We have  the
following existence and approximation result (cf.  \cite[Thm. 5.9]{Rossi-Savare04}).
\begin{theorem}
\label{teor:exist-w-coerc} Assume
\eqref{eq:hilb-triplet}-\eqref{eq:gflows_RosSav_2004:39} and
\eqref{eq:weakly-coercive-case}, and let $F$
 fulfil either \eqref{eq:chain-rule-1}
or~\eqref{eq:intermediate}-\eqref{eq:stima-intermediate}. Let
$\{(u_\lambda, \chi_\lambda)\}_\lambda$ be the sequence of solution
pairs to \eqref{probl:lambda}. Then,
 there exists a
subsequence $\lambda_k \down 0$ as $k \up +\infty$ and a
pair $(u,\chi)$ such that $u \in H^1 (0,T;V') \cap L^\infty
(0,T;H)$,
 $\chi\in L^\infty(0,T;H)$,
  $u-\chi\in L^2(0,T;V)$,
and the following convergences hold:
\begin{equation}
    \label{e:convergenze-semif}
\begin{gathered}
     u_{\lambda_k} \to u \quad \text{strongly in $C^0 ([0,T];V')$,}
\\
\phi \circ u_{\lambda_k} \to \phi \circ u \quad \text{uniformly on
$[0,T]$.}
\end{gathered}
\end{equation}
Moreover, the pair $(u,\chi)$ fulfils $u(0)=u_0$, the system
\eqref{eq:energy-solution},  and the
 {energy identity}~\eqref{eq:gflows_RosSav_2004:86}.
\end{theorem}



\subsection{Long-time behaviour for
general quasi-stationary evolution systems}
 \label{sez:6}
 \noindent
This section is devoted to the investigation of the long-time
behaviour of the solutions of the evolution problem
 \begin{equation}
  \label{qstat-infto}
  \left\{
    \begin{aligned}
    &  u'(t)+A(u(t)-\chi(t))=0\quad &\text{in }V' \quad \forae \, t \in (0,+\infty),
      \\
     & \chi(t)+\FrSbd F(\chi(t))\ni u(t)\quad& \text{in $H$} \quad \forae \, t \in (0,+\infty).
    \end{aligned}
  \right.
\end{equation}
In doing so,  we maintain the distinction between the two  cases: 1. the form $a$ is \emph{coercive} and 2. the form $a$ is \emph{weakly coercive}.

In the coercive case, we shall keep to the abstract gradient flow
approach of \cite{Rossi-Savare04} (cf. Section \ref{s:5.1.1}),
 and analyze the long-term behaviour of the
 solutions of \eqref{qstat-infto} derived from the
related gradient flow equation \eqref{eq:nostro-cauchy}. We shall
refer to such solutions as {\itshape
  energy solutions} (cf. Definition \ref{ene-sol} below).
More precisely, by using the abstract results presented in the
former Section \ref{sez:4}, we will show that the set of the
 {\itshape
  energy solutions}
   of \eqref{qstat-infto}
    is a generalized
   semiflow,  which possesses  a  Lyapunov stable global attractor.
On the other hand, in
 the weakly coercive case we shall  follow the approximation
     approach outlined in Section \ref{s:5.1.2}.
Specifically, we will only consider  the solutions of
\eqref{qstat-infto}  which are limits of \emph{energy solutions} of
the approximate  coercive problem \eqref{probl:lambda-l} below.
These \emph{limiting energy solutions} form a weak generalized semiflow
 (in the sense of Section \ref{sez:2}), which possesses a weak global attractor.
\subsubsection{The coercive case}
\label{s:5.2.1}
\begin{definition}[Energy solutions]
\label{ene-sol} We say that a function $u\in H^1(0,T;V')\cap
L^\infty(0,T;H)$ $ \forall\, T>0
$
 is an \emph{energy solution} of
\eqref{qstat-infto} in the coercive case if $u$ solves the gradient
flow equation
\begin{equation}
\label{e:ns-semiflow}
\begin{gathered}
u'(t)+ \slmSbd \phi (u(t))\ni 0
    \quad \forae \ t \in (0,+\infty),
    \\
    \text{in the Hilbert space $\Hilbert:=V'$, for the functional}
    \\
 \phi(u):=
\begin{cases}
 \inf_{\chi\in H}  \left( \frac12\|u-\chi\|_H^2+F(\chi) \right)
   &  u \in H,
\\
+\infty & u \in V'\setminus H.
\end{cases}
\end{gathered}
\end{equation}
We denote by $\esseci$ the set of all energy solutions.
 \end{definition}
Note that this definition focuses on the role of the solution
component $u$, rather than on $\chi$. In order to study the long-time behaviour of
the energy solutions of \eqref{qstat-infto},
 we  shall apply our abstract results Theorem \ref{teor:1} and
Theorem \ref{teor:2} in the framework of  the phase space (cf. with
\eqref{eq:phase-space})
\begin{gather}
\label{phase-space-coercive} \cx:=D(\phi)=H, \quad \text{with} \quad
\dcx(u,v):= \sqrt{a(A^{-1}(u-v))} + |\phi(u)-\phi(v)|\nonumber\\
 \quad \forall
u,v \in H,\nonumber
\end{gather}
where as usual we have used the notation $a(w):=a(w,w)$ for $w \in
V.$

As we have recalled in Section~\ref{s:5.1.1} (cf.
Proposition~\ref{prop:subdif-repre}), under the assumption
\eqref{eq:gflows_RosSav_2004:39} the potential $\phi$ in
(\ref{e:ns-semiflow}) is lower semicontinuous on $\Hilbert$ and
complies with \eqref{eq:coerc-ass} and with  the coercivity
condition (\ref{ass:bounded-below}) (since it takes positive
values). On the other hand, the chain rule
\eqref{eq:chain-rule-energy} holds true for $\phi$ once we assume
\eqref{eq:chain-rule-1} or
\eqref{eq:intermediate}-\eqref{eq:stima-intermediate}. Hence,
Theorem \ref{teor:1}  guarantees that $\esseci$~is a generalized
semiflow.

In order to apply Theorem  \ref{teor:2},  we shall  check that
$\phi$ complies with \eqref{ass:strong-coercivity} and with
\eqref{eq:rest-points}, with the choice
$
\mathcal{D}=X=H,
$
 cf.
\eqref{phase-space-coercive}. Preliminarily, we need the following
lemma (in fact, a direct corollary of Proposition
\ref{prop:subdif-repre}), which sheds light on the set
${Z}(\esseci)$ of the rest points of the semiflow $\esseci$.
\begin{lemma}
\label{lemma:stationary-point} Assume
\eqref{eq:hilb-triplet}-\eqref{eq:gflows_RosSav_2004:39}
and
\eqref{eq:coercive-case}.
Then,
\begin{gather}
\forall \, \bar{u} \in {Z}(\esseci)=\{ u \in H \, : \ \slmSbd\phi(u) \ni 0 \}
\quad\nonumber\\
    \exists\, \bar{\chi}\in M(\bar u):\quad
    \bar u-\bar{\chi}\in V,\quad
 A(\bar u-\bar{\chi})=0.\label{eq:repre-rest}
\end{gather}
\end{lemma}
\begin{proposition}
\label{prop:rest-point} Under the  assumptions of Lemma
\emph{\ref{lemma:stationary-point}}, suppose further that the functional
$F:H \to [0,+\infty]$ fulfils:
\begin{enumerate}
\item
there exist constants $ \kappa_1, \ \kappa_2 >0$ such that for all $
\chi\in D(F)$
\begin{equation}\label{eqn:coercivity}
F(\chi) \geq \kappa_1 \|\chi\|_{H}^{2} -\kappa_2,
\end{equation}
\item and either one of the following
\begin{enumerate}
\item the proper domain of $F$
\begin{equation}
\label{eq:bounded-domainF}
 \text{ $D(F)$ is bounded in the metric
space $(\cx, \dcx)$,}
\end{equation}
 \item there
exist two constants $L_1, $ $L_2 > 0$ such that for all $\chi \in
D(\FrSbd F)$
\begin{equation}
\label{eq:coercivity-subdif-F} (\xi, \chi)_H \geq L_1 \| \chi\|_H -L_2
\quad \forall \xi \in \FrSbd F(\chi).
\end{equation}
\end{enumerate}
\end{enumerate}
Then, the potential $\phi$ in \eqref{e:ns-semiflow} satisfies the
coercivity condition \eqref{ass:strong-coercivity}. Furthermore,
 the set
${Z}(\esseci)$ of the rest points for  $\esseci$ fulfils
\begin{equation}
\label{bounded-rest-points} \text{${Z}(\esseci)$  is bounded
in $(\cx, \dcx)$.}
\end{equation}
\end{proposition}
\emph{Proof.} Preliminarily, let us recall the representation
formula \eqref{decomp-phi} for $\phi$, and let us fix an  element
$\overline{\chi} \in D(F).$ Noting that
 \begin{equation}
 \label{eq:stima-g}
  \sup_{\chi\in H}\left((u,\chi)_H-
  \Big(\frac12\|\chi\|_H^2+F(\chi)\Big)\right)
 \geq -\frac14 \| u \|_H^2 -\frac32\|\overline{\chi} \|_H^2
 - F(\overline{\chi}),
\end{equation}
we deduce from \eqref{decomp-phi} that
 there exists a   constant $J_3
\geq 0$, only depending on the chosen $\overline \chi$, such that
\begin{equation}
\label{eq:growth-phi}
 \phi(u) \leq \frac34 \| u \|_H^2  + J_3
\quad \forall u \in H,
\end{equation}
i.e., $\phi$ has at most a quadratic growth.
 In order to show \eqref{ass:strong-coercivity}, let us
note that, by elementary computations and (\ref{eqn:coercivity}),
 there holds
 \begin{equation}
\label{eqn:provacoerc2}
\begin{gathered}
\frac12 \| u -\chi\|_H^2 +F(\chi) \geq \frac12 \|u\|_H^2 +  \frac12
\|\chi\|_H^2 - (u,\chi)_H +F(\chi)
\\
\geq \frac{\kappa_1}{1+2\kappa_1} \| u \|_H^2 - \kappa_1 \|\chi
\|_H^2 +F(\chi) \geq \frac{\kappa_1}{1+2\kappa_1} \| u \|_H^2 -
\kappa_2
 \;\;\;\;\forall \chi\in L^2(\Omega).
 \end{gathered}
\end{equation}
Hence, by taking the infimum with respect to $\chi$ and recalling
the
 definition (\ref{eq:right-functional}) of $\phi$,  we deduce that $\phi$ controls the $H$-norm and \eqref{ass:strong-coercivity} ensues.

Now,  we have to prove the boundedness of the set ${Z}(\esseci)$
under either the assumption (\ref{eq:bounded-domainF}) or
(\ref{eq:coercivity-subdif-F}). We start by showing that
${Z}(\esseci) \subset D(F)$. Indeed, let $\bar{u} $ be an arbitrary
element of $ {Z}(\esseci)$. It follows from Lemma
\ref{lemma:stationary-point} and from the coercivity of $A$ that
there exists $\chi \in M(\bar{u} )$ such that $ \chi=\bar{u} $. In
particular, $\bar{u} \in M(\bar{u} ) \subset D(F)$. Thus, if
\eqref{eq:bounded-domainF} holds, \eqref{bounded-rest-points} is
trivially proved. Let us alternatively  assume
\eqref{eq:coercivity-subdif-F}. From $\bar{u} \in M(\bar{u} )$ we
infer $\,0 \in \FrSbd F(\bar{u})$.  Then,
\eqref{eq:coercivity-subdif-F}  yields $\,\| \bar u \|_H \leq
{L_2}/{L_1}\,$, whence we deduce  \eqref{bounded-rest-points} owing
to \eqref{eq:growth-phi}. \fin
\\
In view of  Proposition \ref{prop:rest-point}, Lemma
\ref{lemma:stationary-point},   and Theorem \ref{thm:ball2},
 we have the
following
\begin{theorem}
\label{teor:exist-attract-coercive}
Let \eqref{eq:hilb-triplet}-\eqref{eq:gflows_RosSav_2004:39},
 \eqref{eq:coercive-case}, \eqref{eqn:coercivity}, and either
\eqref{eq:bounded-domainF} or \eqref{eq:coercivity-subdif-F} hold.
 Further, assume  that $F$  complies   either with \eqref{eq:chain-rule-1}, or
 with
\eqref{eq:intermediate}-\eqref{eq:stima-intermediate}.
 Then, the
set $\esseci $ of  the energy solutions of the evolution problem
\eqref{qstat-infto} is a generalized semiflow in the phase space
$X=D(\phi)=H$, endowed with the metric \eqref{phase-space-coercive},
and  $\esseci$
  satisfy the
continuity property \emph{(C4)}.
 Moreover, $\esseci$
possesses a unique  global attractor $\att_\esseci$, which is
Lyapunov stable.
Finally, for any trajectory $u \in \esseci\,$ and all $\, u_\infty
\in \omega (u)$, there holds $\, 0 \in \FrSbd F
(u_\infty)$.
\end{theorem}

\subsubsection{The weakly coercive case}
\label{s:5.2.2}
In the setting of \eqref{eq:hilb-triplet}-\eqref{eq:operator-A} and
\eqref{eq:weakly-coercive-case}, we shall work in the phase space
\begin{equation}
\label{ps:weak-coercive} X=D(\phi)=H, \quad \dcxw(u,v):= \|
u-v\|_{V'} + |\phi(u)-\phi(v)| \quad \forall u,v \in H,
\end{equation}
where $\phi$ is defined by \eqref{e:ns-semiflow}. Along the lines of
the approximation procedure outlined in Section~\ref{s:5.2.2},
 for any
$\lambda >0$  we consider the set   $\esseci_\lambda $   of the
\emph{energy solutions} (cf. Definition \ref{ene-sol}) of the
approximate problems (cf. with \eqref{probl:lambda})
\begin{equation}
\label{probl:lambda-l}
    \begin{cases}
      u'_\lambda(t)
      +A_\lambda(u_\lambda(t)-\chi_\lambda(t))=0\quad\text{in }V' \quad \forae
      \ t \in (0,+\infty),
      \\
      \chi_\lambda(t)\in M(u_\lambda(t)) \quad\text{in $H$}
\quad \forall t \in (0,+\infty).
    \end{cases}
\end{equation}
 Now, we  may introduce the class of solutions of
 \eqref{qstat-infto} to which we shall restrict our investigation.
\begin{definition}[Limiting energy solutions.]
\label{def:limiting-energy-solu} We say that a function $u \in
H^1(0,T;V') \cap L^\infty (0,T;H) $ for all $T>0$
 is a \emph{limiting energy
solution} to the evolution  problem \eqref{qstat-infto}  in the
weakly coercive case, if $u$ fulfils the system
\eqref{eq:energy-solution} a.e. on $(0, +\infty)$, the energy
identity \eqref{eq:gflows_RosSav_2004:86} for all $0 \leq s \leq t
<+\infty, $ and there exists a sequence $\{\lambda_k\}$, $\lambda_k
\down 0$ as $k \up +\infty$, and a sequence $u_{\lambda_k}\in
\esseci_{\lambda_k} $ for all $k $, such that
\begin{equation}
\label{e:approx} u_{\lambda_k} \to u \quad \text{in $X$  locally
uniformly on} \ \ [0,+\infty).
\end{equation}
We denote by $\ovesse$ the set of all limiting energy solutions.
\end{definition}
\noindent Once again,
in this definition we only
focus on the role of the variable $u$. In fact, as it will be clear
from the sequel,
 for any $u \in \ovesse$ there exists a
function $\chi \in L^\infty (0,T;H)$ for all $T>0$ such that $u-\chi
\in L^2(0,T;V) $ for all $T>0$ and \eqref{eq:energy-solution},
\eqref{eq:gflows_RosSav_2004:86} hold on $[0,+\infty)$, cf. the
proof of Proposition \ref{lemma:prep-1}. Of course,
Definition~\ref{def:limiting-energy-solu} has been inspired by the
existence Theorem~\ref{teor:exist-w-coerc}, ensuring that the set
$\ovesse$ is non-empty and indeed complies with the axiom (H1) of
the definition of a generalized semiflow. In the forthcoming
Propositions~\ref{lemma:prep-1}, \ref{lemma:prep-2} we shall get
further insight into the semiflow properties of $\ovesse$.
\begin{proposition}
\label{lemma:prep-1}
 Assume
\eqref{eq:hilb-triplet}-\eqref{eq:gflows_RosSav_2004:39} and
\eqref{eq:weakly-coercive-case}, and let $F$ fulfil either
\eqref{eq:chain-rule-1} or~\emph{(\ref{eq:intermediate})-
(\ref{eq:stima-intermediate})}. Then, $\ovesse$ is a weak
generalized semiflow complying with \emph{(C4)}, and its elements
are continuous functions on $[0,+\infty)$.
\end{proposition}
\par\noindent \emph{Proof.}
Axiom (H2) can be trivially checked. The elements of $\ovesse$ are continuous on $[0,+\infty)$ since  $u \in C^0 ([0,T];V')$ for all $T>0$ and the
energy identity \eqref{eq:gflows_RosSav_2004:86} ensures that
$\,\phi \circ u \,$ is locally absolutely continuous on
$[0,+\infty).$

 In order to verify  (C4) (which obviously yields
(H4)), let us fix a sequence $\{u_n\} \subset \ovesse$ such that
$u_{n}(0) \to u_0$ in $X$, i.e. $u_{n}(0) \to u_0$ in $V'$ and
$\phi(u_{n}(0)) \to \phi({u}_0)$. We aim to show that there exists
$u \in \ovesse$ such that, up to a subsequence,
\begin{equation}
\label{remains}
 \text{$u_n$ converges to $u$ in $X$ locally
uniformly on $[0,+\infty)$.}
\end{equation}
To this purpose, we note that, by definition of $\ovesse$,  for all
$\, n \,$ there exists a sequence $\, \{ u_n^{\lambda_k}\}_{k}
\subset \esseci_{\lambda_k}\,$
 such that $\, u_n^{\lambda_k} \to u_n $ as $k \up +\infty$
 locally uniformly on  $\, [0,+\infty)$. In particular, we can
choose some increasing sequence $\,  \{\lambda_{k_n}\} \,$ (in
short: $\{\lambda_{n}\}$)
 in such a
way that
\begin{equation}
\label{e:convergenze-0-n}
 \sup_{t \in [0,n]}\dcxw\big(u_n(t), u_n^{\lambda_n}(t)\big)\leq \frac1n.
\end{equation}
Whence, in particular, $\,  u_n^{\lambda_n}(0) \to {u}_0$ in $X$.
Thus we have that $\phi(u_n^{\lambda_n}(0)) \leq C$ for a constant
independent of $n \in \N$. The energy identity
\eqref{eq:identity-lambda} for  the pair $(u_n^{\lambda_n},
\chi_n^{\lambda_n})$ reads on the interval $[0,n]$:
\begin{gather}
\label{ide-lambda-n}
 \int_s^t a_{\lambda_n} (u_n^{\lambda_n}(r)-\chi_n^{\lambda_n}(r)) \,dr+
       \mathscr F(u_n^{\lambda_n}(t),\chi_n^{\lambda_n}(t))  =
           \mathscr F(u_n^{\lambda_n}(s),\chi_n^{\lambda_n}(s))
           \end{gather}
 for all $\, 0\leq s \leq t \leq
      n.$
Using that $ \mathscr F(u_n^{\lambda_n},\chi_n^{\lambda_n})\ge
 \frac12\|u_n^{\lambda_n}- \chi_n^{\lambda_n}\|_H^2 $,
 that the sublevels of $F$
 are bounded in $H$ and the first of \eqref{probl:lambda-l},
we deduce the a priori estimates
$$
\|u_n^{\lambda_n} \|_{H^1 (0,n;V')} + \|u_n^{\lambda_n}-
\chi_n^{\lambda_n} \|_{L^2 (0,n;V) \cap L^\infty(0,n;H)} + \|
\chi_n^{\lambda_n}\|_{L^\infty(0,n;H)} \leq C
$$
for a constant independent of $n \in \N.$ Thus, suitable compactness
results and a diagonal argument
 yield that there exist subsequences
 $\{u_{n_j}^{\lambda_{n_j}}\}$ and
$\{\chi_{n_j}^{\lambda_{n_j}}\}$  (we will use the short-hand
notation $\{\lambda_j\}$, $\{u_j\},$ and $\{\chi_j\}$),
 and a pair of functions $(u, \chi_*)$, with $u \in
H^1 (0,T;V') \cap L^\infty (0,T;H) $, $\chi_* \in L^\infty (0,T;H)$
and $u-\chi_* \in L^2(0,T;V) $ for all $T>0$,
 for which the following
convergences hold as $j \up \infty$:
\begin{equation}
\label{e:convergenze-intermedie}
   \begin{gathered}
  u_j \weaksto u \quad  \text{ in $H^1(0,T;V') \cap L^\infty (0,T;H)$,}
    \\
    u_j \to u  \quad  \text{ in $C^0([0,T];V')$,}
    \ \ \forall\,T>0,
\\
u_j(t) \weakto u(t) \ \  \text{in $H$ for any $t \in (0,+\infty)$,}\\
 \chi_j
 \weaksto \chi_*  \quad \text{in $L^\infty
(0,T;H)$} \ \ \\ \text{\
  and \ $u_j-\chi_j \weakto
u-\chi_*$ in $L^2 (0,T;V)$} \quad \forall \, T>0.
 \end{gathered}
 \end{equation}
Note that the pointwise weak convergence of $u_j$ follows from the
generalized Ascoli
 theorem \cite[Cor. 4]{Simon86}. In particular, $u(0)=u_0$.
Hence, $(u,\chi_*)$ fulfils
$$
u'(t)+A(u(t)-\chi_{*}(t))=0 \quad \text{and} \quad \chi_{*}(t)\in
\conv(M(u(t))) \quad\forae \ t \in (0,+\infty).
$$
Moreover, taking the $\limsup$ as $j \up +\infty$ of the energy
identity \eqref{ide-lambda-n} with $s=0$, we get for all $T>0$
\begin{equation}
\label{eq:disug-lambda}
\begin{gathered}
    \int_0^t a(u(r)-\chi_{*}(r))\,dr+ \phi(u(t))\\ \le
      \limsup_{j \up +\infty}
      \int_0^t a_{\lambda_j} (u_j(r)-\chi_j(r)) \,dr+
\phi(u_j(t)) \leq
      \phi(u_0)
    \\
     = \phi(u(t))+ \limsup_{j \up +\infty}
      \int_0^t a_{\lambda_j}(u(r)-\chi_{*}(r))\,dr
\\
 =\phi(u(t))+ \int_0^t a(u(r)-\chi_{*}(r))\,dr
  \end{gathered}
  \end{equation}
     $ \forall t \in
      [0,T].$
Indeed,  in \eqref{eq:disug-lambda}  we have  used that, thanks to
either \eqref{eq:chain-rule-1} or to
\eqref{eq:intermediate}-\eqref{eq:stima-intermediate}, for any $T>0$
the map $\phi \circ u \in AC(0,T)$, and that, for any fixed $j \in
\N$, the following chain rule holds:
$$
 \frac d{dt}(\phi \circ u) =   \la u',u- \chi_{*}\ra  =  \la A_{\lambda_j}
  (u-\chi_{*}),u-\chi_{*}\ra= a_{\lambda_j} (u-\chi_{*}) \quad
    \aein \ (0,T),
$$
see also the proof of \cite[Thm. 5.9]{Rossi-Savare04}.
 Finally,  the last passage in
\eqref{eq:disug-lambda} follows from the trivial convergence
$\lambda_j (u_j -\chi_j) \to 0$ in $L^2(0,T;H)$ as $j \up \infty$.
 Thanks to the  lower semicontinuity argument also  exploited in  the final
part of the proof of Theorem \ref{teor:1}, we easily infer from
\eqref{eq:disug-lambda} that for all $T>0$
\begin{gather}
\label{quasi-conv-V}
 A(u_j -\chi_j) \to A(u- \chi_*) \ \  \text{strongly in $L^2 (0,T;V')$,}
 \nonumber\\
 \phi(u_j(t)) \to \phi(u(t)) \  \forall t \in [0,T].
 \end{gather}

 By a careful measurable selection argument, detailed in the proof
of \cite[Thm.5.9]{Rossi-Savare04}, it is possible to show that there
exists a function $\chi \in L^\infty (0,T;H)$ fulfilling
\begin{gather}
\label{e:chi-minimiz}
 \chi(t) \in M(u(t)) \quad \forall t \in (0,T),
\qquad u-\chi \in L^2 (0,T;V),
\\
\label{e:sono-uguali}
 A(u(t)- \chi_{*}(t))=A(u(t)-\chi(t)) \quad \forae \ t \in (0,T).
\end{gather}
 Being $T$ arbitrary, we conclude that  the pair $(u,\chi)$ fulfils
\eqref{eq:energy-solution} a.e.  on $(0,+\infty)$. Furthermore, from
the energy identities \eqref{ide-lambda-n} and
\eqref{eq:disug-lambda} we infer
 for all $t>0$
 $$
\begin{aligned}
&\left| \phi(u_j(t)) -\phi(u(t)) \right|\\
 &\leq |\phi(u_j (0))
-\phi(u_0)|+
  \int_0^t \left |  \|
A(u_j (s)-\chi_j (s))\|_{V'}^2 - \| A(u(s)-\chi(s))\|_{V'}^2 \right|
ds
\\
&
\leq |\phi(u_j (0)) -\phi(u_0)| + \Big( \| A(u_j -\chi_j)\|_{L^2
(0,t;V')}
 \\
& + \| A(u-\chi)\|_{L^2 (0,t;V')}\Big) \|A(u_j
-\chi_j)-A(u-\chi)\|_{L^2 (0,t;V')}.
\end{aligned}
$$
Hence, in view of \eqref{e:convergenze-0-n}
 and of \eqref{quasi-conv-V},
  we easily   conclude  (cf. \eqref{eq:conv-unif-cpt}),
that $\phi(u_j) \to \phi(u)$ locally uniformly on $[0, +\infty)$.
Combining the latter convergence with the first of
\eqref{e:convergenze-intermedie}, we find that
\begin{equation}
\label{conve-u-j}
 u_j \to u \quad \text{in $X$ locally uniformly  on $[0,+\infty)$.}
\end{equation}
Finally, owing to  \eqref{quasi-conv-V}-\eqref{conve-u-j}, we pass
to the limit in the energy identity \eqref{ide-lambda-n},  and we
deduce that the pair $(u,\chi)$ fulfils the energy identity
\eqref{eq:gflows_RosSav_2004:86} for all $0 \leq s \leq t <+\infty$.
By the previous construction, $u$ is approximated in the sense of
\eqref{e:approx}, whence $u \in \ovesse$.

In the end,
  one directly  checks that,
for all $\, T < +\infty\,$ and $\, {n_j} > T$,
\begin{gather}
  \sup_{t \in [0,T]}\dcxw\big(u(t),u_{n_j}(t)\big)\leq \sup_{t \in
[0,T]}\dcxw\big(u(t), u_j(t)\big) + \sup_{t \in [0,T]}\dcxw\big(
u_j(t),u_{n_j}(t) \big)\nonumber\\ \leq \sup_{t \in
[0,T]}\dcxw\big(u(t), u_j(t)\big)+ \frac{1}{n_j},\nonumber
\end{gather}
also in view of \eqref{e:convergenze-0-n}. Owing to
\eqref{conve-u-j}, we conclude the convergence \eqref{remains},
 and (C4) ensues.
 \fin
\begin{proposition}
\label{lemma:prep-2} Under the same hypotheses of Proposition
\emph{\ref{lemma:prep-1}}, assume further that $F$ complies with
\eqref{eqn:coercivity}. Then, $\ovesse$ is compact and eventually
bounded.
\end{proposition}
\emph{Proof.} Let us point out that, by Definition
\ref{def:limiting-energy-solu}, the limiting energy solutions of
Problem \ref{quasi-stationary} comply with the energy identity
\eqref{eq:gflows_RosSav_2004:86} just like the energy solutions
deriving from the gradient flow equation \eqref{e:ns-semiflow}.
Thus, the eventually boundedness of $\ovesse$ follows exactly by the
same argument developed in the proof of our abstract Theorem
\ref{teor:2} (cf.
\eqref{e:ad-event-bounded-1}-\eqref{e:ad-event-bounded-2}), since
assumption  \eqref{eqn:coercivity} provides the sufficient
coercivity (cf. the proof of Proposition \ref{prop:rest-point}).

In order to prove that $\ovesse$ is compact, we fix a sequence $u_n
\in \ovesse$ such that $u_{n}(0)$ is bounded in $X$. The same
computations as in the proof of Proposition \ref{lemma:prep-1} yield
that there exists an increasing sequence $\{\lambda_{n}\}$ and
$u_{\lambda_{n}} \in \esseci_{\lambda_{n}}$ for which
\eqref{e:convergenze-0-n} holds. In particular,  note that
$\{u_{n}^{\lambda_{n}}(0)\}$
 is
bounded in $X$. Hence, again exploiting
 the energy identity \eqref{ide-lambda-n}  for the
pair $(u_n^{\lambda_n}, \chi_n^{\lambda_n}) $, we infer that there
exists a subsequence (which we do not relabel) and a limit pair $(u, \bar
\chi)$ for which the convergences \eqref{e:convergenze-intermedie}
 hold true on $(0,+\infty).$
However, since  in this  case we cannot conclude anymore that
$\{u_{n}^{\lambda_{n}}(0)\}$
 converges, we cannot
exploit
  the proof of Proposition \ref{lemma:prep-1}
 in order to conclude that $u_n^{\lambda_n}$
 converges to $u$ locally uniformly on $[0,+\infty)$.
  Instead, we
 will argue in the same way as in the proof  of the compactness property
 in  Theorem \ref{teor:2}. Let us sketch this procedure. First,
the energy identity
  \eqref{eq:identity-lambda} yields
  that the map $t \mapsto \phi(u_n^{\lambda_n}(t))$
  is non-increasing. By Helly's Theorem,
    for all $t >0$ the function
  $\varphi(t):= \lim_{n \up +\infty} \phi(u_n^{\lambda_n}(t))$ is
  well-defined. Moreover, \eqref{eq:identity-lambda} and
  Fatou's Lemma entail that
  \begin{gather}
\liminf_{n \up+ \infty} \| A(u_n^{\lambda_n}(t) -
\chi_n^{\lambda_n}(t)) \|_{V'}^2 \nonumber\\ + \sup_{n} \left(\frac12 \|
u_n^{\lambda_n}(t)- \chi_n^{\lambda_n}(t) \|_{H}^2 +
F(\chi_n^{\lambda_n}(t)) \right) <+\infty \quad\nonumber
  \end{gather}
  (where $\chi_n^{\lambda_n} \in M(u_n^{\lambda_n})$)
for almost every $\, t >0$. Also using the compactness of the
sublevels of $F$ \eqref{eq:gflows_RosSav_2004:39}, one easily infers
that for almost any $t>0$ there exist a subsequence $j \mapsto n_j$,
possibly depending on $t$, and a pair $(\hat{u}(t),\hat{\chi}(t) )$
for which (using short-hand notation) $\chi_j (t) \to \hat{\chi}(t)$
and $u_j(t)-\chi_j(t) \to \hat{u}(t)-\hat{\chi}(t) $ strongly in
$H$. Thus, $u_j(t) \to \hat{u}(t)$ in $H$, whence necessarily
$\hat{u}(t)=u(t)$ for a.e. $t \in (0,T)$ thanks to
\eqref{e:convergenze-intermedie}. Finally, it is not difficult to
check that $\hat{\chi}(t) \in M(u(t))$, and that
$$
\lim_{j \up +\infty} \phi(u_j(t))= \phi(u(t)) \quad \forae \, t \in
(0,T),
$$
cf. with \eqref{eqn:lim2}. Arguing as in \eqref{eqn:lim3}, we
finally deduce that $\varphi(t)=\phi(u(t))$ for a.e.
$t~\in~(0,+\infty).$ Thus, exactly as in the proof of Theorem
\ref{teor:2} we may  pass to the limit in \eqref{eq:identity-lambda}
for all $t>0$ and  for a.e. $s \in (0,t)$ for which
$\varphi(s)=\phi(u(s)).$ We can now develop the same energy identity
argument of \eqref{eq:disug-lambda}-\eqref{quasi-conv-V} (of course
replacing $u_0$ with $u(s)$), and  we deduce $ \phi(u_j(t)) \to
\phi(u(t)) $ $ \forall t >0. $ Then, exploiting
\eqref{e:convergenze-0-n}, we
  complete the proof of the compactness property.
 \fin
\par\noindent{\bf Long-time behaviour of the limiting energy solutions.}
 We
shall prove that the weak generalized semiflow $\ovesse$ of the limiting
energy solutions of \eqref{qstat-infto}
 possesses a weak
global attractor in the particular case (which is however meaningful
in view of the applications):
\begin{gather}
 V=H^1 (\Omega), \quad H=L^2 (\Omega), \nonumber\\
\quad
\pairing{H^1 (\Omega)'}{H^1 (\Omega)}{Au}{v}= \int_\Omega {\sf
A}_1\nabla u \nabla v \quad \forall \, u,v \in H^1 (\Omega).\label{e:parti-case}
\end{gather}
Here,  ${\sf A}_1:\Omega\to \M^{m\times m}$ is a field of  symmetric
matrices, with bounded and measurable coefficients, satisfying the
usual uniform ellipticity condition
\begin{equation}
  \label{eq:ellipticity}
   {\sf A}_{1}(x)\eta\cdot\eta\ge \rho>0\quad
  \forall\, x\in\Omega,\ \eta\in \R^m,\ |\eta|=1.
\end{equation}
Let us point out that,  according to Definition
\ref{def:limiting-energy-solu} and to \eqref{e:parti-case},
 any limiting  energy solution $u $ of \eqref{qstat-infto} fulfils
the system
\begin{equation}
\label{qstat-infto-part}
\begin{cases}
&
 u'(t)-\div{\sf A}_1\nabla (u(t)-\chi(t))=0 \ \ \text{in $\Omega \times (0,+\infty),$}
      \\
      &
      \chi(t)\in M(u(t)) \ \ \text{in $\Omega \times (0,+\infty),$}
      \\
 &
      {\sf A}_1\nabla(u-\chi)\cdot\nn =0 \ \ \text{in $\partial\Omega \times (0,+\infty).$}
    \end{cases}
\end{equation}
Note that $u$ is  a conserved parameter. Indeed, taking the
integral in space of the first equation in \eqref{qstat-infto-part},
one finds that the map $\, t \mapsto \int_\Omega u(t) \,$ is
constant along the evolution. This in particular implies that the
semiflow corresponding to the  limiting energy solutions of
\eqref{qstat-infto-part} is not point dissipative. In other words,
the set of stationary solutions of \eqref{qstat-infto-part} is
unbounded in $\,H^1(\Omega)'$. Eventually, no global attractor in
the phase space $\, X = H^1(\Omega)' \,$ is to be expected (this
kind of difficulty is well-known and is, for instance, discussed in
\cite[Chapter 3]{Temam88} in connection with the long-time analysis
of the
 Cahn-Hilliard equation). Hence, we shall
consider some modification of the phase space by fixing explicit
bounds on the conserved quantity $\, \int_\Omega u$. To this aim, we
use the notation
\begin{equation}\label{dm}
 m(u):= \disp\frac{1}{|\Omega|}\pairing{H^1 (\Omega)'}{H^1
(\Omega)}{u}{1},      \quad
 \mathcal{D} (\bar{m}):= \left\{u \in H^1 (\Omega)' \,:\,  m(u) \leq
\bar{m}\right\},
\end{equation}
for   given $u \in H^1 (\Omega)'$ and $\bar{m}>0$ (here
$\,|\Omega|\,$ stands for the volume of $\, \Omega$).

 Note that the energy identity
\eqref{eq:gflows_RosSav_2004:86} suggests that another choice for
the invariant region $\mathcal{D}$ could be, for a given positive
$C_\phi>0$,
\begin{equation}\label{dphi}
\mathcal{D}_\phi =\left\{ v\in X: \phi(v)\leq C_\phi\right\}.
\end{equation}
In the next Theorem, we apply the abstract results of Theorem \ref{teor:Lyapun-w-attractor} to the set $\ovesse$ of the limiting energy
solutions of \eqref{qstat-infto-part}. Although we give the proof in the case in which $\mathcal{D}$ is as in \eqref{dm}, the same results hold also when we choose $\mathcal{D}$ in \eqref{dphi}.
\begin{theorem}
\label{teor:exist-att-weak-coercive} In the setting of
\eqref{e:parti-case}, let $F$ comply with
\eqref{eq:gflows_RosSav_2004:39} and either with
\eqref{eq:chain-rule-1} or
with~\eqref{eq:intermediate}-\eqref{eq:stima-intermediate}.
Further, suppose that
\begin{equation}
\label{ci-vuole-proprio} D(\partial F) \  \ \text{is bounded in $L^2
(\Omega)$.}
\end{equation}
Then, for any $\bar{m}>0$ the set $\ovesse$ of the limiting energy
solutions of \eqref{qstat-infto-part}
 admits
 the  weak global attractor $\geneatt_{\ovesse}$
in the set $\mathcal{D} (\bar{m})$.
  Moreover, for any trajectory $u \in \gamma^+
  (\mathcal{D}(\bar{m}))$
  and for any $u_\infty
\in \omega(u)$ there exists $\chi_\infty \in M(u_\infty) $ such that
\begin{equation}
\label{e:omega-limit-traiettoria}
\begin{cases}
& -\div{\sf A}_1\nabla (u_\infty-\chi_\infty) = 0 \quad \text{in
$\Omega$},
\\
& \chi_\infty \in M(u_\infty) \quad \text{in $\Omega$}
      \\
&
      {\sf A}_1\nabla(u_\infty-\chi_\infty)\cdot\nn =0 \ \ \text{in $\partial\Omega .$}
\end{cases}
\end{equation}
\end{theorem}\noindent
Note that the assumptions \eqref{eqn:coercivity} and
\eqref{eq:bounded-domainF}-\eqref{eq:coercivity-subdif-F} of
Theorem \ref{teor:exist-attract-coercive}
 have been replaced by the stronger coercivity condition
 \eqref{ci-vuole-proprio}.
 \smallskip\\
 \emph{Proof. }
Preliminarily, it is easy to see that
assumption~\eqref{eqn:coercivity} in Proposition~\ref{lemma:prep-2}
may be replaced by \eqref{ci-vuole-proprio}. Then,
 relying on Propositions \ref{lemma:prep-1},
\ref{lemma:prep-2} we conclude that the weak generalized
 semiflow $\ovesse$ is eventually bounded and compact.
 Furthermore,
  since any $u \in \ovesse$ complies with
the energy identity \eqref{eq:gflows_RosSav_2004:86} and with
\eqref{eq:energy-solution}, we have that  $\phi$ is a Lyapunov
function for $\ovesse$, in fact arguing as in the proof of
Theorem~\ref{teor:2}. Then, in view of
Theorem~\ref{teor:Lyapun-w-attractor} it is sufficient to see that
for any $\bar{m}>0$ the set $\mathcal{D}(\bar{m})$ complies with
 conditions
   \eqref{eq:invariance-2}-\eqref{eq:rest-points-2}.

As already observed, for any trajectory  $u$ starting from  the set
$\mathcal{D}(\bar{m})$
we have $m (u'(t))=0$ for a.e. $t \in
  (0,+\infty)$. Thus,
  the invariance condition
   \eqref{eq:invariance-2} ensues.
 In order to check \eqref{eq:rest-points-2}, let us fix $\bar{u} \in
 {Z}(\ovesse) \cap \mathcal{D}(\bar{m})$.
Recalling \eqref{qstat-infto-part}, we easily see that there
 exists $\bar{\chi}\in M(\bar{u}) $ such that the pair $(\bar{u}, \bar{\chi})$
 fulfils the system \eqref{e:omega-limit-traiettoria}. In
 particular,
 $\bar{\chi} \in D(\partial F)$, so by \eqref{ci-vuole-proprio}
 there exists a constant $\bar{r}>0$
 such that
$|m(\bar \chi)| \leq \bar{r}. $ Thus, $|m(\bar{u}-\bar{\chi})|\leq
\bar m +\bar r. $ Combining this with the first of
\eqref{e:omega-limit-traiettoria} and with Poincar\'e's inequality,
we infer that there exists a positive constant $C$ independent of
$\bar u$ and $\bar \chi$ such that
$
\|\bar{u}-\bar{\chi} \|_{H^1 (\Omega)} \leq C.
$
Since $\bar{\chi}$ is bounded in $L^2 (\Omega)$ by
\eqref{ci-vuole-proprio}, we conclude that $\bar{u}$ is bounded in
$L^2 (\Omega)$. Thus, the set ${Z}(\ovesse) \cap
\mathcal{D}(\bar{m})$ is bounded
 in the phase space \eqref{ps:weak-coercive}, as $\phi $  is controlled by
  the norm on $L^2 (\Omega)$, cf. the growth estimate
  \eqref{eq:growth-phi}.

Therefore,  the existence of a weak attractor $\geneatt_{\ovesse}$
in the set ${Z}(\ovesse) \cap \mathcal{D}(\bar{m})$
 is established, and
\eqref{e:omega-limit-traiettoria} follows from the last part of the
statement of Theorem~\ref{teor:Lyapun-w-attractor}.
 \fin
\subsubsection{Approximation of the weak global attractor}\label{sec:approx-attractor}
In this section we discuss  the approximation of the weak global
attractor of the limiting energy solutions with the global attractor
of the weal generalized semiflow $\esseci_\lambda $, generated by the
solutions of the approximating scheme  \eqref{probl:lambda-l}.
 We shall denote by  $X_\phi$ the subset  $X\cap \mathcal{D}_\phi$ of the phase space
$X=D(\phi)$,  endowed with the  distance $\dcx$
\eqref{phase-space-coercive}. For any $\lambda>0$, let
$\geneatt_\lambda$ be the global attractor of the generalized
semiflow $\esseci_\lambda $ in the phase space $(X_\phi, \dcx)$,
  whose existence is
ensured by Theorem \ref{teor:exist-attract-coercive}. Further,
 let $\geneatt_{\ovesse}$  be the weak global attractor of the
set $\ovesse$ of the limiting energy solutions of
\eqref{qstat-infto-part} in the phase space $(X_\phi, \dcxw)$
\eqref{ps:weak-coercive}. Finally, we denote by
  $e_{\phi}$   the Hausdorff semidistance (or
 \emph{excess}) associated with the distance $\dcxw$.
 We  have the following
\begin{theorem}\label{upper-semicont-attr} In the setting of
\eqref{e:parti-case}, let  $F$ comply with
\eqref{eq:gflows_RosSav_2004:39} and  either
 with \eqref{eq:chain-rule-1}, or
 with
\eqref{eq:intermediate}-\eqref{eq:stima-intermediate}. Further,
assume \eqref{ci-vuole-proprio}.
 Then,
\begin{equation}\label{convergenza attrattori}
\lim_{\lambda\down
0}~e_{\phi}(\geneatt_\lambda,\geneatt_{\ovesse})=0.
\end{equation}
\end{theorem}
\emph{Proof. } In order to prove  \eqref{convergenza attrattori}
 we argue by contradiction along the lines of {\sc Hale \& Raugel}, cf. \cite{hale-raugel}.
Assume  that  \eqref{convergenza attrattori} does not hold: then, we
can find $r_0>0$ and sequences $\left\{\lambda_n\right\}_{n\in
\mathbb{N}}$ and $\left\{\xi_n\right\}_{n\in \mathbb{N}}$ such that
$\lambda_n\down 0$ and for all $n\in \mathbb{N}$
\begin{equation}\label{convattr1}
\xi_n \in \geneatt_{\lambda_n}, \quad  \inf_{\xi\in
\geneatt_{\ovesse}}\dcxw(\xi_n,\xi)\ge r_0.
\end{equation}
Now, the invariance  of $\geneatt_{\lambda_n}$ (but actually the
sole quasi-invariance  would be sufficient, see
 the proof of Theorem
\ref{teor:exist-w-att}) entails that there exists a complete orbit
$u_n$ with $u_n(0)=\xi_n$ and $u_n(t)\in \geneatt_{\lambda_n}$ for
all $t\in \mathbb{R}$. It is not difficult to see that this orbit is
bounded independently of $\lambda_n$ with respect to $\dcxw$. In
fact, the energy identity (recall that $u_n$ is in particular an
energy solution), \eqref{dphi} and the translation invariance of the
complete orbit $u_n$ entail that
\begin{equation}
\label{bound-orbita-lambda} \int_{-T}^{T}\vert u'_{n}(s)\vert^2 \,
ds +\phi(u_n(T))\le C_\phi \quad \;\;\forall T>0.
\end{equation}
The proof of Theorem \ref{teor:exist-w-coerc} in
\cite{Rossi-Savare04} (see also Propositions \ref{lemma:prep-1} and
\ref{lemma:prep-2})
 shows that this
estimate is sufficient to pass to the limit as $\lambda_n\down 0$,
obtaining in the limit a complete and bounded orbit $u$ of the set
$\ovesse$ of the limiting energy solutions  to \eqref{qstat-infto}.
In particular, there holds $\xi_n=u_n(0)\to u(0)$ in $X$. Now, since
by \eqref{struttura-attrattore} the weak global attractor is
generated by the complete and bounded orbits of $\ovesse$, we
conclude that $u(0)\in \geneatt_{\ovesse}$. This leads to
 contradiction with
\eqref{convattr1}.
 \fin

\subsubsection{Applications to  quasi-stationary phase field models}
\label{s:5.1.3}
 Let us consider
 the following quasi-stationary system, which
generalizes the quasi-stationary phase field model (cf. system
\eqref{e:corazon}-\eqref{cagi}):
\begin{equation}
 \label{BBinfto}
\begin{cases}
\partial_{t} u -\div {\sf A}_1\nabla(u-\chi)=0,
\\
 -\div {\sf A}_2\nabla \chi + \FrSbd {\mathcal{W}}(\chi) \ni u,
\end{cases}
\quad \text{in $\Omega \times (0, +\infty)$.}
\end{equation}
Here,   ${\sf A}_2:\Omega\to \M^{m\times m}$ is a  field of
symmetric matrices, with bounded and measurable coefficients,
satisfying the uniform ellipticity condition \eqref{eq:ellipticity}.
 On the other hand,  $\mathcal{W}$ is either an arbitrary  $C^1$ real function with superlinear
growth (in this case $\FrSbd F$ reduces to $\mathcal{W}'$), or a
semi-convex function valued in $\R \cup \{+\infty\}$. Meaningful
examples of $\mathcal{W}$
 are:
 \begin{align}
 \label{doppio-pozzo}
&
  \mathcal{W}(\chi):= \frac{(\chi^2 -1)^2}{4},
 \\
&
 \label{visintin} \mathcal{W}(\chi):=I_{[-1,1]}(\chi) +  (1-\chi)^2;
\\
&
 \label{logaritmico} \mathcal{W}(\chi):= c_1\left((1+\chi) \ln(1+\chi) +(1-\chi)
\ln(1-\chi)\right)-c_2 \chi^2 + c_3 \chi +c_4,
\end{align}
with  $c_1, c_2 >0$ and $c_3, c_4 \in \R$ (see e.g. \cite[4.4,
p.170]{Brokate-Sprekels96} for \eqref{logaritmico},
\cite{Blowey-Elliott91}, \cite{Visintin96} for \eqref{visintin}).
 The symbol $I_{[-1,1]}$  denotes the indicator
function of $[-1,1],$ which forces the constraint  $-1 \leq \chi
\leq 1$. In the sequel, we shall employ the notation
$
 D(\mathcal{W}):=\left\{\chi\in L^2(\Omega): \mathcal{W}(\chi)\in
L^2(\Omega)\right\}.
$

In \cite[Sec.~5]{Rossi-Savare04}, existence results have been
obtained for some   initial boundary-value problems  for
\eqref{BBinfto} on a finite time interval. Specifically,
\eqref{BBinfto} has been supplemented with
   the natural homogeneous
Neumann boundary condition on $\chi$, and with  homogeneous, either
Dirichlet or
 Neumann, boundary conditions on $u-\chi$, and the existence
results of \cite{Plotnikov-Starovoitov93} and of \cite{Schatzle00}
have been respectively recovered. Here, we shall focus on the
long-time behaviour of  \eqref{BBinfto}, supplemented with both
kinds of boundary conditions. In fact, we shall
  apply the
abstract results of Sections~\ref{s:5.2.1}, \ref{s:5.2.2} to
 suitable families of
solutions of the related boundary value problems.
\par\noindent{\bf Attractor for the  quasi-stationary phase field model
with
 Dirichlet-Neumann boundary conditions.}
We supplement \eqref{BBinfto} with the  boundary conditions
\begin{equation}
\label{dirichlet-neumann}
 u-\chi=0,\quad A_2\nabla\chi\cdot \nn=0\quad\text{in }\partial\Omega
  \times(0,+\infty).
\end{equation}
Note that the system \eqref{BBinfto}, \eqref{dirichlet-neumann}
 may
be reformulated as  the abstract evolution system
\eqref{qstat-infto} with the choices
  $V:=H^1_0(\Omega),$ $H:=L^2 (\Omega)$, $V':=H^{-1} (\Omega) $,
  $A:=-\div({\sf A}_1\nabla \cdot)$,
and with $F: L^2 (\Omega) \to [0,+\infty]$ given by
\begin{equation}
\label{F-phase-field}
 F(\chi):\begin{cases}
\disp\int_\Omega \left( \frac12 {\sf A}_2 (x) \nabla \chi(x) \cdot \nabla
\chi(x) + \mathcal{W}(\chi(x)) \right) dx \quad & \chi \in H^1
(\Omega)\cap D(\mathcal{W}),
\\
+\infty \quad & \text{otherwise}.
\end{cases}
\end{equation}
As $A$ is coercive  on $V$, we will focus on the \emph{energy
solutions} of (\ref{BBinfto})-(\ref{dirichlet-neumann}). They stem
from the gradient flow  equation \eqref{e:ns-semiflow},  in the
space $\Hilbert=H^{-1}(\Omega)$,  for the functional $\phi:H^{-1}
(\Omega) \to (-\infty, +\infty]$
\begin{equation}
\label{eq:phi-quasistat} \phi(u):=\begin{cases}
 \inf_{\chi \in H^1 (\Omega)}
\Big\{\disp \int_\Omega \frac12 |u(x)-\chi(x)|^2   +F(\chi)\Big\}, \ \\
\text{for} \ \ u \in L^2 (\Omega),
\\
+\infty \ \text{otherwise},
\end{cases}
\end{equation}
with  $\mathcal{W}$ as in \eqref{doppio-pozzo}-\eqref{logaritmico},
for instance. Hence, let us  check that the assumptions of~Theorem
\ref{teor:exist-attract-coercive} are fulfilled within this
framework. Since the matrix field ${\sf A_2}$ is uniformly elliptic,
$F$ has strongly compact sublevels in $L^2 (\Omega)$ for  all the
examples (\ref{doppio-pozzo})-(\ref{logaritmico}).
 Concerning condition~(\ref{eqn:coercivity}),  it is sufficient to show that there exist
constants $\kappa_1, \ \kappa_2
>0$ such that {
$ \int_\Omega \big(\mathcal{W}(\chi(x))-\kappa_1\vert
\chi(x)\vert^{2}\big)dx\ge -\kappa_2, $ } which is satisfied in all
cases (\ref{doppio-pozzo})-(\ref{logaritmico}).
  Also note that $F$
complies with \eqref{eq:intermediate}-\eqref{eq:stima-intermediate}
(with the choice $W=H^1 (\Omega)$). Instead, the validity of
(\ref{eq:bounded-domainF})
  (or (\ref{eq:coercivity-subdif-F})) depends on the
particular choice of the potential $\mathcal{W}$. More precisely, if
we choose the singular potentials (\ref{visintin}) or
 (\ref{logaritmico}),
  then (\ref{eq:bounded-domainF}) is easily satisfied,
   since the domain of $F$ fulfils
   \begin{equation}
   \label{dominio-limitato}
    D(F)\subseteq H^{1}(\Omega)\cap \left\{v\in L^2(\Omega): \;-1\le v(x)\le 1,
    \ \  \forae\, x\in \Omega\right\}
\end{equation}
    (the two sets coincide
     if we choose the potential $\mathcal{W}$
     in (\ref{visintin})). Thus, $D(F)$
      is clearly bounded in $L^2(\Omega)$.
       On the other hand, it is not difficult
        to control that the usual double
         well potential (\ref{doppio-pozzo}) complies
 with (\ref{eq:coercivity-subdif-F}).
Eventually, we conclude
 that the set of the energy solutions of
  \eqref{BBinfto}, \eqref{dirichlet-neumann} is  a generalized semiflow.
Such a semiflow
   possesses a
Lyapunov stable global attractor in the phase space $D(\phi)=L^2
(\Omega)$, endowed  with the distance
defined by the functional
$\phi$ \eqref{eq:phi-quasistat}.
\par\noindent{\bf Attractor for the  quasi-stationary phase field model
with Robin-Neumann boundary conditions.} We supplement
\eqref{BBinfto} with the conditions
\begin{equation} \label{robin-neumann} {\sf
A}_1\nabla(u-\chi)\cdot\nn + \omega (u-\chi)=0, \quad {\sf
A}_2\nabla\chi\cdot \nn =0 \quad \text{in
}\partial\Omega\times(0,+\infty),
\end{equation}
where $\omega>0$. This problem may be recast in the form
\eqref{qstat-infto} by setting $V:=H^1 (\Omega),$ $H:=L^2 (\Omega),
$
$$
\pairing{V'}{V}{Au}{v}:=\int_\Omega {\sf A}_1 (x) \nabla(u(x)) \cdot
\nabla v(x) dx + \omega \int_{\partial \Omega} u(s)v(s) ds,
$$
and choosing $F$ as in \eqref{F-phase-field}. Since $A$ is coercive
on $H^1(\Omega)$,
 we may again
consider the energy solutions of \eqref{BBinfto},
\eqref{robin-neumann} in the sense of  Definition \ref{ene-sol}. In
this setting, the ambient space $\Hilbert$ is $(H^1 (\Omega))'$,
with $\phi$ defined by \eqref{eq:phi-quasistat}. Hence, we may argue
exactly in the same way as for the Dirichlet-Neumann problem, with
the sole difference that now $F$ complies with
\eqref{eq:chain-rule-1}. Therefore,
  Theorem
\ref{teor:exist-attract-coercive}
 applies and we conclude the existence of a global attractor for the
 semiflow of the energy solutions of \eqref{BBinfto},
\eqref{robin-neumann}.
This gradient flow approach could also be extended to tackle more
general boundary conditions on $u-\chi$, such as homogeneous
Dirichlet (or Robin) on a portion of $\partial \Omega$, and
non-homogeneous Neumann on the remaining part.
\par\noindent{\bf Attractor for the  quasi-stationary
 phase field model with Neumann-Neumann boundary conditions.}
We supplement the system \eqref{BBinfto}
  with the  boundary conditions
\begin{equation}
\label{neumann-neumann} {\sf A}_1\nabla(u-\chi)\cdot\nn =0 , \quad
{\sf A}_2\nabla\chi\cdot \nn =0 \quad \text{in
}\partial\Omega\times(0,+\infty).
\end{equation}
Problem \eqref{BBinfto}, \eqref{neumann-neumann}
 can be rephrased in the form of  Problem
\ref{qstat-infto} by setting
\begin{equation}
\label{e:setting-neumann} V:=H^1 (\Omega), \ \ H:=L^2 (\Omega),
\quad \pairing{V'}{V}{Au}{v}:= \int_\Omega {\sf A}_1 (x)
\nabla(u(x)) \cdot \nabla v(x) dx,
\end{equation}
and $F$ as in \eqref{F-phase-field}. Note that $A$ is only weakly
coercive on $H^1 (\Omega)$. Following the outline of Section
\ref{s:5.2.2},
 we shall  focus on the  long-time behaviour of the set $\ovesse_{\text{neu}}$ of
the
 limiting energy
solutions of \eqref{BBinfto}, \eqref{neumann-neumann}. Let us  check
the conditions of Theorem~\ref{teor:exist-att-weak-coercive}. First,
note that
 $F$
satisfies \eqref{eq:intermediate}-\eqref{eq:stima-intermediate},
with $W=H^1 (\Omega)$, for the potential $\mathcal{W}$  as in
\eqref{doppio-pozzo}-\eqref{logaritmico}. On the other hand, in view
of \eqref{dominio-limitato}, condition
 \eqref{ci-vuole-proprio} holds  true only in the cases of
\eqref{visintin}-\eqref{logaritmico}. Arguing as for the
Dirichlet-Neumann and Robin-Neumann cases, it is not difficult to
see that $F$ complies with the remaining assumptions of
Theorem~\ref{teor:exist-att-weak-coercive}. Thus, we conclude that
for all $\bar{m}>0$
 $\ovesse_{\text{neu}}$ admits a unique  weak global attractor
$\mathcal{A}_{\ovesse_{\text{neu}}}$
 in
the set $\mathcal{D}(\bar{m})$,  and  that
\eqref{e:omega-limit-traiettoria} holds for
 $\omega-$limit points
of the  trajectories. Finally, referring to the notation of Section
\ref{sec:approx-attractor} (with $\phi $ defined by
\eqref{eq:phi-quasistat}),  we have that the sequence
$\{\geneatt_\lambda\}$ of the global attractors of the solutions of
the approximate problems
 \eqref{probl:lambda-l} converges to the
weak global attractor $\mathcal{A}_{\ovesse_{\text{neu}}}$ in the
sense that $\,\lim_{\lambda\down
0}~e_\phi(\geneatt_\lambda,\mathcal{A}_{\ovesse_{\text{neu}}})=0.$

\begin{acknowledgement}
 The authors would like to thank Prof.
Giuseppe Savar\'{e} for some va\-lua\-ble and inspiring
conversations.
 \end{acknowledgement}

\def\cprime{$'$}
\providecommand{\bysame}{\leavevmode\hbox to3em{\hrulefill}\thinspace}
\providecommand{\MR}{\relax\ifhmode\unskip\space\fi MR }
\providecommand{\MRhref}[2]{%
  \href{http://www.ams.org/mathscinet-getitem?mr=#1}{#2}
}
\providecommand{\href}[2]{#2}

\address{
Dipartimento di Matematica, Universit\`a di  Brescia,\\ via Valotti 9, I--25133 Brescia, Italy\\
email: {\tt riccarda.rossi\,@\,ing.unibs.it}
\and
Weierstrass-Institut f\"ur Angewandte Analysis und Stochastik,\\
Mohrenstrasse 39, D--10117 Berlin, Germany\\
email: {\tt segatti\,@\,wias-berlin.de}
\and
Istituto di Matematica Applicata e
Tecnologie Informatiche -- CNR,\\ via Ferrata 1, I--27100 Pavia,
Italy\\
email: {\tt ulisse.stefanelli\,@\,imati.cnr.it}}

\begin{thebibliography}{10}



\bibitem{Ambrosio95}
{\sc L.~Ambrosio}: {Minimizing movements}. {\it Rend. Accad. Naz. Sci. XL
Mem. Mat. Appl. (5)} {\bf 19}, 191--246 (1995).

\bibitem{Ambrosio-Gigli-Savare04}
{\sc L.~Ambrosio, N~Gigli, and G.~Savar\'e}; {\it Gradient flows in
metric spaces
  and in the {W}asserstein spaces of probability measures}, Lectures in Mathematics ETH Z\"urich.
  Birkh\"auser Verlag, Basel, 2005.


\bibitem{Balder84}
{\sc E.~J. Balder}: {A general approach to lower semicontinuity and
lower
  closure in optimal control theory}. SIAM J. Control Optim. {\bf 22}
   4:570--598 (1984).





\bibitem{Ball89}
{\sc J. M. Ball}: {A version of the fundamental theorem for {Y}oung
measures}.
  {\it PDEs and continuum models of phase transitions (Nice 1988)}, Lecture Notes in
  Phys., vol. 344. pp.~207--215, Springer, Berlin, 1989.

\bibitem{Ball97}
\bysame: {Continuity properties and global attractors of
generalized semiflows and the Navier-Stokes equations}. {\it J. Nonlinear
Sci.} {\bf 7} 475--502 (1997).


\bibitem{Ball04}
   \bysame:
  {Global attractors for damped semilinear wave equations}.
   {\it Discrete Contin. Dyn. Syst.}
 {\bf 10} 1:31--52 (2004).

\bibitem{Blowey-Elliott91}
{\sc J. F~Blowey and C. M.~Elliott}: { The {C}ahn-{H}illiard gradient theory
 for the phase separations with nonsmooth free energy. {I}}.
{\it European J. Appl. Math.} {\bf 2}
 3:233--280 (1991).


\bibitem{Brezis71}
{\sc H.~Br\'ezis}: {Monotonicity methods in {H}ilbert spaces and some
  applications to nonlinear partial differential equations}. {\it Contribution to
  {N}onlinear {F}unctional {A}nalysis, Proc. {S}ympos. {M}ath. {R}es. {C}enter,
  {U}niv. {W}isconsin, {M}adison, 1971}. Academic Press, New York, 1971,

\bibitem{Brezis73}
\bysame: {\it Op\'erateurs maximaux monotones et semi-groupes de
  contractions dans les espaces de {H}ilbert}. North-Holland Publishing Co.,
  Amsterdam, 1973,

\bibitem{Brezis83}
\bysame: {\it Analyse fonctionnelle - {T}h\'eorie et applications},
Masson,
  Paris, 1983.



\bibitem{Brokate-Sprekels96} {\sc M. Brokate and J. Sprekels}:
 {\it Hysteresis and phase transitions}. Appl. Math.Sci., 121,
  Springer, New York, 1996.

\bibitem{Caginalp86}{\sc G. Caginalp}: {An analysis of a phase field model of a free boundary}. {\it Arch. Rational Mech. Anal.} {\bf 92}
   205--245 (1986).



\bibitem{Caraballo03}
 {\sc T. Caraballo, P. Marin-Rubio, and J. C. Robinson}:
{ A comparison between two theories for multi-valued semiflows
and their asymptotic behaviour}. {\it Set-Valued Anal.} {\bf 11}
 3:297--322 (2003).

\bibitem{Cardinali-Colombo-Papalini-Tosques97}
{\sc T.~Cardinali, G.~Colombo, F.~Papalini, and M.~Tosques}:
  {On a class of evolution equations without convexity}. {\it Nonlinear Anal.}
  {\bf 28} 2:217--234 (1997).



\bibitem{Chep-Vish95}
{\sc V. V. Chepyshoz and M. I. Vishik}: {\it Attractors for equations of
mathematical physics}. American Mathematical Society Colloquium
Publications, 49, American Mathematical Society, Providence, RI,
2002.



\bibitem{Crandall-Pazy69}
{\sc M.~G. Crandall and A.~Pazy}:
{Semi-groups of nonlinear
contractions and
  dissipative sets}. {\it J. Functional Analysis} {\bf 3} 376--418 (1969).

\bibitem{DeGiorgi93}
{\sc E.~De~Giorgi}: {New problems on minimizing movements}. {\it Boundary
Value
  Problems for PDE and Applications} (Claudio Baiocchi and Jacques~Louis Lions,
  eds.), Masson, Paris, 1993.

\bibitem{DeGiorgi-Marino-Tosques80}
{\sc E.~De~Giorgi, A.~Marino, and M.~Tosques}: {Problems of evolution
  in metric spaces and maximal decreasing curve}. {\it Atti Accad. Naz. Lincei Rend.
  Cl. Sci. Fis. Mat. Natur. (8)} {\bf 68}  3:180--187 (1980).

\bibitem{elmounir-simondon00}
{\sc A.~Ould Elmounir and F.~Simondon}:
 {Attracteurs compacts pour des probl\`{e}mes d'evolution sans unicit\'e}.
{\it  Ann. Fac. Sci. Toulouse Math. (6)} {\bf 9}  4:631--654 (2000).

\bibitem{evans-gariepy}
{\sc L. C.~Evans and R.~Gariepy}: {\it Measure theory and fine properties of functions},
Studies in Advanced Mathematics, CRC Press, Boca Raton, FL, 1992.


\bibitem{hale-raugel}
{\sc J. K.~Hale and G.~Raugel}: {Upper semicontinuity of the attractor for a singularly perturbed hyperbolic equation}. {\it J. Differential Equations} {\bf 73} 197-214 (1988).

\bibitem{Kapp-Mell-Vall03}
{\sc  A.~V. Kapustyan,  V. S. Melnik, and J. Valero}:
{Attractors of multivalued dynamical processes generated by phase-field equations}.  {\it Internat. J. Bifur. Chaos Appl. Sci. Engrg.} {\bf 13} 7:1969--1983 (2003).


\bibitem{Komura67}
{\sc Y.~K{\=o}mura}: {Nonlinear semi-groups in {H}ilbert space}. {\it J.
Math. Soc.
  Japan} {\bf 19} 493--507 (1967).



\bibitem{Marino-Saccon-Tosques89}
{\sc A.~Marino, C.~Saccon, and M.~Tosques}: {Curves of maximal slope
and
  parabolic variational inequalities on nonconvex constraints}. {\it Ann. Scuola
  Norm. Sup. Pisa Cl. Sci. (4)} {\bf 16}  2:281--330 (1989).

\bibitem{Mell-Vall98}
  {\sc  V. S. Melnik and J. Valero}:
{On attractors of multivalued semi-flows and differential inclusions}. {\it Set-Valued Anal. (4)} {\bf 6}
 83--111 (1998).

\bibitem{Mell-Vall00}
 {\sc V. S. Melnik and J. Valero}:
{On global attractors of multivalued semiprocesses and
nonautonomous evolution inclusions}. {\it Set-Valued Anal. (4)} {\bf 8}
 375--403 (2000).





\bibitem{Plotnikov-Starovoitov93}
{\sc P.~I. Plotnikov and V.~N. Starovoitov}: {The {S}tefan problem
with surface
  tension as the limit of a phase field model}. {\it Differential Equations} {\bf 29}  395--404 (1993).


\bibitem{Rocca-Schimperna04} {\sc E.~Rocca and G.~Schimperna}:
{Universal attractor for some singular phase transition
systems}. {\it Phys. D} {\bf 192} 3-4:279--307 (2004).


\bibitem{Rossi-Savare-Proc}
{\sc R.~Rossi and G.~Savar\'e}:{Existence and approximation results
for gradient flows}. {\it Rend. Mat. Acc. Lincei}  {\bf 15} 
183--196 (2004).


\bibitem{Rossi-Savare04}
{\sc R.~Rossi and G.~Savar\'e}: {Gradient flows of non convex
functionals in Hilbert spaces and applications}. {\it ESAIM Control Optim. Calc. Var.} {\bf 12} 3:564--614 (2006).

\bibitem{rossi-segatti-stefanelli-preprint}
{\sc R.~Rossi, A.~Segatti,  and U.~Stefanelli}: {Attractors for
gradient flows of non convex functionals and applications}. Preprint
IMATI-CNR n. 6-PV (2006), 1-47.


\bibitem{Schatzle00}
{\sc R.~Sch{\"a}tzle}: {The quasistationary phase field equations
with
  {N}eumann boundary conditions} {\it J. Differential Equations} {\bf 162}
  2:473--503 (2000).


\bibitem{Sell73}
{\sc G. R. Sell}: {Differential equations without uniqueness and
classical
              topological dynamics}. {\it J. Differential Equations}
{\bf 14} 42--56 (1973),.

\bibitem{Sell96}
\bysame: {Global attractors for the three-dimensional
{N}avier-{S}tokes
              equations}. {\it J. Dynam. Differential Equations}
       {\bf 8} 1:1--33 (1996).





\bibitem{Segatti04} {\sc A.~Segatti}: {Global attractor for a class
 of doubly nonlinear abstract evolution equations}. {\it Discrete Contin. Dyn. Syst.} {\bf 14} 4:801--820 (2006).

\bibitem{Segatti06} \bysame: {On the hyperbolic relaxation of the Cahn-Hilliard equation in 3-D: approximation and long time behaviour}. {\it Math. Models Methods Appl. Sci.} {\bf 17} 3:411-437 (2007).


\bibitem{siyk98}  {\sc K. Shirakawa, A. Ito,  N. Yamazaki, and N. Kenmochi}:
     {Asymptotic stability for evolution equations governed by
              subdifferentials}.
{\it Recent developments in domain decomposition methods and flow
              problems (Kyoto, 1996; Anacapri, 1996)}.
  GAKUTO Internat. Ser. Math. Sci. Appl.,
   {\bf 11},
  {Gakk\=otosho},
    {Tokyo},
    {1998}.

\bibitem{Simon86}
{\sc J. Simon}: {Compact sets in the space $L\sp p(0,T;B)$},  {\it Ann. Mat. Pura Appl. (4)} {\bf 146} 65--96, (1987).

\bibitem{Temam88} {\sc R. Temam}: {\it Infinite dimensional
mechanical systems in mechanics and physics}.
Applied Mathematical Sciences 68, Springer-Verlag, New York, 1988.




\bibitem{Valero01} {\sc J. Valero}: {
 Attractors for parabolic equations without uniqueness}. {\it J. Dynam. Differential Equations} {\bf 13} 711--744 (2001).

\bibitem{Visintin96}
{\sc A.~Visintin}: {\it Models of phase transitions}, Progress in
Nonlinear
  Differential Equations and Their Applications, vol.~28, Birkh\"auser, Boston,
  1996.


 \end{thebibliography}
\end{document}